\documentclass[11pt,a4paper]{article}
\usepackage[table]{xcolor}
\usepackage{amsmath,amssymb,amsthm,mathtools}
\usepackage{geometry}
\usepackage{booktabs}
\usepackage{array}
\usepackage{hyperref}
\usepackage{enumitem}
\usepackage{booktabs}
\usepackage{booktabs}
\usepackage{array}
\usepackage[table]{xcolor}
\usepackage{ragged2e}
\newtheorem{definition}{Definition}[section]
\newtheorem{theorem}{Theorem}[section]

\usepackage{cite}

\theoremstyle{definition}

\theoremstyle{Remark}

\theoremstyle{lemma}

\newtheorem{corollary}{Corollary}[section]
\theoremstyle{example}
\begin{document}
	\title{\textbf{Projection-Based and Modulus-Based Equivalent Reformulations of the Linear Complementarity Problem and Corresponding Iterative Methods: A Review}}
\author{Bharat Kumar$^{a}$, Abhishek Kumar Singh Sengar $^{b}$, Deepmala $^{c}$ \\
	\emph{\small $^{a}$Economic Sciences, Indian Institute of Technology, Kanpur 208016 (UP), India}\\
	\emph{\small $^{b,c}$Department of Mathematics, PDPM-IIITDM Jabalpur 482005 (MP), India}\\
	%		\emph{\small $^{b}$Economic Sciences, Indian Institute of Technology, Kanpur- 208016 (UP), India}\\
	%\emph{\small $^{b}$}\\
	\emph{\small   $^a$Email: bharatnishad.kanpu@gmail.com}\\
	\emph{\small   $^b$Email: abhishek47sengar@gmail.com}\\
		\emph{\small   $^c$Email: dmrai23@gmail.com}\\
}

\date{}
\maketitle

\begin{abstract}
	
	The Linear Complementarity Problem (LCP) is a fundamental
	mathematical model arising in optimization, equilibrium
	theory, and various scientific and engineering applications.
	Its complementarity structure admits several equivalent
	mathematical representations, which provide a foundation
	for developing a wide range of iterative solution methods.
	This review presents a systematic and critical analysis
	of the principal reformulations of the LCP, with particular
	emphasis on minimum-function, projection-based,
	positive-negative part, and modulus-based formulations.
	The mathematical equivalence of these representations
	to the standard LCP is examined, highlighting their
	underlying relationships and the preservation of the
	complementarity conditions.
	Building upon these formulations, the review surveys
	the corresponding iterative methods developed in the
	existing literature, focusing on their mathematical
	foundations, algorithmic structures, and convergence
	requirements. Particular attention is given to the
	influence of coefficient matrix properties, matrix-splitting
	strategies, and algorithmic parameters on the convergence
	behavior and computational characteristics of these
	methods. The reviewed approaches are further examined
	in terms of their theoretical convergence guarantees,
	computational cost, iteration requirements, numerical
	accuracy, and parameter sensitivity.
	A comparative perspective is developed to discuss the
	relative advantages and limitations of the different
	approaches, emphasizing that the suitability of an
	iterative method depends on the structural properties
	of the underlying LCP and the specific algorithm employed.
	In addition, relevant performance measures, including
	iteration count, CPU time, complementarity residual,
	memory requirements, and parameter sensitivity, are
	discussed as criteria for numerical evaluation. Finally,
	the review identifies important challenges and potential
	research directions, particularly in adaptive parameter
	selection, effective matrix-splitting strategies, scalable
	implementations, and reproducible numerical benchmarking,
	with the aim of supporting the development of reliable
	and computationally efficient methods for solving LCPs.
	
\end{abstract}

\section{Introduction}

The Linear Complementarity Problem  is a fundamental
mathematical problem that arises in a wide range of fields,
including mathematical programming, optimization, equilibrium
theory, economics, engineering, and numerical analysis.
It provides a general framework for expressing equilibrium
and optimality conditions and has therefore attracted
considerable attention in both theoretical investigations
and computational studies. A variety of important
optimization and equilibrium problems can be formulated
as, or transformed into, an LCP, making the study of its
mathematical structure, solvability, and numerical solution
methods particularly significant.

Let $A=(a_{ij})\in\mathbb{R}^{n\times n}$, where
$i,j=1,\ldots,n$, and let $q\in\mathbb{R}^{n}$ be a
given vector. The linear complementarity problem associated
with $(A,q)$, denoted by $\operatorname{LCP}(A,q)$, consists
of finding a vector $z\in\mathbb{R}^{n}$ such that
\begin{equation}
	\begin{aligned}
		w &= Az+q,\\
		z &\geq 0,\qquad w\geq 0,\qquad z^{T}w=0.
	\end{aligned}
	\label{eq:LCP}
\end{equation}
Equivalently, the complementarity condition can be expressed
componentwise as
\begin{equation*}
	z_iw_i=0,\qquad i=1,2,\ldots,n.
\end{equation*}
Thus, each pair of corresponding components must satisfy
nonnegativity and complementarity: at least one of $z_i$
and $w_i$ is zero, while both are nonnegative. Consequently,
the LCP combines a linear system with componentwise
inequality and complementarity constraints.

Unlike a conventional system of linear equations, the LCP
involves additional feasibility and complementarity
requirements that make its solution structure considerably
more challenging. Although the relation $w=Az+q$ is linear,
the complementarity constraints introduce a nonconvex
structure and can lead to nonsmooth mathematical
representations. These features complicate the analysis
of solution existence, uniqueness, and numerical
convergence. In particular, the properties of the
coefficient matrix $A$ play a central role in determining
the solvability of the problem and the convergence behavior
of various iterative algorithms. Consequently, different
classes of matrices and their associated structural
properties have been extensively investigated in the
complementarity literature.

An important approach to the theoretical analysis and
numerical solution of LCPs is to replace the original
complementarity conditions with mathematically equivalent
representations. Such reformulations express the
complementarity constraints through projection operators,
piecewise-defined functions, or absolute-value equations.
Although these representations describe the same
complementarity problem, they provide different mathematical
structures that can be exploited using tools from fixed-point
theory, matrix-splitting techniques, nonsmooth analysis,
and numerical optimization. In particular, they form
the basis for constructing a variety of iterative methods
for solving LCPs.

For the purpose of this review, the principal equivalent
reformulations considered in the literature are organized
into two broad frameworks: projection-based and
modulus-based formulations. The projection-based framework
includes projection, minimum-function, maximum-function,
and positive-negative part representations. These
formulations characterize complementarity through the
nonnegative orthant, its associated projection operator,
or equivalent componentwise constructions. In particular,
the positive-negative part representation expresses
complementarity through the decomposition of an
unrestricted vector into its positive and negative parts.
This provides an alternative mathematical representation
closely related to projection-based formulations and
facilitates the development of projection-type and
nonlinear-equation-based iterative schemes.

In contrast, modulus-based formulations encode the
complementarity conditions through absolute-value
transformations involving an auxiliary unrestricted
variable. By introducing appropriate scaling factors
and diagonal parameters, the original LCP can be
represented by an equivalent absolute-value equation.
This representation provides flexibility in constructing
iterative schemes and has motivated the development
of several classes of modulus-based algorithms.
Although projection-based and modulus-based formulations
use different mathematical representations, their
underlying equivalence to the original LCP establishes
a common foundation for studying their associated
numerical methods.

The availability of these alternative formulations has
led to the development of diverse iterative techniques.
For projection-based formulations, the existing literature
includes projected fixed-point methods, relaxation
methods, matrix-splitting and multisplitting schemes,
preconditioned approaches, and methods based on
positive-negative part representations. Similarly,
modulus-based formulations have motivated a broad
range of algorithms, including modulus-based relaxation,
matrix-splitting, accelerated, two-sweep, preconditioned,
inexact, and double-splitting methods. These developments
illustrate how different equivalent representations
can lead to distinct algorithmic structures and
computational strategies for solving the same class
of complementarity problems.

The convergence and computational performance of these
iterative methods depend on several interrelated factors,
including the structural properties of the coefficient
matrix, the choice of matrix-splitting strategy, and
the selection of algorithmic parameters. Accordingly,
the theoretical convergence conditions and computational
characteristics of the different approaches require
careful examination. Their relative merits can be
assessed in terms of convergence guarantees, iteration
requirements, computational cost, numerical accuracy,
and sensitivity to parameter choices. Such an assessment
also helps clarify the trade-offs between the cost
of an individual iteration and the total computational
effort required to obtain an approximate solution.

Motivated by these considerations, this review provides
a systematic and critical examination of the principal
equivalent reformulations of LCPs and the iterative
methods developed from them. It first establishes
the mathematical foundations and relationships of
the different reformulations and subsequently reviews
the associated algorithms, emphasizing their construction,
convergence requirements, and computational properties.
A comparative discussion is presented to identify
the advantages and limitations of the reviewed methods
and to highlight problem-dependent considerations
for selecting suitable iterative approaches. Relevant
performance measures, including iteration count,
CPU time, complementarity residual, memory requirements,
and parameter sensitivity, are also discussed as
criteria for numerical evaluation. Finally, important
open challenges and potential research directions
are identified, particularly in adaptive parameter
selection, effective matrix-splitting strategies,
scalable implementations, and reproducible numerical
benchmarking.

\subsection{Objectives of the Review}

The main objective of this review is to develop a systematic
and critical perspective on the principal equivalent
reformulations of the linear complementarity problem (LCP)
and the iterative solution techniques arising from these
formulations. Rather than treating the reformulations and
their associated algorithms as separate topics, this review
seeks to establish their mathematical connections and examine
how these connections influence the design, convergence,
and computational performance of iterative methods.
The specific objectives are as follows:

\begin{enumerate}
	\item To examine the principal mathematical reformulations
	of the LCP, including minimum-function, projection-based,
	positive-negative part, and modulus-based formulations,
	and to clarify their equivalence to the original
	complementarity conditions.
	
	\item To trace and critically review the development of
	iterative algorithms derived from these reformulations,
	emphasizing their underlying mathematical principles,
	computational procedures, and theoretical convergence
	requirements.
	
	\item To analyze how the structural properties of the
	coefficient matrix, the choice of matrix splitting, and
	the selection of algorithmic parameters affect the
	convergence conditions and computational behavior of
	the resulting iterative schemes.
	
	\item To assess the computational considerations relevant
	to the implementation and comparison of these methods,
	including iteration requirements, computational cost,
	numerical accuracy, robustness, and their applicability
	to different classes of LCPs.
\end{enumerate}

\subsection{Contributions of the Review}

The principal contribution of this review is to provide an
integrated mathematical and algorithmic perspective on the
major reformulation-based approaches for solving LCPs.
By systematically synthesizing the findings reported in
the existing literature, this review establishes connections
between different mathematical representations of the
complementarity conditions and the iterative methods
derived from them. It further examines the theoretical
properties and computational characteristics of these
methods to provide a basis for their critical comparison.
The main contributions are summarized as follows:

\begin{enumerate}
	\item A systematic classification of the principal
	equivalent reformulations of the LCP is presented,
	encompassing minimum-function, projection-based,
	positive-negative part, and modulus-based approaches.
	Their mathematical relationships and underlying
	equivalence properties are discussed to establish
	a coherent framework for understanding these
	different representations of the complementarity
	conditions.
	
	\item The major iterative methods associated with
	these reformulations are reviewed and critically
	analyzed in terms of their mathematical foundations,
	matrix-splitting strategies, convergence conditions,
	and algorithmic structures. Particular emphasis is
	placed on distinguishing theoretical convergence
	guarantees from computational performance reported
	in the existing literature.
	
	\item A comparative framework is established for
	examining the theoretical and computational
	characteristics of the reviewed methods. The
	comparison considers relevant performance measures,
	including iteration count, CPU time, complementarity
	residual, memory requirements, and sensitivity to
	algorithmic parameters. It also discusses the
	relationship between per-iteration computational
	cost and the overall effort required to obtain an
	approximate solution.
	
	\item The review identifies problem-dependent
	considerations for selecting appropriate iterative
	methods based on their theoretical and computational
	characteristics. In addition, it highlights important
	open problems and potential research directions,
	particularly in adaptive parameter selection,
	effective matrix-splitting strategies, scalable
	implementations, and standardized and reproducible
	numerical benchmarking.
\end{enumerate}

The remainder of this paper is organized as follows.
Section~2 introduces the notation, fundamental definitions,
and preliminary results required for the subsequent
discussion. Section~3 presents the principal equivalent
reformulations of the LCP within the projection-based
and modulus-based frameworks, together with a review
of the corresponding iterative methods developed in
the existing literature. Section~4 provides a critical
comparative assessment of these methods, emphasizing
their convergence properties, computational complexity,
dependence on matrix structure, parameter sensitivity,
and numerical robustness. It also examines their
respective advantages and limitations, discusses
performance measures and criteria for numerical
benchmarking and method selection, and identifies
important open challenges and potential directions
for future research. Finally, Section~5 summarizes
the main findings of the review and presents concluding
remarks on the development of efficient and robust
methods for solving LCPs.

\section{Preliminaries}
\label{sec:preliminaries}

This section introduces the notation, definitions, and fundamental
properties required for the subsequent discussion of equivalent
reformulations and iterative methods for the LCP. Throughout this paper, $\mathbb{R}^n$ denotes the
$n$-dimensional real Euclidean space, and $\mathbb{R}_+^n$ denotes
its nonnegative orthant. The principal notations and preliminary
concepts used in the subsequent sections are presented below.

\begin{enumerate}
	
	\item Let $A=(a_{ij})\in\mathbb{R}^{n\times n}$ and
	$B=(b_{ij})\in\mathbb{R}^{n\times n}$. We use
	$A\geq B$ (respectively, $A>B$) to denote
	$a_{ij}\geq b_{ij}$ (respectively, $a_{ij}>b_{ij}$)
	for all $1\leq i,j\leq n$.
	
	\item $(\star)^T$ denotes the transpose of the given matrix or vector.
	
	\item We use $A=O\in\mathbb{R}^{n\times n}$ to denote
	$a_{ij}=0$ for all $i,j$.
	
	\item For a vector or matrix, $|A|$ denotes the element-wise absolute
	value, i.e.,
	$
	|A|=(|a_{ij}|),~ \forall\, i,j.
	$
	Similarly, for a vector $x=(x_i)\in\mathbb{R}^n$,
$
	|x|=(|x_i|),~ \forall\,i.
	$
	
	\item $A^{-1}$ represents the inverse of the nonsingular matrix $A$.
	
	\item Let $a=(a_1,a_2,\ldots,a_n)\in\mathbb{R}^n$ and
	$b=(b_1,b_2,\ldots,b_n)\in\mathbb{R}^n$. The componentwise maximum
	of $a$ and $b$ is defined by
	\[
	\max\{a,b\} = \left(\max\{a_i,b_i\}\right)_{i=1}^{n}.
	\]
	
	\item Let $a=(a_1,a_2,\ldots,a_n)\in\mathbb{R}^n$ and
	$b=(b_1,b_2,\ldots,b_n)\in\mathbb{R}^n$. The componentwise minimum
	of $a$ and $b$ is defined by
\[
\min\{a,b\} = \left(\min\{a_i,b_i\}\right)_{i=1}^{n}.
\]
	
	\item For $x=(x_1,x_2,\ldots,x_n)\in\mathbb{R}^n$, the positive
	part of $x$ is defined componentwise by
	\[
	x_+ = \max\{0,x\} = \left(\max\{0,x_i\}\right)_{i=1}^{n}.
	\]
	\item For $x=(x_1,x_2,\ldots,x_n)\in\mathbb{R}^n$, the negative
	part of $x$ is defined componentwise by
	\[
	x_- = \max\{0,-x\} = \left(\max\{0,-x_i\}\right)_{i=1}^{n}.
	\]
	Consequently,
	$$
	x=x_+-x_-,
	~
	|x|=x_++x_-.
	$$
	\item $P_+(\cdot)$ denotes the Euclidean projection onto the nonnegative orthant, defined componentwise by
	\[
	[P_+(x)]_i=\max\{x_i,0\},\qquad i=1,\ldots,n.
	\]
	\item $\|\star\|_2$ denotes the Euclidean norm. In particular, for
	$x=(x_i)\in\mathbb{R}^n$,
	\[
	\|x\|_2
	=
	\sqrt{\sum_{i=1}^n x_i^2}.
	\]
	
\end{enumerate}

\begin{definition}\cite{BermanPlemmons1994}
	Let $A=(a_{ij})\in\mathbb{R}^{n\times n}$. The comparison matrix
	$\langle A\rangle=(\langle a_{ij}\rangle)$ of $A$ is defined by
	\[
	\langle a_{ij}\rangle
	=
	\begin{cases}
		|a_{ii}|, & i=j,\\
		-|a_{ij}|, & i\neq j,
	\end{cases}
	~ i,j=1,2,\ldots,n.
	\]
\end{definition}

\begin{definition}\cite{BermanPlemmons1994}
	Let $A\in\mathbb{R}^{n\times n}$. The matrix $A$ is called a
	$Z$-matrix if
$
	a_{ij}\leq 0,~ i\neq j.
	$
	A $Z$-matrix $A$ is called a nonsingular $M$-matrix if
	$A^{-1}\geq0$. The matrix $A$ is called an $H$-matrix if its
	comparison matrix $\langle A\rangle$ is a nonsingular $M$-matrix.
	Moreover, an $H$-matrix $A$ is called an $H_+$-matrix if
	$
	a_{ii}>0,~ i=1,2,\ldots,n.
	$
\end{definition}

\begin{definition}\cite{CottlePangStone2009}
	Let $A\in\mathbb{R}^{n\times n}$. The matrix $A$ is called a
	$P$-matrix if all its principal minors are positive, i.e.,
	$
	\det(A_{\alpha\alpha})>0
	$
	for every nonempty
	$\alpha\subseteq\{1,2,\ldots,n\}$.
\end{definition}

\begin{definition}\cite{HornJohnson2012}
	Let $A=(a_{ij})\in\mathbb{R}^{n\times n}$. The matrix $A$ is called
	strictly diagonally dominant (SDD) by rows if
$
	|a_{ii}|
	>
	\sum_{\substack{j=1\\j\neq i}}^n |a_{ij}|,
	~ i=1,2,\ldots,n.
$
\end{definition}

\begin{definition}\cite{HornJohnson2012}
	Let $A\in\mathbb{R}^{n\times n}$. The matrix $A$ is called symmetric
	positive definite (SPD) if
	$
	A=A^T
	$
	and
	$
	x^T A x>0,
	~
	\forall\,x\in\mathbb{R}^n,\quad x\neq0.
	$
\end{definition}

\begin{definition}\cite{BaiGolubNg2003}
	Let $A\in\mathbb{R}^{n\times n}$ and let
$
	A=M-N,
	$
	where $M$ is nonsingular. The splitting $A=M-N$ is called an
	$M$-splitting if $M$ is a nonsingular $M$-matrix and
	$
	N\geq0.
	$
	It is called an $H$-splitting if
	$
	\langle M\rangle-|N|
	$
	is a nonsingular $M$-matrix. Moreover, the splitting is called an
	$H$-compatible splitting if
	$
	\langle A\rangle
	=
	\langle M\rangle-|N|.
	$
\end{definition}

The matrix classes and matrix splittings introduced above play an
important role in the analysis of iterative methods for solving linear
complementarity problems. In the literature, different iterative
schemes derived from equivalent reformulations of the LCP are developed
under different structural assumptions on the coefficient matrix $A$.
These assumptions are not merely technical; rather, they are closely
related to the solvability of the LCP and to the convergence properties
of the corresponding iterative algorithms.

In particular, the assumption that $A$ is a $P$-matrix is fundamental
in the theory of LCPs. A $P$-matrix guarantees that the LCP has a unique
solution for every $q\in\mathbb{R}^n$, and hence this class of matrices
frequently appears in the analysis of iterative methods for LCPs. In
contrast, $H_+$-matrices constitute an important class characterized
through their comparison matrices. Since the comparison matrix of an
$H_+$-matrix is a nonsingular $M$-matrix, the properties of $M$-matrices
can be exploited to derive convergence conditions for various
iterative schemes. Consequently, $H_+$-matrix assumptions are commonly
employed in the convergence analysis of both projection-based and
modulus-based methods.

Strict diagonal dominance is another frequently used assumption in the
development and analysis of iterative methods. In particular, strictly
diagonally dominant matrices provide readily verifiable sufficient
conditions that can be used to establish bounds on the associated
iteration matrices and, consequently, to guarantee convergence under
suitable conditions. Similarly, symmetric positive definite (SPD)
matrices form an important class in the analysis of LCPs and related
iterative methods. The symmetry and positive definiteness of $A$ provide
favorable spectral properties that can be exploited in establishing
convergence and stability of the resulting iterative schemes. Thus,
several methods in the literature are developed under the assumption
that the coefficient matrix is symmetric positive definite.

Moreover, a number of iterative methods, particularly splitting-based
methods, are constructed by introducing a matrix splitting of the form
$
A=M-N,
$
where $M$ is nonsingular. The choice of the splitting has a direct
influence on the resulting iteration matrix and hence on the convergence
of the method. In this context, $H$-splittings and $H$-compatible
splittings are particularly useful because they establish a connection
between the splitting and the comparison matrix $\langle A\rangle$,
thereby allowing the theory of nonsingular $M$-matrices to be
incorporated into the convergence analysis. Such splittings provide
convenient sufficient conditions under which the spectral properties of
the associated iteration matrices can be controlled.

Therefore, the assumptions appearing in the existing iterative methods
are not identical. Depending on the equivalent reformulation and the
construction of the corresponding iteration, a method may require $A$
to be a $P$-matrix, an $H_+$-matrix, a strictly diagonally dominant
matrix, a symmetric positive definite matrix, or a matrix admitting a
particular type of splitting, such as an $H$-splitting or an
$H$-compatible splitting. These matrix properties provide different
sufficient frameworks for establishing the well-posedness, convergence,
and stability of the resulting iterative schemes. In this review, we
therefore report the matrix assumptions associated with each method in
order to provide a clear picture of the conditions under which the
method is theoretically justified and computationally applicable.

\section{Equivalent Reformulation of LCP }
The LCP admits several equivalent
reformulations that express its complementarity conditions in alternative
mathematical forms. These reformulations provide different perspectives
on the underlying problem and serve as a foundation for the development,
analysis, and implementation of various numerical methods. In this
section, we present the principal equivalent formulations of the LCP.

\subsection{Equivalent reformulation of LCP using minimum function}

	The complementarity conditions in the LCP can
	also be represented by a single nonsmooth equation using the
	componentwise minimum function. This reformulation provides a compact
	characterization of the nonnegativity and complementarity conditions and
	is useful for converting the LCP into a system of nonlinear equations.
	In particular, the LCP~\eqref{eq:LCP} is equivalently expressed as
	\begin{equation*}
		\min\{z,Az+q\}=0.
		\label{minequivalent}
	\end{equation*}
	This minimum-function formulation also serves as a convenient starting
	point for deriving other equivalent formulations, such as
	projection-based reformulations.
	
	The following theorem characterizes the solutions of the LCP~\eqref{eq:LCP}
	through an equivalent minimum-function equation.
\begin{theorem}\cite{Mangasarian1976}
	Let $A\in\mathbb{R}^{n\times n}$ and $q\in\mathbb{R}^n$. A vector $z\in\mathbb{R}^n$ solves the LCP~\eqref{eq:LCP} if and only if it solves the minimum-function equation
	\begin{equation}
		\min\{z,Az+q\}=0.
		\label{eq:min-equivalent}
	\end{equation}
\end{theorem}
\begin{proof}
	Let $w=Az+q$. Suppose first that $z$ solves the LCP~\eqref{eq:LCP}. Then
	$
	z\geq0,~ w\geq0,~ z^Tw=0.
	$
	Hence $z_iw_i=0$ for all $i$, and therefore $\min\{z_i,w_i\}=0$. Thus,
	$
	\min\{z,w\}=0,
	$
	which gives
	$
	\min\{z,Az+q\}=0.
	$
	
	Conversely, suppose that
	$
	\min\{z,Az+q\}=0.
	$
	Then $\min\{z_i,w_i\}=0$ for all $i$. This implies $z_i\geq0$ and $w_i\geq0$; otherwise, if either $z_i<0$ or $w_i<0$, then $\min\{z_i,w_i\}<0$. Since both are nonnegative and their minimum is zero, $z_iw_i=0$. Hence,
	$
	z\geq0,~ w\geq0,~ z^Tw=0.
	$
	
	Using $w=Az+q$, we obtain
	$
	z\geq0,~ Az+q\geq0,~ z^T(Az+q)=0,
	$
	so $z$ solves the LCP~\eqref{eq:LCP}.
\end{proof}

\begin{corollary}
	Let $A\in\mathbb{R}^{n\times n}$, $q\in\mathbb{R}^n$, and let $\alpha>0$. The minimum-function equation
	$
		\min\{z,Az+q\}=0
	$
	is equivalent to the projection equation
$
		z=P_+\left(z-\alpha(Az+q)\right),
	$
	where $P_+(x)=\max\{0,x\}$ denotes the projection onto the nonnegative orthant. Moreover, for $\alpha=1$, it reduces to the maximum-function equation
	\begin{equation}\label{maxpro}
		z=\max\{0,(I-A)z-q\}.
	\end{equation}
\end{corollary}
\begin{proof}
	Using $\min\{a,b\}=a-\max{0,a-b}$, we obtain
	
	$$
	\min\{z,Az+q\}=0
	\Longleftrightarrow
	z=\max\{0,z-(Az+q)\}
	=\max\{0,(I-A)z-q\},
	$$
	
	which gives \eqref{maxpro}. By the preceding theorem, the LCP is equivalent to
	
	$$
	z=P_+\left(z-\alpha(Az+q)\right),~ \alpha>0.
	$$
	
	Hence, the minimum-function, projection, and maximum-function formulations are equivalent, with the maximum-function formulation corresponding to $\alpha=1$.
\end{proof}

\subsection{Equivalent Reformulations of the LCP Using Projection Functions}

 Projection-based approaches provide an alternative characterization of
the complementarity conditions in the LCP. By exploiting the projection
onto the nonnegative orthant, the complementarity conditions can be
expressed as a fixed-point equation, which forms the basis for a variety
of projection and projection-type iterative methods.

\begin{theorem}\cite{Mangasarian1977}
	Let $A\in\mathbb{R}^{n\times n}$ and $q\in\mathbb{R}^n$. For every $\alpha>0$, a vector $z\in\mathbb{R}^n$ solves the LCP~\eqref{eq:LCP} if and only if it solves the fixed-point equation
	\begin{equation}\label{pro+}
		z=P_+\left(z-\alpha(Az+q)\right),
	\end{equation}
	where $P_+(x)=\max\{0,x\}$ denotes the projection onto the nonnegative orthant.
\end{theorem}
\begin{proof}
	Let $w=Az+q$. Suppose first that $z$ solves the LCP. Then, for each $i=1,\ldots,n$,
	$
	z_i\geq0,~ w_i\geq0,~ z_iw_i=0.
	$
	
	If $z_i>0$, then $w_i=0$, and hence $P_+(z_i-\alpha w_i)=z_i$. If $z_i=0$, then $w_i\geq0$, and thus $P_+(z_i-\alpha w_i)=P_+(-\alpha w_i)=0=z_i$. Therefore,
	$
	z=P_+(z-\alpha w)=P_+\left(z-\alpha(Az+q)\right).
	$

	Conversely, suppose that $z=P_+\left(z-\alpha(Az+q)\right)$. Then, componentwise,
	
	$$
		z_i=P_+(z_i-\alpha w_i).
	$$

	Hence $z_i\geq0$. If $z_i>0$, then $z_i-\alpha w_i=z_i$, so $w_i=0$. If $z_i=0$, then $P_+(-\alpha w_i)=0$, which implies $w_i\geq0$ since $\alpha>0$. Thus, $z_i\geq0$, $w_i\geq0$, and $z_iw_i=0$ for all $i=1,\ldots,n$. Consequently,
	$
	z\geq0,~ Az+q\geq0,~ z^T(Az+q)=0,
	$
	so $z$ solves $\operatorname{LCP}(A,q)$.
\end{proof}

The projection reformulation of the LCP can also be expressed through
the positive and negative parts of an auxiliary vector. In particular,
the complementarity structure of the LCP can be represented using these
positive and negative parts. This representation provides an alternative
way of incorporating the nonnegativity and complementarity conditions
within the projection-based framework. Specifically, the relations
\begin{equation*}
	x=x_+-x_-,
	~
	|x|=x_++x_-,
\end{equation*}
together with
\begin{equation*}
	x_+\geq0,~
	x_-\geq0,~
	(x_+)^{T}x_-=0,
\end{equation*}
allow the complementary variables to be expressed in terms of a single
vector. This observation leads to the following equivalent reformulation
of the LCP.

The following theorem provides an equivalent reformulation of the LCP
using the positive and negative parts of an auxiliary vector.

\begin{theorem}\cite{Fang2021}
	Let $A\in\mathbb{R}^{n\times n}$ and $q\in\mathbb{R}^n$. Define the positive and negative parts of $x\in\mathbb{R}^n$ by
	$
	x_+=\max\{0,x\},~ x_-=\max\{0,-x\},
	$
	then $z$ solves $\operatorname{LCP}(A,q)$ if and only if there exists $x\in\mathbb{R}^n$ such that
	$
		x_-=Ax_++q,
$
	with
	$
	z=x_+,~ w=x_-.
	$
	
\end{theorem}
\begin{proof}
	Suppose first that $z$ solves $\operatorname{LCP}(A,q)$ and let $w=Az+q$. Then
	$$
	z\geq0,~ w\geq0,~ z^Tw=0.
	$$
	Define $x=z-w$. Since $z$ and $w$ are nonnegative and complementary,
	$
	x_+=z,~ x_-=w.
	$
	Hence,
	$$
	x_-=w=Az+q=Ax_++q.
	$$
	Conversely, suppose that $x_-=Ax_++q$ and define
	$
	z=x_+,~ w=x_-.
$
	Then $z,w\geq0$ and, by the properties of the positive and negative parts,
	$
	z^Tw=(x_+)^Tx_-=0.
	$
	
	Moreover,
	$$
	w=x_-=Ax_++q=Az+q.
	$$
	
	Thus,
	$$
	z\geq0,~ Az+q\geq0,~ z^T(Az+q)=0,
	$$
	so $z$ solves $\operatorname{LCP}(A,q)$.
\end{proof}
\subsubsection{Existing Iterative Methods Based on the Projection Reformulation}
Projection-based methods constitute one of the classical approaches for
solving LCP. Their fundamental idea is
to reformulate the complementarity conditions as a fixed-point equation
involving the projection onto the nonnegative orthant. Since the projection
is applied componentwise, the resulting iteration naturally preserves
nonnegativity of the iterates. This basic principle has led to a wide
variety of projected relaxation, stationary, matrix-splitting,
multisplitting, preconditioned, and fixed-point methods for LCPs
\cite{Cryer1971,Mangasarian1977,Ahn1981,Pang1984,HadjidimosZhang2018}.

From Equations  (\ref{maxpro}),  and (\ref{pro+})  the basic projected fixed-point iteration is
\begin{equation*}
	z^{(k+1)}
	=
	P_+\left(
	z^{(k)}-\alpha(Az^{(k)}+q)
	\right)
	=
	\max\left\{
	0,(I-\alpha A)z^{(k)}-\alpha q
	\right\},
	~ k=0,1,2,\ldots .
	\label{eq:basic_PFI}
\end{equation*}
For $\alpha=1$, this reduces to
\begin{equation*}
	z^{(k+1)}
	=
	\max\left\{
	0,(I-A)z^{(k)}-q
	\right\}.
	\label{eq:basic_PFI_alpha1}
\end{equation*}
The principal developments of projection-based iterative methods are
summarized chronologically below.
\begin{enumerate}
\item
\textbf{Cryer (1971)}~\cite{Cryer1971} proposed a systematic overrelaxation method for solving the LCP in (\ref{eq:LCP}), where $A$ is a symmetric positive definite matrix. The LCP was formulated
equivalently as the strictly concave quadratic programming problem
\begin{equation*}
	z_i^{(k+1)}
	=
	\max\left\{
	0,\,
	z_i^{(k)}
	+
	\omega\frac{r_i^{(k+1)}}{a_{ii}}
	\right\},
	~ 0<\omega<2.
\end{equation*}
which has a unique solution. Cryer proposed a modified systematic
overrelaxation (SOR) method in which, for $i=1,\ldots,n$, the residual is
updated as
\begin{equation*}
	r_i^{(k+1)}
	=
	-q_i-\sum_{j<i}a_{ij}z_j^{(k+1)}
	-\sum_{j\geq i}a_{ij}z_j^{(k)},
\end{equation*}
followed by the projected relaxation step
\begin{equation*}
	z_i^{(k+1)}
	=
	\max\left\{
	0,\,
	z_i^{(k)}
	+\omega\frac{r_i^{(k+1)}}{a_{ii}}
	\right\},
	~ 0<\omega<2.
\end{equation*}
and
\begin{equation*}
	w_i^{(k+1)}
	=
	r_i^{(k+1)}+a_{ii}z_i^{(k+1)}.
\end{equation*}
Thus, the method combines sequential SOR updates with projection onto the
nonnegative orthant. Cryer proved convergence of the method for
$0<\omega<2$ and analyzed its asymptotic convergence rate in terms of the
SOR iteration matrix. The optimal relaxation parameter was related to the
SOR-optimal parameter of the principal submatrix corresponding to the
positive components of the solution, and procedures for estimating this
parameter were discussed.
	\item
\textbf{Mangasarian (1977)}~\cite{Mangasarian1977} proposed a unified class of projected iterative methods for solving the symmetric LCP in (\ref{eq:LCP}), where $A=A^T$. The basic idea is to compute a relaxed correction of the current iterate and then project it onto the nonnegative orthant. The general iteration is
\[
z^{k+1}=\lambda\left[z^k-\omega E^k
\left(Az^k+q+K^k(z^{k+1}-z^k)\right)\right]_+
+(1-\lambda)z^k,
\]
where $E^k$ is a positive diagonal matrix, $\omega>0$, 
$0<\lambda\leq1$, and $[\cdot]_+$ denotes projection onto the
nonnegative orthant.

\textit{Projected Jacobi Overrelaxation:}
For this method, $K^k=0$, so every component is updated using only
the previous iterate. The projection guarantees nonnegativity,
while $\omega$ provides relaxation of the iterative correction:
\[
z^{k+1}=\lambda[z^k-\omega E(Az^k+q)]_+
+(1-\lambda)z^k.
\]
Thus, the method extends the classical Jacobi iteration to the LCP
by incorporating projection and relaxation.

\textit{Projected SOR:}
Let $A=L+D+U$, where $L$ and $U=L^T$ are the strictly lower and
upper triangular parts of $A$. By choosing $K^k=L$ (or $U$), the
method uses newly computed components during the iteration:
\[
z^{k+1}=\lambda\left[z^k-\omega E
\left(Az^k+q+K^k(z^{k+1}-z^k)\right)\right]_+
+(1-\lambda)z^k.
\]
Hence, projected SOR combines the sequential updating principle of
SOR with projection onto the nonnegative orthant.

\textit{Projected Gauss-Seidel:}
Projected Gauss-Seidel is obtained as a special case of projected
SOR by taking $\lambda=1$, $\omega=1$, $E=D^{-1}$, and $K=L$.
Consequently, the latest available components are used immediately,
without an additional relaxation factor:
\[
z^{k+1}=\left[z^k-D^{-1}
\left(Az^k+q+L(z^{k+1}-z^k)\right)\right]_+.
\]
\textit{Projected Symmetric SOR:}
The symmetric version alternates between lower and upper triangular updates in successive iterations, using the lower triangular part (L) in one iteration and the upper triangular part (U) in the next, as follows:
\[
z^{k+1}=\lambda\left[z^k-\omega E
\left(Az^k+q+L(z^{k+1}-z^k)\right)\right]_+
+(1-\lambda)z^k.
\]
and \[
z^{k+1}=\lambda\left[z^k-\omega E
\left(Az^k+q+U(z^{k+1}-z^k)\right)\right]_+
+(1-\lambda)z^k.
\]
so that a forward projected SOR step is followed by a backward
projected SOR step. This alternating procedure combines both
directions of sequential updating while maintaining the
nonnegativity required by the LCP.
\item 
\textbf{Ahn (1981)}~\cite{Ahn1981} extended the iterative framework for LCPs from symmetric to nonsymmetric system matrices.
For the LCP in (\ref{eq:LCP}), with $A\in\mathbb{R}^{n\times n}$ and
$q\in\mathbb{R}^{n}$, Ahn considered a positive diagonal matrix $E$ and
a parameter $\omega>0$ and transformed the LCP into the fixed-point
formulation
$$
z=\left[z-\omega E(Az+q)\right]_+.
$$
 Based on this formulation, a general iterative scheme
was proposed using a relaxation parameter $\lambda\in(0,1]$ and bounded
matrix sequences $\{E^{(k)}\}$ and $\{K^{(k)}\}$:
\[
z^{(k+1)}
=
\lambda\left[
z^{(k)}-\omega E^{(k)}
\left(
Az^{(k)}+q+
K^{(k)}(z^{(k+1)}-z^{(k)})
\right)
\right]_+
+(1-\lambda)z^{(k)}.
\]
As a particular case, by taking $\lambda=1$, $E^{(k)}=E$, and
$K^{(k)}=K$, Ahn obtained the projected SOR-type iteration
\[
z^{(k+1)}
=
\left[
z^{(k)}-\omega E
\left(
Az^{(k)}+q+
K(z^{(k+1)}-z^{(k)})
\right)
\right]_+.
\]
Writing $A=D+L+U$, where $D$, $L$, and $U$ denote the diagonal, strictly
lower-triangular, and strictly upper-triangular parts of $A$, respectively,
the choices $K=L$, $K=U$, and $K=0$ lead to projected SOR-type and
projected Jacobi-type schemes; moreover, $E=D^{-1}$ and $\omega=1$
yield the projected Gauss-Seidel method. The convergence analysis was
established through a recursive inequality between successive iterates.
In particular, if
\[
\rho\!\left(
(I-\omega E|K|)^{-1}
\left|I-\omega E(A-K)\right|
\right)<1,
\]
then the generated sequence converges to a solution of the nonsymmetric
LCP. Ahn further established convergence for important classes of matrices,
including $M$-matrices and strictly quasi-dominant diagonal matrices.

	\item
\textbf{He (1992)}~\cite{He1992} proposed a projection and contraction (PC) method for solving the LCP (\ref{eq:LCP}) when the coefficient matrix $A$ is positive semidefinite, but not necessarily symmetric. The method is based on the equivalent projection formulation
\[
z=P_+[z-(Az+q)]=P_+(z-w),
\]
where $w=Az+q$ and $P_+(\cdot)$ denotes the componentwise projection onto the nonnegative orthant, defined by $[P_+(x)]_i=\max\{x_i,0\}$. He defined the merit function and search direction as
\[
\begin{aligned}
	\phi(z)&=\{z-P_+(z-w)\}^{T}w,\\
	g(z)&=A^{T}\{z-P_+(z-w)\}+w.
\end{aligned}
\]
For the construction of the iteration, He defined the index sets
\[
\begin{aligned}
	N(z)&=\{i\in\{1,\ldots,n\}:z_i=0,\ g_i(z)\geq0\},\\
	B(z)&=\{1,\ldots,n\}\setminus N(z).
\end{aligned}
\]
Here, $g_B(z)$ denotes the vector whose $i$th component is $g_i(z)$ for $i\in B(z)$ and zero for $i\in N(z)$; that is,
\[
[g_B(z)]_i=
\begin{cases}
	g_i(z),&i\in B(z),\\
	0,&i\in N(z).
\end{cases}
\]
The PC iteration is given by
\[
\begin{aligned}
	\alpha_k&=\frac{\phi(z^{(k)})}
	{\|g_B(z^{(k)})\|_2^2},\\
	z^{(k+1)}&=P_+\left[
	z^{(k)}-\alpha_k g_B(z^{(k)})\right].
\end{aligned}
\]
Regarding convergence, He established that, for every solution $z^*\in\mathcal{Z}^*$, the iterates satisfy
\[
\begin{aligned}
	\|z^{(k+1)}-z^*\|_2^2
	&\leq \|z^{(k)}-z^*\|_2^2-\frac{\phi^2(z^{(k)})}
	{\|g_B(z^{(k)})\|_2^2}.
\end{aligned}
\]
Consequently, the sequence is Fejér-monotone with respect to the solution set $\mathcal{Z}^*$. Moreover, if $\mathcal{Z}^*\neq\emptyset$, the method is globally convergent in the sense that
\[
\lim_{k\to\infty}
\operatorname{dist}(z^{(k)},\mathcal{Z}^*)=0.
\]
\item
\textbf{He (1996)}~\cite{He1996} proposed a modified projection and contraction (PC) method for solving a class of linear complementarity problems with positive semidefinite, not necessarily symmetric, coefficient matrices. For the LCP in (\ref{eq:LCP}), let $w=Az+q$. The method is based on the projection residual
\[
e(z)=z-P_+(z-w),
\]
where $P_+(\cdot)$ denotes the componentwise projection onto the nonnegative orthant, defined by $[P_+(x)]_i=\max\{x_i,0\}$. He defined the merit function and search direction as
\[
\phi(z)=e(z)^T w,\qquad
g(z)=A^T e(z)+w.
\]
The modified PC method determines the step-size parameter
\[
p_k=
\frac{\|e(z^{(k)})\|^2}
{\|(I_n+A^T)e(z^{(k)})\|^2},
\]
where $I_n$ is the $n\times n$ identity matrix. The iterates are updated according to
\[
z^{(k+1)}
=P_+\left[z^{(k)}-\gamma_k p_k g(z^{(k)})\right],
\]
where $0<\gamma_1<\gamma_2<2$ and $\gamma_k\in[\gamma_1,\gamma_2]$. Regarding convergence, He established that, for every solution $z^*\in\mathcal{Z}^*$, the generated sequence satisfies
\[
\begin{aligned}
	\|z^{(k+1)}-z^*\|^2
	&\leq \|z^{(k)}-z^*\|^2-\gamma_k(2-\gamma_k)p_k
	\|e(z^{(k)})\|^2,
\end{aligned}
\]
where $\mathcal{Z}^*$ denotes the solution set of the LCP. Consequently, the distance of the iterates to the solution set is nonincreasing. If $\mathcal{Z}^*\neq\emptyset$, the method is globally convergent, and the paper further establishes global linear convergence.

	\item
\textbf{Bai and Evans (1997)}~\cite{BaiEvans1997} developed synchronous and asynchronous matrix multisplitting relaxation methods for solving LCP (\ref{eq:LCP}). The coefficient matrix \(A\) is decomposed as
$$
A=D+L+U,
$$
with \(D=\operatorname{diag}(A)\), \(L\) strictly lower triangular, and \(U\) strictly upper triangular. More generally, Bai and Evans introduced a multisplitting of \(A\) of the form
$$
A=D+L_i+U_i=D+B_i,~ i=1,\ldots,\alpha,
$$
where \(D=\operatorname{diag}(A)\), \(L_i\) is strictly lower triangular, \(U_i\) is a zero-diagonal matrix, and \(E_i\) are nonnegative diagonal weighting matrices satisfying
$$
\sum_{i=1}^{\alpha}E_i=I.
$$
First, they proposed the synchronous matrix multisplitting AOR method (Method 3.1). Given an initial approximation \(z^{(0)}\in\mathbb{R}^n\), for \(p=0,1,2,\ldots\), compute
$$
z^{p+1}=\sum_{i=1}^{\alpha}E_i z^{p,i},
$$
where
$$
z^{p,i}
=
\left(
z^p-\gamma E
\left[
rL_i z^{p,i}
+
(\omega A-rL_i)z^p
+\omega q
\right]
\right)_+,
~ i=1,\ldots,\alpha.
$$
Here, \((\cdot)_+\) denotes the projection onto the nonnegative orthant, and \(r\) and \(\omega\) are relaxation parameters. For \(r=\omega\), this method reduces to the matrix multisplitting SOR method, whereas for \(\alpha=1\), it reduces to the projected SOR method.

A practical form of Method 3.1  is obtained by choosing
$
\gamma E=D^{-1}.
$
Thus,
$$
z^{p+1}=\sum_{i=1}^{\alpha}E_i z^{p,i},
$$
where
$$
z^{p,i}
=
\left(
z^p-D^{-1}
\left[
rL_i z^{p,i}
+
(\omega A-rL_i)z^p
+\omega q
\right]
\right)_+,
~ i=1,\ldots,\alpha.
$$
For \(r=\omega\), Method 3.2 gives the corresponding matrix multisplitting SOR method.

To obtain an asynchronous scheme, Bai and Evans introduced index sets \(J(p)\) and delay sequences \(s_i(p)\). The asynchronous matrix multisplitting AOR method (Method 4.1) is given by
$$
z^{p+1}=\sum_{i=1}^{\alpha}E_i z^{p,i},
$$
where
$$
z^{p,i}
=
\begin{cases}
	\left(
	z^{s_i(p)}
	-\gamma E
	\left[
	rL_i z^{p,i}
	+
	(\omega A-rL_i)z^{s_i(p)}
	+\omega q
	\right]
	\right)_+,
	& i\in J(p),\\[2mm]
	z^p,
	& i\notin J(p),
\end{cases}
~ i=1,\ldots,\alpha.
$$
Here, \(J(p)\) specifies the processors that perform an update at iteration \(p\), while \(s_i(p)\) specifies the possibly delayed iterate used by the \(i\)th processor. If
$$
J(p)=\{1,2,\ldots,\alpha\},
~
s_i(p)=p,
$$
then Method 4.1 reduces to Method 3.1. Moreover, for \(r=\omega\), it yields an asynchronous matrix multisplitting SOR method.

Finally, by choosing
$
\gamma E=D^{-1},
$
Bai and Evans obtained the practical form (Method 4.2):
$$
z^{p+1}=\sum_{i=1}^{\alpha}E_i z^{p,i},
$$
where
$$
z^{p,i}
=
\begin{cases}
	\left(
	z^{s_i(p)}
	-D^{-1}
	\left[
	rL_i z^{p,i}
	+
	(\omega A-rL_i)z^{s_i(p)}
	+\omega q
	\right]
	\right)_+,
	& i\in J(p),\\[2mm]
	z^p,
	& i\notin J(p),
\end{cases}
~ i=1,\ldots,\alpha.
$$
For \(r=\omega\), Method 4.2 corresponds to a particular form of the asynchronous matrix multisplitting SOR method. Under the assumptions established by Bai and Evans, these methods converge to the unique solution of the LCP.
	\item
\textbf{Bai (1998)}~\cite{Bai1998Monotone} further investigated the 
matrix multisplitting accelerated overrelaxation (MMAOR) methods introduced 
by Bai and Evans \cite{BaiEvans1997}. Rather than introducing a fundamentally new iterative scheme, the author established the monotone convergence properties of the existing MMAOR method for solving LCP (\ref{eq:LCP}). The coefficient matrix is decomposed according to the multisplitting
\[
A=D+L_i+U_i=D+B_i,\qquad i=1,\ldots,\alpha,
\]
where the weighting matrices satisfy
$
\sum_{i=1}^{\alpha}E_i=I.
$
The MMAOR iteration considered in the paper is
\begin{equation*}
	z^{p+1}=\sum_{i=1}^{\alpha}E_i z^{p,i}.
\end{equation*}
and 
\begin{equation*}
	z^{p,i}
	=
	\left(
	z^p-\gamma E
	\left[
	rL_i z^{p,i}
	+
	(\omega A-rL_i)z^p
	+\omega q
	\right]
	\right)_+,
\end{equation*}
where $r \in [0, \infty)$ is a relaxation factor, $\omega \in (0, \infty)$ is an acceleration factor $\gamma \in (0, \infty)$ is an arbitrary and $E$ an arbitrary non singular matrix.

Bai established monotone convergence results for this method and derived
comparison results describing the effects of the multiple splittings,
relaxation parameters, and auxiliary parameters on the convergence rate.
Thus, the main contribution of Bai was the development of the
monotone convergence theory and comparison analysis for the previously
proposed MMAOR framework, rather than the introduction of a new iterative
method.

\item \textbf{Koulisianis and Papatheodorou (2003)}
~\cite{KoulisianisPapatheodorou2003} proposed an improved projected
successive over-relaxation (IPSOR) method for LCP arising from moving-boundary problems and American option
valuation. The LCP  can be formulated in an infinite-dimensional setting as finding a function (u) satisfying
$$
u\geq 0,\quad Lu-b\geq 0,\quad u^T(Lu-b)=0,
$$
where $(L)$ is a linear operator and $(b)$ is a given function. After discretization, the problem reduces to the discrete LCP
$$
U\geq 0,\quad AU-\beta\geq 0,\quad U^T(AU-\beta)=0,
$$
where (A) is the coefficient matrix and $(\beta)$ is the right-hand-side vector. An extension is the implicit LCP, which incorporates a constraint vector (G) and is written as
$$
U\geq G,\quad AU-\beta\geq 0,\quad
(U-G)^T(AU-\beta)=0,
$$
Writing
$
A=D+L+U,
$
where $D=\operatorname{diag}(A)$, $L$ is strictly lower triangular, and
$U_A$ is strictly upper triangular, the projected successive over-relaxation
(PSOR) method computes
\begin{equation*}
	y^{k+1}_{(i)}
	=
	(1-\omega)X^k_{(i)}
	+
	\frac{\omega}{a_{ii}}
	\left[
	\beta_{(i)}
	-\sum_{m=1}^{i-1}a_{im}X^{k+1}_{(m)}
	-\sum_{m=i+1}^{n}a_{im}X^k_{(m)}
	\right],
\end{equation*}
followed by
\begin{equation*}
	X^{k+1}_{(i)}
	=
	\max\left\{y^{k+1}_{(i)},G_{(i)}\right\},
	~ \omega\in(0,2).
\end{equation*}
To reduce the computational cost, the authors exploited the moving-boundary
structure and proposed Algorithm~5.1. For an index set $s$ containing
$n-k$ consecutive equations $k,k+1,\ldots,n-1$, the first $k-1$
components are set equal to their corresponding constraint values, and the
reduced problem
$
	\operatorname{DLCP}(\beta^s,A^s,G^s)
$
is solved by PSOR. The condition
$
	A_{k-1}U\geq\beta_{k-1}
$
is then checked. If this condition is not satisfied, the dimension of the
reduced problem is increased by one and the reduced problem is solved again.

The authors further proposed Algorithm~5.2, in which the dimension of the
reduced problem is increased by a prescribed integer $\nu$ instead of one.
Thus, when
$
	A_{k-1}U<\beta_{k-1},
$
the equations
$
	k-\nu,\ldots,k-1
$
are added to the reduced subsystem and the resulting DLCP is solved again
by PSOR. Algorithm~5.2 constitutes the proposed IPSOR method. Numerical
experiments showed that IPSOR reduces the computational time while
maintaining the accuracy of the standard PSOR method.
	\item
\textbf{Li and Dai (2007)}~\cite{LiDai2007} proposed a generalized accelerated over-relaxation (GAOR) method for solving LCP (\ref{eq:LCP}), based on the matrix splitting
\[
A=D+L+U,
\]
where \(D=\operatorname{diag}(A)\), and \(L\) and \(U\) are the strictly lower and strictly upper triangular parts of \(A\), respectively. The LCP is first
reformulated as the projected fixed-point equation
\begin{equation*}
	z=\left(z-D^{-1}(Az+q)\right)_+,
\end{equation*}
where $(x)_+=\max\{0,x\}$ componentwise. Based on this reformulation, the
authors introduced the GAOR iteration
\begin{equation*}
	z^{(k+1)}
	=
	\left(
	z^{(k)}
	-
	D^{-1}
	\left[
	\alpha\Omega Lz^{(k+1)}
	+
	\left(\Omega A-\alpha\Omega L\right)z^{(k)}
	+
	\Omega q
	\right]
	\right)_+,
\end{equation*}
where
$
\Omega=\operatorname{diag}(\omega_1,\ldots,\omega_n)
$
is a diagonal relaxation matrix and $\alpha$ is an acceleration parameter.
For $\alpha=1$, the GAOR method reduces to the generalized SOR (GSOR)
method. The authors established sufficient convergence conditions for the
proposed methods under suitable matrix assumptions, including $H$-matrices,
$M$-matrices, and strictly or irreducibly diagonally dominant matrices.
Numerical experiments were also provided to illustrate the effectiveness
of the proposed GAOR and GSOR methods.
	\item
\textbf{Dehghan and Hajarian (2009)}~\cite{DehghanHajarian2009} studied successive over-relaxation (SSOR) iterative methods for solving LCP (\ref{eq:LCP}). The coefficient matrix \(A\) is decomposed as
$
A=D-L-U,
$
where \(D\), \(L\), and \(U\) denote the diagonal, strictly lower triangular, and strictly upper triangular parts of \(A\), respectively. By applying
the SSOR splitting, the authors proposed two iterative methods for the
LCP. For $0<\omega<2$, the first method is given by
\begin{equation*}
	z^{(k+1)}
	=
	\left(
	z^{(k)}
	-
	D^{-1}
	\left[
	-\omega Lz^{(k+1)}
	+
	\{\omega(2-\omega)A+\omega L\}z^{(k)}
	-
	\omega(2-\omega)q
	\right]
	\right)_{+},
\end{equation*}
where $(\cdot)_{+}$ denotes the componentwise projection onto the
nonnegative orthant. The second method is obtained by interchanging
the roles of $L$ and $U$, and is given by
\begin{equation*}
	z^{(k+1)}
	=
	\left(
	z^{(k)}
	-
	D^{-1}
	\left[
	-\omega Uz^{(k+1)}
	+
	\{\omega(2-\omega)A+\omega U\}z^{(k)}
	-
	\omega(2-\omega)q
	\right]
	\right)_{+}.
\end{equation*}
The convergence of these two SSOR methods was established when the
coefficient matrix $A$ is an $H$-matrix with positive diagonal
entries; in particular, the results also apply to $M$-matrices.
The authors further presented two numerical examples to illustrate
the efficiency of the proposed methods.
	\item
\textbf{Li, Jiang, Cheng, and Liao (2011)}~\cite{LiJiangChengLiao2011} proposed an interval generalized accelerated over-relaxation (IGAOR) method, an interval extension of the GAOR method, for solving LCP (\ref{eq:LCP}). The proposed method is based
on a combination of the GAOR method and interval arithmetic. In
addition, the authors developed a class of multisplitting IGAOR
methods that can be implemented in parallel.

For the underlying GAOR framework, let
$
A=D+L+U,
$
where $D$, $L$, and $U$ denote the diagonal, strictly lower triangular,
and strictly upper triangular parts of $A$, respectively. The
corresponding GAOR projected iteration is
\begin{equation*}
	z^{(k+1)}
	=
	\left[
	z^{(k)}
	-
	D^{-1}
	\left\{
	\alpha\Omega Lz^{(k+1)}
	+
	\left(\Omega A-\alpha\Omega L\right)z^{(k)}
	+
	\Omega q
	\right\}
	\right]_{+},
\end{equation*}
where $\alpha$ is a real parameter and
$\Omega=\operatorname{diag}(\omega_1,\omega_2,\ldots,\omega_n)$ is a
real diagonal relaxation matrix. The IGAOR method is obtained by
combining this GAOR framework with interval arithmetic.

Furthermore, the authors developed multisplitting IGAOR methods by
extending the IGAOR framework to a multisplitting setting. These
methods are designed so that the computations associated with the
different splittings can be carried out in parallel. For $H$-matrices
with positive diagonal elements, the authors established convergence
of the proposed IGAOR and multisplitting IGAOR methods. Numerical
experiments were also presented to illustrate the efficiency of the
proposed algorithms.
	\item
\textbf{Saberi Najafi and Edalatpanah (2013)}~\cite{SaberiNajafi2013} proposed preconditioned generalized accelerated over-relaxation (GAOR) methods for solving LCP (\ref{eq:LCP}), based on the matrix splitting
$
A=D+L+U,
$
where \(D\), \(L\), and \(U\) denote the diagonal, strictly lower triangular, and strictly upper triangular parts of \(A\), respectively. By choosing
$\alpha E=D^{-1}$ in the fixed-point formulation, they obtained
\begin{equation*}
	z=\left(z-D^{-1}(Az+q)\right)_+.
\end{equation*}
Based on this formulation, the GAOR iteration is given by
\begin{equation*}
	z^{(k+1)}
	=
	\left(
	z^{(k)}
	-
	D^{-1}
	\left[
	\alpha\Omega Lz^{(k+1)}
	+
	\left(\Omega A-\alpha\Omega L\right)z^{(k)}
	+
	\Omega q
	\right]
	\right)_+,
\end{equation*}
where $\alpha$ is a real parameter and
$
\Omega=\operatorname{diag}(w_1,\ldots,w_n)
$
is a real diagonal relaxation matrix. For an $H$-matrix $A$ with positive
diagonal elements, the authors recalled the convergence result for the GAOR
method. In particular, if
\[
G=I-\alpha\Omega D^{-1}|L|,
~
F=\left|I-D^{-1}\left(\Omega A-\alpha\Omega L\right)\right|,
\]
then
\begin{equation*}
	\rho(G^{-1}F)
	\leq
	\max_i\left\{
	|1-w_i|+w_i\rho(|J|)
	\right\}
	<1,
\end{equation*}
provided that
\[
0<w_i<\frac{2}{1+\rho(|J|)},
~
0\leq\alpha\leq1,
\]
where
$
J=D^{-1}(L+U)
$
is the Jacobi iteration matrix.

To accelerate the convergence of the GAOR method, the authors then
introduced $(I+K_i)$-type preconditioners and considered the preconditioned
system
\begin{equation*}
	\overline{A}=(I+K_i)A,
	~
	\overline{q}=(I+K_i)q,
\end{equation*}
with
$
	\overline{A}
	=
	\overline{D}
	+
	\overline{L}
	+
	\overline{U}.
$
The corresponding preconditioned GAOR iteration is
\begin{equation*}
	z^{(k+1)}
	=
	\left(
	z^{(k)}
	-
	\overline{D}^{-1}
	\left[
	\alpha\Omega\overline{L}z^{(k+1)}
	+
	\left(
	\Omega\overline{A}
	-
	\alpha\Omega\overline{L}
	\right)z^{(k)}
	+
	\Omega\overline{q}
	\right]
	\right)_+.
\end{equation*}
For an $H$-matrix $A$ with positive diagonal elements, the authors proved
that the preconditioned GAOR methods converge and established the comparison
\begin{equation*}
	\rho\left(\overline{G}^{-1}\overline{F}\right)
	\leq
	\rho\left(G^{-1}F\right)
	<1,
\end{equation*}
where $\overline{G}$ and $\overline{F}$ are the corresponding matrices
associated with the preconditioned GAOR method. Thus, the preconditioned
GAOR methods have a convergence rate no worse than that of the corresponding
GAOR method. Numerical experiments were also presented to demonstrate the
effectiveness of the proposed preconditioners.
	\item
\textbf{Hadjidimos and Tzoumas (2016)}~\cite{HadjidimosTzoumas2016} introduced the projected matrix analogue of the accelerated over-relaxation (MAAOR) method for solving LCP (\ref{eq:LCP}) with an irreducible \(H_+\)-matrix. This method extends the GAOR approach by employing diagonal parameter matrices. The coefficient matrix is decomposed as
$
A=D-L-U,
$
where \(D\) is the positive diagonal part of \(A\). By scaling the problem with \(D^{-1}\), the authors obtained the normalized matrix
$
\widetilde{A}=D^{-1}A=I-\widetilde{L}-\widetilde{U},
$
where \(\widetilde{L}=D^{-1}L\), \(\widetilde{U}=D^{-1}U\), and \(\widetilde{q}=D^{-1}q\). Using a positive diagonal matrix $\Omega$
and a diagonal parameter matrix $R$, they proposed the projected
MAAOR iteration
\[
\begin{aligned}
	z^{(k+1)}=\Big[ z^{(k)}-\big(&-R\widetilde{L}z^{(k+1)}
	+\Omega(I-\widetilde{L}-\widetilde{U})z^{(k)}\\
	&+R\widetilde{L}z^{(k)}+\Omega\widetilde{q}\big)\Big]_+,
\end{aligned}
\]
where $[\cdot]_+$ denotes projection onto the nonnegative orthant.
For the convergence analysis, they established an error bound
involving the nonnegative matrix
\[
G=
(I-|R|\,|\widetilde{L}|)^{-1}
\left(
|I-\Omega|+
|\Omega-R|\,|\widetilde{L}|+
|\Omega|\,|\widetilde{U}|
\right).
\]
They showed that $\rho(G)<1$ is a sufficient condition
for convergence. Using Perron--Frobenius theory and regular-splitting
arguments, the authors derived sufficient convergence intervals
for the diagonal parameters, considering both strictly diagonally
dominant $H_+$-matrices and general irreducible $H_+$-matrices.

For the optimal choice $R=\Omega=I$, the projected MAAOR method
reduces to the projected Gauss--Seidel method, given by
\[
z^{(k+1)}
=\left[\widetilde{L}z^{(k+1)}
+\widetilde{U}z^{(k)}-\widetilde{q}\right]_+.
\]
To compare this method with the best modulus-based splitting
Gauss--Seidel (MBSGS) method, the authors considered their
associated nonnegative majorizer matrices. For the MBSGS method,
the majorizer is
$
T_1=(2I-|\widetilde{L}|)^{-1}
\bigl(|\widetilde{L}|+2|\widetilde{U}|\bigr),
$
whereas for the projected Gauss--Seidel method, it is
$
T_2=(I-|\widetilde{L}|)^{-1}|\widetilde{U}|.
$
Here, $\widetilde{L}$ and $\widetilde{U}$ are the strictly lower
and upper triangular parts of the normalized coefficient matrix,
respectively. The authors established that
\[
\rho(T_2)\leq\rho(T_1)<1,
\]
with strict inequality between the spectral radii for irreducible
matrices. Thus, the projected Gauss--Seidel method is superior
to the MBSGS method in terms of the spectral radii of their
associated majorizer matrices, and $\rho(T_2)<1$ provides a
sufficient convergence condition for the projected method.
		\item
\textbf{Shi, Yang, and Huang (2016)}~\cite{ShiYangHuang2016} proposed an explicit fixed-point iterative method for solving LCP (\ref{eq:LCP}) arising from American option pricing. For a vector \(z\in\mathbb{R}^n\), its positive and negative parts are defined componentwise as
\[
z=x_+,\qquad
w=x_-,
\]
such that
$
x=x_+-x_-.
$
Introducing a positive relaxation parameter \(\alpha\), the authors obtained the fixed-point equation
\[
x=x_+-\alpha(Ax_++q).
\]
Accordingly, the proposed explicit fixed-point iteration is
\begin{equation*}
	x^{(k+1)}
	=
	x_+^{(k)}
	-
	\alpha\left(Ax_+^{(k)}+q\right),
\end{equation*}
where
\begin{equation*}
	x_+^{(k)}=P_{+}\left(x^{(k)}\right).
\end{equation*}
For the LCP arising from American option
pricing, the coefficient matrix is considered in the splitting
\begin{equation*}
	A=D-U-L,
\end{equation*}
where $D$, $U$, and $L$ denote the diagonal, strictly upper triangular,
and strictly lower triangular parts of $A$, respectively. Using the
Gauss-Seidel idea, the authors considered the componentwise iteration
\begin{equation*}
	x_i^{(k+1)}
	=
	{x_i}_+^{(k)}
	-
	\alpha
	\left[
	(D-U)x_+^{(k)}
	-
	Lx_+^{(k+1)}
	+
	q
	\right]_i,
	~ i=1,\ldots,n,
\end{equation*}
where the components of $x^{(k+1)}$ are computed successively.
Equivalently, the most recently computed positive components are used
in the lower triangular part of the iteration.

For the American option pricing problem, the coefficient matrix is a
symmetric tridiagonal matrix of the form
\begin{equation*}
	a_{ii}=|1-\alpha2\lambda\theta|,
\end{equation*}
\begin{equation*}
	a_{i,i+1}=a_{i+1,i}=\alpha\lambda\theta,
	~ i=1,\ldots,n-1,
\end{equation*}
with all other entries equal to zero. The authors established that the
sequence generated by the above fixed-point iteration converges to the
unique solution of the fixed-point equation provided that
\begin{equation*}
	0<\alpha<
	\frac{2}{1+4\lambda\theta}.
\end{equation*}
Numerical experiments for American option pricing were presented, and
the proposed fixed-point method was reported to be more accurate than
the PSOR method for the test problems considered.

\item 
\textbf{Edalatpanah (2020)}~\cite{Edalatpanah2020} investigated preconditioned projective iterative methods for solving LCP (\ref{eq:LCP}) and showed that several previously proposed preconditioned projected iterative (PPI) models do not, in general, preserve the complementarity condition of the original problem. In particular, for a left preconditioner \(P\), the preconditioned matrix and vector are given by
\[
\overline{A}=PA,\qquad \overline{q}=Pq,
\]
but the resulting preconditioned problem is not necessarily equivalent to the original LCP.

The paper first recalled the projected fixed-point formulation
\begin{equation*}
	z=\left(z-\alpha E(Az+q)\right)_+,
\end{equation*}
where $\alpha>0$ and $E$ is a diagonal matrix with positive diagonal
elements. For the splitting
$
	A=D-L-U,
$
where $D$, $L$, and $U$ denote the diagonal, strictly lower triangular,
and strictly upper triangular parts of $A$, respectively, the
Generalized Accelerated Overrelaxation (GAOR) iteration is given by
\begin{equation*}
	z^{(k+1)}
	=
	\left(
	z^{(k)}
	-
	D^{-1}
	\left[
	-\alpha\Omega Lz^{(k+1)}
	+
	\left(\Omega A-\alpha\Omega L\right)z^{(k)}
	+
	\Omega q
	\right]
	\right)_+,
\end{equation*}
where $\alpha$ is a real parameter and
$\Omega=\operatorname{diag}(w_1,\ldots,w_n)$ is a real diagonal
relaxation matrix.

The authors then examined existing preconditioned projective
iterative models and demonstrated, by a counterexample, that the
solution obtained from such a preconditioned LCP need not satisfy the
complementarity condition of the original LCP. To overcome this
difficulty, they proposed a new preconditioned GAOR strategy based on
a symmetric preconditioner. Given a symmetric preconditioner $P$, the
transformed quantities are defined by
\begin{equation*}
	\widetilde{z}^{(0)}=P^{-1}z^{(0)},
	~
	\widetilde{q}=Pq,
	~
	\widetilde{A}=PAP,
\end{equation*}
where
$
	\widetilde{A}
	=
	\widetilde{D}
	-
	\widetilde{L}
	-
	\widetilde{U}.
$
The resulting preconditioned GAOR iteration is
\begin{equation*}
	\widetilde{z}^{(k+1)}
	=
	\left(
	\widetilde{z}^{(k)}
	-
	\widetilde{D}^{-1}
	\left[
	-\alpha\Omega\widetilde{L}\widetilde{z}^{(k+1)}
	+
	\left(
	\Omega\widetilde{A}
	-
	\alpha\Omega\widetilde{L}
	\right)
	\widetilde{z}^{(k)}
	+
	\Omega\widetilde{q}
	\right]
	\right)_+,
\end{equation*}
and, after convergence,
\begin{equation*}
	z=P\widetilde{z}.
\end{equation*}

The authors further proposed a new preconditioner $P$ of monomial
type. Let $\gamma,d\in\mathbb{R}\setminus\{0\}$ and, for each
$r=1,\ldots,n$, define
\begin{equation*}
	k_r=\min\{j:a_{rj}=a_{jr}\},
	~ j=r+1,\ldots,n,
\end{equation*}
with $k_r=n+1$ when the corresponding set is empty. They then set
$k_v=\min\{k_r\}$ and define
\begin{equation*}
	p_{v k_v}=p_{k_v v}=-\gamma a_{v k_v},
\end{equation*}
while
\begin{equation*}
	p_{ii}=d,
	~
	i\in\{1,\ldots,n\}\setminus\{v,k_v\},
\end{equation*}
with all remaining entries equal to zero. Thus, the proposed
preconditioner is a symmetric monomial matrix.

For an $M$-matrix $A$, the authors proved that the preconditioned
matrix
$
	\widetilde{A}=PAP
$
is also an $M$-matrix when
\begin{equation*}
	\gamma>0,~ d>0.
\end{equation*}
Consequently, for any initial vector $z^{(0)}\in\mathbb{R}^n$, the
iterative sequence generated by the proposed preconditioned GAOR
algorithm converges to the unique solution of the original
$\operatorname{LCP}(A,q)$. Numerical experiments were also presented.
For the test problem considered, the preconditioned GAOR method
required fewer iterations and less CPU time than the corresponding
unpreconditioned GAOR method.
	\item
\textbf{Fang (2021)}~\cite{Fang2021} introduced a general fixed-point (GFP) method for solving LCP (\ref{eq:LCP}). By introducing a positive diagonal parameter matrix \(\Omega\), the problem was reformulated as the fixed-point equation
\[
z=z_+-\Omega(Az_++q),
\]
where \(\Omega\) is a positive diagonal matrix. Using \(z=z_+-z_-\), this equation can equivalently be expressed as
\[
z_-=\Omega(Az_++q).
\]
The equivalence between the fixed-point equation and the original LCP
was proved in the paper. In particular, if $A$ is a $P$-matrix, then
for any positive diagonal matrix $\Omega$, the fixed-point equation has
a unique solution.

Based on this formulation, Fang proposed the general fixed-point
(GFP) method in the basic form
\begin{equation*}
	z^{(k+1)}
	=
	z_+^{(k)}
	-
	\Omega\left(Az_+^{(k)}+q\right)
	=
	(I-\Omega A)z_+^{(k)}-\Omega q,
	\quad k=0,1,2,\ldots .
\end{equation*}
For a $P$-matrix $A$, this iteration converges to the unique solution
of the fixed-point equation for any initial vector
$z^{(0)}\in\mathbb{R}^{n}$ if
\begin{equation*}
	\rho(|I-\Omega A|)<1.
\end{equation*}

The paper also considered a second, converted form of the GFP method
based on the idea of Gauss-Seidel. Let
\[
A=D-L-U,
\]
where $D$, $-L$, and $-U$ denote the diagonal, strictly lower
triangular, and strictly upper triangular parts of $A$, respectively.
Then the second GFP iteration is given by
\begin{equation*}
	z^{(k+1)}
	=
	(I-\Omega(D-U))z_+^{(k)}
	+
	\Omega Lz_+^{(k+1)}
	-
	\Omega q,
	\quad k=0,1,2,\ldots .
\end{equation*}
For a $P$-matrix $A$, the GFP iterations were shown to converge under
suitable spectral-radius conditions. In particular, for $H_{+}$-matrices,
the choice $\Omega=\omega D^{-1}$ guarantees convergence for
\[
0<\omega<
\frac{2}{1+\rho(D^{-1}|B|)},
\quad B=L+U,
\]
with $\omega=1$ being optimal for the corresponding bound. For symmetric
positive definite matrices, the first GFP iteration converges when
\[
0<\omega<\frac{2}{\lambda_{\max}(A)}.
\]
Numerical experiments were presented to illustrate the convergence
properties and effectiveness of the two GFP iterations. The numerical
results also showed that the converted GFP iteration based on the
Gauss-Seidel idea generally required fewer iterations than the basic
GFP iteration when the same parameter value was used. The paper further
compared the GFP method with the fixed-point method of Shi, Yang, and
Huang, showing the effectiveness of the proposed generalization.

\item \textbf{Kumar, Deepmala, and Das (2022, 2024, 2025)}~\cite{KumarDeepmalaDas2022,KumarDeepmalaDas2023,DasDeepmalaKumar2025} proposed a class of fixed-point and projected-type matrix-splitting methods for solving large and sparse LCP (\ref{eq:LCP}). They developed generalized fixed-point formulations using positive diagonal parameter matrices and matrix splittings. In the  work \cite{KumarDeepmalaDas2023}, two 
positive diagonal matrices $\Omega_1$ and $\Omega_2$ were introduced with
\[
z=\Omega_1x_+,~ w=\Omega_2x_-,
\]
leading to the fixed-point iteration
\[
x^{(k+1)}
=
\left(I-\Omega_2^{-1}A\Omega_1\right)(x^{(k)})_+
-\Omega_2^{-1}q,
~
z^{(k+1)}=\Omega_1(x^{(k+1)})_+.
\]
A modified general fixed-point method (MGFP) \cite{KumarDeepmalaDas2022} was further obtained by 
introducing a strictly lower triangular matrix $\Omega$ and the splitting
\[
A=(D+\Omega)-(L+U+\Omega),
\]
which gives
\[
x^{(k+1)}
=
\left[I-\Omega_2^{-1}(D+\Omega-U)\Omega_1\right](x^{(k)})_+
+\Omega_2^{-1}(L+\Omega)\Omega_1(x^{(k+1)})_+
-\Omega_2^{-1}q.
\]
For $P$-matrices, convergence was established under suitable 
spectral-radius conditions, while sufficient convergence domains were 
also derived for $H_+$-matrices. In their subsequent work, the authors 
extended this framework to a class of new projected-type matrix-splitting 
methods \cite{DasDeepmalaKumar2025} based on
\[
A=(M+I+D)-(N+I+D),
\]
and proposed the general projected-type iteration
\[
x_{+}^{(k+1)}
=
(M+2I+D)^{-1}
\left[
(N+I+D)x_{+}^{(k)}
+\left|(A-I)x_{+}^{(k)}+q\right|-q
\right].
\]
Different choices of the splitting matrices yield the new 
projected-type Gauss-Seidel (NPGS), successive overrelaxation (NPSOR), 
accelerated overrelaxation (NPAOR), and Jacobi (NPJ) methods, with NPAOR 
reducing to NPSOR, NPGS, and NPJ for
\[
(\alpha,\beta)=(\alpha,\alpha),~ (1,1),~ (1,0),
\]
respectively. For $P$-matrices, convergence was guaranteed under a 
spectral-radius condition, while for $H_+$-matrices sufficient 
convergence conditions were obtained in terms of $M$-matrix properties. 
Numerical experiments showed that the proposed methods generally require 
fewer iterations than the corresponding existing GFP, MGS, and MSOR 
methods, demonstrating their effectiveness for solving large and sparse 
LCPs.
\end{enumerate}
\paragraph{}

The literature on projection-based approaches to LCPs has evolved through several stages, with successive studies extending the basic projection idea to stationary, relaxation, splitting, multisplitting, preconditioned, and fixed-point schemes. This development has been accompanied by efforts to establish broader convergence conditions and to improve the computational efficiency of the resulting algorithms. More recent contributions have further generalized the framework through positive-negative part representations and projected matrix-splitting techniques. A chronological summary of the major developments is given in Table~\ref{tab:lcp_chronology}.
\begin{table*}[ht]
	\centering
	\caption{Chronological development of iterative methods for LCPs.}
	\label{tab:lcp_chronology}
	\small
	\begin{tabular}{ccc}
		\hline
		\textbf{Year} & \textbf{Authors} & \textbf{Main Contribution} \\
		\hline
		
		1971 & \textbf{Cryer} &
		Projected relaxation methods for LCPs. \\
		
		1977 & \textbf{Mangasarian} &
		Projected stationary iterative methods. \\
		
		1981 & \textbf{Ahn} &
		Methods for nonsymmetric LCPs . \\
		
			1992 & \textbf{He} &
		Projection and contraction method for LCPs. \\
			
			1996 & \textbf{He} &
			A modified PC method  for a class of LCPs. \\
		
		1997-1998 & \textbf{Bai and Evans} &
		Synchronous and asynchronous   methods. \\
		
		1998 & \textbf{ Bai} &
	MMAOR  methods. \\
		
		2003 & \textbf{Koulisianis and Papatheodorou} &
	Improved 	projected SOR  methods. \\
		
		2007 & \textbf{Li and Dai} &
		Generalized  AOR methods. \\
		
		2009 & \textbf{Dehghan and Hajarian} &
		Projected SSOR methods. \\
		
		2011 & \textbf{Li et al.} &
		IGAOR and multisplitting IGAOR methods. \\
		
		2013 & \textbf{Saberi et al.} &
		Preconditioned iterative methods. \\
		
		2016 & \textbf{Hadjidimos and Tzoumas} &
		PMAAOR and PMAGS methods. \\
		
		2016 & \textbf{Shi et al.} &
		Explicit fixed point methods. \\

		2020 & \textbf{Edalatpanah} &
	Preconditioned	projected fixed-point methods. \\
		
		2021 & \textbf{Fang} &
		Generalized fixed point methods. \\
		
2022,2024,2025 & \textbf{Kumar et al.} &
		Projected matrix-splitting methods. \\
		
		\hline
	\end{tabular}
	
\end{table*}
Methods based on the scaled modulus transformation
\[
z=\gamma^{-1}(|x|+x),
~
w=\gamma^{-1}\Omega(|x|-x),
\]
and the corresponding modulus-based matrix-splitting literature are
treated separately in the following subsection.
\subsection{ Equivalent formulations of LCP using modulus-based functions}
\begin{enumerate}
	\item 
Modulus-based reformulations provide an alternative representation of the
LCP by introducing an auxiliary variable and
expressing the complementary variables in terms of its positive and
negative parts. This transformation converts the complementarity
conditions into an implicit equation involving the componentwise absolute
value. Such a formulation provides a useful framework for the development
of modulus-based iterative methods.

To establish the modulus-based reformulation, let
$\Omega=\operatorname{diag}(\omega_1,\ldots,\omega_n)$ be a positive
diagonal matrix and let $\gamma>0$. The following theorem gives the
equivalence between the LCP and its modulus-based formulation.
\begin{theorem}\cite{Bai2010}
	Let $A\in\mathbb{R}^{n\times n}$, $q\in\mathbb{R}^n$, let $\Omega=\operatorname{diag}(\omega_1,\ldots,\omega_n)$ be a positive diagonal matrix, and let $\gamma>0$. A pair $(z,w)$ solves the LCP (\ref{eq:LCP})
	if and only if there exists $x\in\mathbb{R}^n$ satisfying
	\begin{equation}\label{modulseq}
	(\Omega+A)x=(\Omega-A)|x|-\gamma q,
		\end{equation}
	with
	$
		z=\frac{1}{\gamma}(|x|+x),~
		w=\frac{1}{\gamma}\Omega(|x|-x).
	$
\end{theorem}
\begin{proof}
	Suppose first that $(z,w)$ solves the LCP. Define
	$
	x=\frac{\gamma}{2}\left(z-\Omega^{-1}w\right).
	$
	Since $z_i,w_i\geq0$ and $z_iw_i=0$, we have
	$
	|x|+x=\gamma z,~ |x|-x=\gamma\Omega^{-1}w.
	$
	Thus,
	$$
	z=\frac{1}{\gamma}(|x|+x),~
	w=\frac{1}{\gamma}\Omega(|x|-x).
	$$
	Using $w=Az+q$ gives
	$
	\Omega(|x|-x)=A(|x|+x)+\gamma q,
	$
	and hence
	$$
	(\Omega+A)x=(\Omega-A)|x|-\gamma q.
	$$
	Conversely, suppose that
	$
	(\Omega+A)x=(\Omega-A)|x|-\gamma q.
	$
	Equivalently,
	$$
	\Omega(|x|-x)=A(|x|+x)+\gamma q.
	$$
	Define
	$
	z=\frac{1}{\gamma}(|x|+x),~
	w=\frac{1}{\gamma}\Omega(|x|-x).
	$
	Since $\gamma>0$ and $\Omega$ is positive diagonal, $z,w\geq0$. Moreover,
	$$
	z_iw_i=\frac{\omega_i}{\gamma^2}(|x_i|+x_i)(|x_i|-x_i)
	=\frac{\omega_i}{\gamma^2}(|x_i|^2-x_i^2)=0.
	$$
	Thus, $z^Tw=0$, and the above equation yields $w=Az+q$. Therefore, $(z,w)$ solves the LCP.
\end{proof}

\item An alternative reformulation of the LCP can be obtained by exploiting
the elementary identity
\begin{equation*}\label{abequivalent}
	|a-b|=a+b,
\end{equation*}
which holds for nonnegative complementary scalars $a$ and $b$.
This identity provides a convenient way to express the complementarity
condition through an absolute-value equation. We first state the
corresponding scalar characterization.
\begin{theorem}
	Let $a,b\in\mathbb{R}^n$ satisfy $a\geq0$ and $b\geq0$. Then the componentwise complementarity condition
	$
	a_i b_i=0,~ i=1,\ldots,n,
	$
	is equivalent to
	$
		|a-b|=a+b.
	$
\end{theorem}
\begin{proof}
	Suppose first that $a_i b_i=0$ for all $i$. Since $a_i,b_i\geq0$, at least one of $a_i$ and $b_i$ is zero, and hence
	$$
	|a_i-b_i|=a_i+b_i,~ i=1,\ldots,n.
	$$
	Thus, $|a-b|=a+b$. Conversely, suppose that $|a-b|=a+b$. Then
	$$
	|a_i-b_i|=a_i+b_i,~ i=1,\ldots,n.
	$$
	Squaring both sides gives
	$$
	(a_i-b_i)^2=(a_i+b_i)^2,
	$$
	and hence $4a_i b_i=0$. Therefore, $a_i b_i=0$ for $i=1,\ldots,n$, which establishes the componentwise complementarity condition.
\end{proof}
In the following result, we derive an equivalent formulation of the
LCP by exploiting the identity
$
|a-b|=a+b,
$
which characterizes complementarity for nonnegative vectors.
\begin{theorem}\cite{MangasarianMeyer2006}
	Let $A\in\mathbb{R}^{n\times n}$ and $q\in\mathbb{R}^n$. A vector $z\in\mathbb{R}^n$ solves the LCP~\eqref{eq:LCP} if and only if it solves the absolute-value equation
	\begin{equation}\label{modequalityeqvalent}
		\left|(I-A)z-q\right|=(I+A)z+q.
	\end{equation}
\end{theorem}
\begin{proof}
	Suppose first that $z$ solves $\operatorname{LCP}(A,q)$ and let $w=Az+q$. Then
	$$
	z\geq0,~ w\geq0,~ z_iw_i=0,~ i=1,\ldots,n.
	$$
	Hence, $|z-w|=z+w$. Since
	$
	z-w=(I-A)z-q,~ z+w=(I+A)z+q,
	$
	we obtain
	$$
	\left|(I-A)z-q\right|=(I+A)z+q.
	$$
	Conversely, suppose that $\left|(I-A)z-q\right|=(I+A)z+q$ and let $w=Az+q$. Then $|z-w|=z+w$. Thus, for each $i$,
	$
	|z_i-w_i|=z_i+w_i.
	$
	Since the left-hand side is nonnegative, $z_i+w_i\geq0$. Squaring gives
	$$
	(z_i-w_i)^2=(z_i+w_i)^2,
	$$
	and hence $z_iw_i=0$. Since at least one of $z_i,w_i$ is zero and their sum is nonnegative, it follows that $z_i,w_i\geq0$. Therefore, $z\geq0$, $w\geq0$, and $z^Tw=0$. Together with $w=Az+q$, these are precisely the LCP conditions.
\end{proof}
\end{enumerate}
\subsection{Existing Iterative Methods Based on the Modulus Reformulation}

Modulus-based methods provide an alternative framework for solving linear 
complementarity problems (LCPs), in which the complementarity conditions are 
incorporated through an absolute-value transformation. An auxiliary 
unrestricted variable is introduced, and the complementary variables are 
represented through the positive and negative parts of this variable. This 
transformation converts the LCP into an equivalent absolute-value equation 
or a related fixed-point formulation and provides the foundation for a broad 
class of modulus-based iterative methods 
\cite{Bai2010,Li2013,WuLi2016}.

From equation~\eqref{modulseq}, the basic modulus-based iteration is given by
\begin{equation*}
	(\Omega+A)x^{k+1}
	=(\Omega-A)|x^k|-\gamma q.
	\label{eq:mb-basic}
\end{equation*}

More generally, consider the matrix splitting
$
A=M-N,
$
where
\[
M=\frac{1}{\alpha}(D-\beta L)
\]
and
\[
N=\frac{1}{\alpha}
\left[(1-\alpha)D+(\alpha-\beta)L+\alpha U\right],
\]
with $\alpha\neq 0$ and $\beta$ being suitable splitting parameters.
Using this splitting, we define the following splitting-based iterative scheme:
\begin{equation*}
	(\Omega+M)x^{k+1}
	=Nx^k+(\Omega-A)|x^k|-\gamma q.
	\label{eq:mb-splitting}
\end{equation*}
Consequently, the corresponding iterative formula is given by
\begin{equation*}
	x^{k+1}
	=(\Omega+M)^{-1}
	\left[Nx^k+(\Omega-A)|x^k|-\gamma q\right].
	\label{eq:mb-iteration}
\end{equation*}

Similarly, for the modulus-based iteration given in equation~\eqref{modequalityeqvalent},
\begin{equation*}
	(I+A)z^{(k+1)}
	=\left|(I-A)z^{(k)}-q\right|-q,
\end{equation*}
using the same matrix splitting $A=M-N$, we define the following splitting-based iterative scheme:
\begin{equation*}
	(I+M)z^{(k+1)}
	=Nz^{(k)}
	+\left|(I-A)z^{(k)}-q\right|-q.
\end{equation*}
Hence, the corresponding iterative formula can be expressed as
\begin{equation*}
	z^{(k+1)}
	=(I+M)^{-1}
	\left[
	Nz^{(k)}
	+\left|(I-A)z^{(k)}-q\right|-q
	\right].
\end{equation*}

This equation provides the basic implicit fixed-point formulation from 
which the matrix-splitting-based modulus iterations are developed. The 
principal developments based on this reformulation are summarized below.

\begin{enumerate}

\item \textbf{Van Bokhoven (1980)}~\cite{vanBokhoven1980} introduced the \emph{modulus algorithm} for solving LCP (\ref{eq:LCP}). The method transforms the complementarity problem into an implicit equation involving the modulus function through the change of variables
\[
w=|x|-x,\qquad z=|x|+x.
\]
Using this transformation, the problem can be expressed as
\[
(A+I)x+(A-I)|x|+q=0.
\]
When $I+A$ is nonsingular, this equation can be written in the
fixed-point form
\begin{equation*}
	x=(I+A)^{-1}(I-A)|x|-(I+A)^{-1}q.
\end{equation*}

Based on the above fixed-point equation, van Bokhoven defined the
\emph{modulus algorithm} as
\begin{equation*}
	x^{(k+1)}
	=
	(I+A)^{-1}(I-A)|x^{(k)}|
	-
	(I+A)^{-1}q,
	~ k=0,1,2,\ldots .
\end{equation*}
If
\begin{equation*}
	\left\|
	(I+A)^{-1}(I-A)
	\right\|<1,
\end{equation*}
then the modulus algorithm converges to the unique solution of the
LCP. The corresponding solution is recovered from the limit
$x^{*}$ through
\[
z^{*}=|x^{*}|+x^{*},
~
w^{*}=|x^{*}|-x^{*}.
\]
Van Bokhoven further considered the case in which $A$ is symmetric
positive definite. Defining
\begin{equation*}
	\delta(A)
	=
	\left\|
	(I+A)^{-1}(I-A)
	\right\|_{2},
\end{equation*}
he established polynomial-time solvability of the LCP by the modulus
algorithm when $A$ is symmetric positive definite and has a bounded
condition number. Thus, the work introduced the modulus algorithm as a
fixed-point approach for LCPs and provided the basis for subsequent
developments of modulus-type algorithms.
 \item \textbf{Kappel and Watson (1986)}~\cite{KappelWatson1986} extended the modulus algorithm to a class of nonsymmetric \(P\)-matrices and developed a block version for solving LCP (\ref{eq:LCP}). By introducing the transformation
 \[
 w=|x|-x,\qquad z=|x|+x,
 \]
 where \(x\in\mathbb{R}^{n}\), they obtained the equivalent fixed-point equation
 \[
 x=-(I+A)^{-1}q+(I+A)^{-1}(I-A)|x|,
 \]
 provided that \(I+A\) is nonsingular. Defining
\begin{equation*}
	T=(I+A)^{-1}(I-A),
	~
	b=-(I+A)^{-1}q,
\end{equation*}
the transformed problem takes the form
\begin{equation*}
	x=T|x|+b.
\end{equation*}

The corresponding fixed-point iteration is
\begin{equation*}
	x^{(k+1)}
	=
	T|x^{(k)}|+b,
	~
	k=0,1,2,\ldots .
\end{equation*}
The authors established convergence of this iteration when $T$ is a
contraction in an appropriate absolute norm. More specifically, their
result applies when, for every principal pivot submatrix $\widetilde A$
of $A$, $I+\widetilde A$ is invertible and there exists a positive
definite diagonal matrix $\Gamma$ and some $p$, $1\leq p\leq\infty$,
such that
\begin{equation*}
	\left\|
	\Gamma
	(I+\widetilde A)^{-1}
	(I-\widetilde A)
	\Gamma^{-1}
	\right\|_{p}
	<1.
\end{equation*}
Under these conditions, $A$ is a $P$-matrix and the LCP is solvable
by the modulus algorithm.

To reduce the computational work of the basic modulus algorithm, the
authors developed a block modulus algorithm that can determine more
than one component in the same outer cycle and consequently reduce the
dimension of the LCP by more than one. Two selection strategies were
also considered: the modulus algorithm \emph{(max)}, which selects the
component with the largest magnitude, and the modulus algorithm
\emph{(test)}, which selects the first component satisfying the
prescribed test value and otherwise chooses the component with the
largest magnitude. The block algorithm applies the same test criterion
to all components satisfying the condition during the same outer cycle.

The authors also implemented the direct fixed-point iteration
\begin{equation*}
	x^{(k+1)}
	=
	T|x^{(k)}|+b,
\end{equation*}
and compared the different algorithms numerically. Numerical
experiments were carried out on several test problems, including
nonsymmetric problems, and showed that the block modulus algorithm can
substantially reduce the number of outer iterations when more than one
component can be determined during an outer cycle.
	\item
\textbf{Hadjidimos and Tzoumas (2009)}~\cite{HadjidimosTzoumas2009} considered LCP (\ref{eq:LCP}) with a symmetric positive-definite coefficient matrix \(A\) and introduced optimal nonstationary extrapolation for the modulus algorithm (MA) and the block modulus algorithm (BMA). Using the modulus transformation
\[
z=|x|+x,\qquad w=|x|-x,
\]
they obtained the fixed-point equation
\[
x=T|x|+b,
\qquad
T=(I+A)^{-1}(I-A),
\qquad
b=-(I+A)^{-1}q.
\]
The corresponding modulus iteration is given by
\[
x^{(k+1)}=T|x^{(k)}|+b,
\qquad k=0,1,2,\ldots.
\]
To accelerate the convergence, the authors introduced a positive
extrapolation parameter $\omega$ and obtained
\begin{equation*}
	x^{(k+1)}
	=
	T_{\omega}|x^{(k)}|+b_{\omega},
\end{equation*}
where
\begin{equation*}
	T_{\omega}
	=
	(I+\omega A)^{-1}(I-\omega A),
	~
	b_{\omega}
	=
	-(I+\omega A)^{-1}\omega q.
\end{equation*}
For every $\omega>0$, the extrapolated iteration converges. If
$\lambda_{\min}$ and $\lambda_{\max}$ are the smallest and largest
eigenvalues of $A$, respectively, the optimal parameter is
\begin{equation*}
	\omega^{*}
	=
	\frac{1}{\sqrt{\lambda_{\min}\lambda_{\max}}},
\end{equation*}
with
\begin{equation*}
	\rho(T_{\omega^{*}})
	=
	\frac{\sqrt{\lambda_{\max}}-\sqrt{\lambda_{\min}}}
	{\sqrt{\lambda_{\max}}+\sqrt{\lambda_{\min}}}.
\end{equation*}
By recalculating the optimal extrapolation parameter at the beginning
of each cycle as the LCP is reduced, the authors obtained the
\emph{Nonstationary Extrapolated Modulus Algorithm (NSEMA)}. The same
nonstationary extrapolation was applied to the BMA, leading to the
\emph{Nonstationary Extrapolated Block Modulus Algorithm (NSEBMA)}.
Numerical experiments showed that the optimal NSEBMA can significantly
improve the convergence performance over the related modulus-based
algorithms.
\item
\textbf{Dong and Jiang (2009)}~\cite{DongJiang2009} proposed a modified modulus method for solving LCP (\ref{eq:LCP}) with a symmetric positive-definite coefficient matrix \(A\). By introducing a positive parameter \(\alpha\) into the classical modulus formulation, they defined the transformation
\[
z=|x|+x,\qquad
w=\alpha\bigl(|x|-x\bigr),
\]
which leads to the equivalent fixed-point equation
\[
(\alpha I+A)x=(\alpha I-A)|x|-q,
\qquad \alpha>0.
\]
Accordingly, the modified modulus iteration is given by
\[
(\alpha I+A)x^{(k+1)}
=(\alpha I-A)|x^{(k)}|-q,
\qquad k=0,1,2,\ldots.
\]
Once the iteration has converged, the solution of the original LCP is
recovered from
\begin{equation*}
	z^{*}=|x^{*}|+x^{*},
	~
	w^{*}=\alpha\bigl(|x^{*}|-x^{*}\bigr).
\end{equation*}
For a symmetric positive-definite matrix $A$, the authors established
convergence of the modified modulus method for every $\alpha>0$ and
analyzed the optimal choice of the parameter. In particular, the
optimal parameter is
\begin{equation*}
	\alpha^{*}
	=
	\sqrt{\lambda_{\min}\lambda_{\max}},
\end{equation*}
where $\lambda_{\min}$ and $\lambda_{\max}$ are the smallest and
largest eigenvalues of $A$, respectively.

To reduce the computational cost of solving the linear system at each
outer iteration, the authors further introduced an
\emph{inexact modified modulus method}, in which the system
\begin{equation*}
	(\alpha I+A)x^{(k+1)}
	=
	(\alpha I-A)|x^{(k)}|-q
\end{equation*}
is solved only approximately by an iterative inner solver. The authors
established global convergence of the resulting inexact method under
suitable conditions on the inner iterations and discussed two specific
implementations. Numerical experiments demonstrated that the modified
and inexact modified modulus methods can be more efficient than the
classical modulus method under suitable conditions.
\item \textbf{Bai (2010)}~\cite{Bai2010} established the \emph{modulus-based matrix splitting (MMS)} framework for solving LCP (\ref{eq:LCP}). Let \(\Omega_1\) and \(\Omega_2\) be nonnegative diagonal matrices, and let \(\Omega\) and \(\Gamma\) be positive diagonal matrices satisfying
\[
\Omega=\Omega_1+\Omega_2.
\]
The author showed that the solution \((w,z)\) yields
\[
x=\frac{1}{2}\left(\Gamma^{-1}z-\Omega^{-1}w\right),
\]
which satisfies the implicit fixed-point equation
\[
\begin{aligned}
	(M\Gamma+\Omega_1)x
	={}(N\Gamma-\Omega_2)x+(\Omega-A\Gamma)|x|-q,
\end{aligned}
\]
where \(A=M-N\) is a matrix splitting of \(A\). Conversely, if $x$ satisfies this equation, then
\begin{equation*}
	z=\Gamma(|x|+x),\qquad
	w=\Omega(|x|-x)
\end{equation*}
provides a solution of the original LCP. Let $\Gamma=\frac{1}{\gamma}I$, as a particular case, the modulus reformulation can be expressed as
\begin{equation*}
	(\Omega+A)x=(\Omega-A)|x|-\gamma q,
\end{equation*}
where $\Omega$ is a positive diagonal matrix and $\gamma>0$. By considering the splitting $A=M-N$, with $M$ nonsingular, the corresponding MMS iteration is obtained as
\begin{equation*}
	\begin{aligned}
		(\Omega+M)x^{(k+1)}
		={}&Nx^{(k)}+(\Omega-A)|x^{(k)}|-\gamma q,\qquad k=0,1,2,\ldots.
	\end{aligned}
\end{equation*}

Based on this fixed-point equation, Bai proposed the
\emph{modulus-based matrix splitting (MMS) iteration}
\begin{equation*}
	(\Omega+M)x^{(k+1)}
	=
	Nx^{(k)}
	+
	(\Omega-A)|x^{(k)}|
	-\gamma q,
	~ k=0,1,2,\ldots,
\end{equation*}
and the solution of the LCP is recovered from
\begin{equation*}
	z^{(k+1)}
	=
	\frac{1}{\gamma}
	\left(
	|x^{(k+1)}|+x^{(k+1)}
	\right).
\end{equation*}
By choosing different matrix splittings of $A$, several
modulus-based stationary iteration methods are obtained as special
cases of the MMS framework. In particular, for
\begin{equation*}
	A=D-L-U,
\end{equation*}
the modulus-based Jacobi (MJ) method follows from
$
M=D,
~
N=L+U,
$
and is given by
\begin{equation*}
	(\Omega+D)x^{(k+1)}
	=
	(L+U)x^{(k)}
	+
	(\Omega-A)|x^{(k)}|
	-\gamma q.
\end{equation*}
For the modulus-based Gauss-Seidel (MGS) method, taking
$
M=D-L,
~
N=U,
$
gives
\begin{equation*}
	(\Omega+D-L)x^{(k+1)}
	=
	Ux^{(k)}
	+
	(\Omega-A)|x^{(k)}|
	-\gamma q.
\end{equation*}
For the modulus-based successive over-relaxation (MSOR) method, the
splitting
\[
M=\alpha^{-1}(D-\alpha L),
~
N=\alpha^{-1}\bigl((1-\alpha)D+\alpha U\bigr)
\]
gives
\begin{align*}
	\left(
	\Omega+\alpha^{-1}(D-\alpha L)
	\right)x^{(k+1)}
	={}&
	\alpha^{-1}
	\bigl((1-\alpha)D+\alpha U\bigr)x^{(k)}
	+
	(\Omega-A)|x^{(k)}|
	-\gamma q.
\end{align*}
Similarly, the modulus-based accelerated over-relaxation (MAOR)
method is obtained from the splitting
\[
M=\alpha^{-1}(D-\beta L),
\]
and
\[
N=
\alpha^{-1}
\bigl(
(1-\alpha)D+(\alpha-\beta)L+\alpha U
\bigr),
\]
which leads to
\begin{align*}
	\left(
	\Omega+\alpha^{-1}(D-\beta L)
	\right)x^{(k+1)}
	={}&
	\alpha^{-1}
	\bigl(
	(1-\alpha)D
	+
	(\alpha-\beta)L
	+
	\alpha U
	\bigr)x^{(k)}
	+
	(\Omega-A)|x^{(k)}|
	-\gamma q.
\end{align*}
Thus, the MJ, MGS, MSOR, and MAOR methods are obtained as particular
instances of the general MMS framework through suitable choices of
the matrix splitting. Bai established convergence results for the MMS
iteration when the system matrix is positive-definite and when it is
an $H_{+}$-matrix. These results yield, in particular, convergence
conditions for the symmetric positive-definite and $M$-matrix cases.
Numerical experiments were presented to compare the proposed
modulus-based relaxation methods with projected relaxation methods
and the modified modulus method, and demonstrated the computational
effectiveness of the modulus-based relaxation approach.
	\item
\textbf{Hadjidimos, Lapidakis, and Tzoumas (2012)}~\cite{HadjidimosLapidakisTzoumas2012} developed two iterative methods for solving LCP (\ref{eq:LCP}) with an \(H_+\)-matrix: the \emph{scaled extrapolated block modulus algorithm (SEBMA)} and an improved modulus-based matrix splitting modified AOR (MAOR) method. For the first approach, they introduced the change of variables
\[
z=|x|+x,\qquad w=|x|-x.
\]
This leads to the modulus fixed-point equation
\begin{equation*}
	x=f(x)=T|x|+b,
	~
	T=(I+A)^{-1}(I-A),
	~
	b=-(I+A)^{-1}q.
\end{equation*}
To scale and extrapolate the LCP, they first introduced
$
	E=(\operatorname{diag}(A))^{-1}
$
and then used a diagonal matrix
$
	\Omega=\operatorname{diag}(\omega_1,\ldots,\omega_n),
	~
	\omega_i>0.
$
With
\begin{equation*}
	A'=\Omega EA\Omega^{-1},
	~
	z'=\Omega z,
	~
	q'=\Omega Eq,
\end{equation*}
the transformed LCP is
\begin{equation*}
	w'=A'z'+q',
\end{equation*}
where
$%
	w'=\Omega Ew.
$
Using the same modulus transformation for the scaled problem, the
corresponding fixed-point equation becomes
\begin{equation*}
	x'=T'|x'|+b',
\end{equation*}
where
\begin{equation*}
	T'=(I+A')^{-1}(I-A'),
	~
	b'=-(I+A')^{-1}\Omega Eq.
\end{equation*}
The resulting scaled extrapolated modulus iteration is
\begin{equation*}
	x'^{(k+1)}
	=
	T'|x'^{(k)}|+b',
	~
	k=0,1,2,\ldots .
\end{equation*}
The authors construct $\Omega$ so that $A'$ is strictly diagonally
dominant with unit diagonal. Under this condition, the above iteration
converges. The scaled extrapolation is then incorporated into
Kappel and Watson's block modulus algorithm, leading to the
\emph{scaled extrapolated block modulus algorithm (SEBMA)}. In each
cycle, all components satisfying the corresponding sign-preserving
condition can be identified together, after which the dimension of the
LCP is reduced.

For the second approach, the authors considered Bai's modulus-based
matrix splitting iteration. Using the substitution
\begin{equation*}
	z=\frac{1}{2}\left(|x|+x\right),
	~
	w=\frac{1}{2}\Omega\left(|x|-x\right),
\end{equation*}
and the splitting
$
	A=M-N,
$
the LCP can be rewritten as
\begin{equation*}
	(\Omega+M)x
	=
	Nx+(\Omega-A)|x|-2q.
\end{equation*}
This gives the modulus-based matrix splitting iteration
\begin{equation*}
	(\Omega+M)x^{(k+1)}
	=
	Nx^{(k)}
	+
	(\Omega-A)|x^{(k)}|
	-
	2q,
	~
	k=0,1,2,\ldots .
\end{equation*}
They considered the modified AOR splitting
\begin{equation*}
	M=\frac{1}{\alpha}(D-\beta L),
\end{equation*}
and
\begin{equation*}
	N=
	\frac{1}{\alpha}
	\left[
	(1-\alpha)D
	+(\alpha-\beta)L
	+\alpha U
	\right],
\end{equation*}
where
$
	A=D-L-U.
$
Consequently, the modulus-based MAOR iteration takes the form
\begin{align}
	\left(
	\Omega+\frac{1}{\alpha}(D-\beta L)
	\right)x^{(k+1)}
	={}&
	\frac{1}{\alpha}
	\left[
	(1-\alpha)D
	+(\alpha-\beta)L
	+\alpha U
	\right]x^{(k)}
	\nonumber\\
	&+
	(\Omega-A)|x^{(k)}|
	-
	2q.
\end{align}
The authors improved the convergence conditions for this MAOR
iteration. In particular, if $A$ is an $H_{+}$-matrix and
\begin{equation*}
	\rho=\rho(D^{-1}|B|)<1,
	~
	B=L+U,
\end{equation*}
then, under
$
	\Omega\geq\frac{1}{\alpha}D,
	~
	\alpha>0,
$
the improved MAOR method (a form of Theorem 4.4 in \cite{Bai2010})converges for
\begin{equation*}
	\alpha\in
	\left(0,\frac{2}{1+\rho}\right),
\end{equation*}
together with the corresponding admissible ranges of $\beta$ given
in the paper. In particular, for
$
	\beta\in[0,\alpha],
$
the above interval for $\alpha$ is sufficient.

The authors further showed that, for fixed admissible $\alpha$ and
$\Omega$, the spectral radius of the improved MAOR iteration decreases
as $\beta$ increases on $(0,\alpha]$. Hence, within this part of the
parameter range, the modified SOR (MSOR) method is the best member
of the MAOR class. Taking
$
	\beta=\alpha=1
$
gives the modified Gauss-Seidel (MGS) method, which is shown to be
the best MSOR method under the stated assumptions. In the numerical
experiments, the improved MGS with
$
	\Omega=\operatorname{diag}(A)
$
was found to be the best among the improved MAOR methods considered.

Finally, the numerical experiments compared SEBMA, MGS, and the
improved MGS method. The results showed that the improved MGS
outperformed MGS in the test problems, while SEBMA was more effective
for larger problem sizes in several of the reported examples. The
authors also observed that the experimentally optimal values of
$\alpha$ and $\beta$ were generally close to, or equal to, $1$.

\item \textbf{Li (2013)}~\cite{Li2013} proposed a general modulus-based matrix splitting iteration method for solving LCP (\ref{eq:LCP}) by introducing two positive diagonal parameter matrices, \(\Omega_1\) and \(\Omega_2\). By introducing the transformation $z=\Omega_1(|x|+x)$, and $w=\Omega_2(|x|-x)$, the method is based on the matrix splitting
\[
A\Omega_1=M_{\Omega_1}-N_{\Omega_1}.
\]
The proposed method is based on the fixed-point equation
\[
(\Omega_2+M_{\Omega_1})x
=
N_{\Omega_1}x
+
(\Omega_2-A\Omega_1)|x|-q.
\]
Accordingly, Method 2.2 computes $x^{(k+1)}$ by solving
\begin{equation*}
	\begin{aligned}
		(\Omega_2+M_{\Omega_1})x^{(k+1)}
		={}&N_{\Omega_1}x^{(k)}
		&+(\Omega_2-A\Omega_1)|x^{(k)}|-q,
	\end{aligned}
\end{equation*}
and then recovers the LCP solution approximation using
\begin{equation*}
	z^{(k+1)}
	=
	\Omega_1\bigl(|x^{(k+1)}|+x^{(k+1)}\bigr) ~~\mbox{and}~~ w^{(k+1)}
	={\Omega_2}
	\left(
	|x^{(k+1)}|-x^{(k+1)}
	\right).
\end{equation*}
This formulation generalizes the standard modulus-based matrix splitting iteration method through the additional parameter matrix $\Omega_1$, allowing the splitting matrices to depend more generally on this parameter.

The authors established convergence when $A$ is an $H_+$-matrix and $A\Omega_1=M_{\Omega_1}-N_{\Omega_1}$ is an $H$-splitting. In particular, the iteration converges to the unique solution of the LCP for any initial vector $x^{(0)}$, provided
\begin{equation*}
	\Omega_2 e
	>
	D\Omega_1 e
	-
	D^{-1}
	\bigl(\langle M_{\Omega_1}\rangle-|N_{\Omega_1}|\bigr)De,
\end{equation*}
for any positive diagonal matrix $D$ such that
$
\bigl(\langle M_{\Omega_1}\rangle-|N_{\Omega_1}|\bigr)D
$
is strictly diagonally dominant. The paper also discusses special choices of the parameter matrices that recover existing methods.  The convergence analysis establishes that the proposed method converges to the unique solution of the LCP when $A$ is an $H_+$-matrix and $A\Omega_1=M_{\Omega_1}-N_{\Omega_1}$ is an $H$-splitting, under the stated sufficient conditions. The numerical example illustrates the effect of the parameter matrix $\Omega_1$ on the spectral radius of the iteration matrix, demonstrating that a suitable choice of this parameter can improve the convergence behavior.

\item \textbf{Zhang and Ren (2013)}~\cite{ZhangRen2013} improved the convergence theorems for the modulus-based matrix splitting (MMS) \cite{Bai2010} and matrix two-stage splitting iterative methods for solving LCP (\ref{eq:LCP}). Their work weakened the assumptions required for convergence, thereby extending the class of matrix splittings for which convergence can be guaranteed. The basic MMS method employs a matrix splitting
\[
A=M-N,
\]
where \(M\) is nonsingular, to construct an iterative scheme based on the modulus reformulation of the LCP. The resulting iteration is given by
\begin{equation*}
	\begin{aligned}
		(\Omega+M)x^{(k+1)}
		={}&Nx^{(k)}
		+(\Omega-A)|x^{(k)}|-\gamma q,
		~ k=0,1,2,\ldots,
	\end{aligned}
\end{equation*}
where \(\Omega\) is a positive diagonal matrix and \(\gamma>0\) is a parameter. The solution of the original LCP is recovered through
\begin{equation*}
	z^{(k+1)}
	=\frac{1}{\gamma}
	\left(|x^{(k+1)}|+x^{(k+1)}\right) ~~\mbox{and}~~ w^{(k+1)}
	=\frac{\Omega}{\gamma}
	\left(
	|x^{(k+1)}|-x^{(k+1)}
	\right).
\end{equation*}

To further develop this approach, the authors introduced the modulus-based matrix two-stage splitting iteration method, in which the splitting of \(M\) is incorporated into the iterative process. Specifically, let
\[
A=M-N,\qquad M=F-G,
\]
where \(M\) is nonsingular. Given an initial vector \(x^{(0)}\), the method computes \(x^{(k+1)}\) through \(\mu_k\) inner iterations:
\begin{equation*}
	(\Omega+F)x^{(k,j+1)}
	=Gx^{(k,j)}+b^{(k)},
\end{equation*}
for \(j=0,1,\ldots,\mu_k-1\), where
$
	b^{(k)}
	=Nx^{(k)}+(\Omega-A)|x^{(k)}|-\gamma q,
	\quad
	x^{(k,0)}=x^{(k)}.
$
After completing the inner iterations, the outer iterate is set to
\(x^{(k+1)}=x^{(k,\mu_k)}\), and the corresponding LCP iterate is obtained from
\begin{equation*}
	z^{(k+1)}
	=\frac{1}{\gamma}
	\left(|x^{(k+1)}|+x^{(k+1)}\right).
\end{equation*}
This two-stage strategy separates the outer modulus-based iteration from the inner splitting iteration, allowing the inner system to be solved iteratively rather than requiring a direct solution at every outer step. The process is repeated until the sequence of LCP iterates converges.

The principal theoretical contribution of Zhang and Ren was to
improve the convergence results for these iterative schemes. In
particular, for an $H_{+}$-matrix $A$, the earlier convergence analyses
often required the splitting
$N
$
to be an $H$-compatible splitting. Zhang and Ren weakened this
requirement by replacing the $H$-compatible splitting assumption with
the more general $H$-splitting assumption. 
Under the corresponding $H$-splitting assumptions, the authors
established convergence results for the MMS iteration and for the
matrix two-stage splitting iteration. Thus, the paper extended the
convergence theory of modulus-based iterative methods for LCPs by
allowing weaker assumptions on the underlying matrix splittings.
\item 
\textbf{Zheng and Yin (2013)}
\cite{ZhengYin2013} proposed a class of accelerated modulus-based matrix splitting iteration methods for solving  LCP (\ref{eq:LCP}). Their approach reformulated the problem as a general implicit fixed-point equation and employed two matrix splittings of the system matrix. Let
\begin{equation*}
	A=M_1-N_1=M_2-N_2
\end{equation*}
be two splittings of $A$. Let $\Omega_1$ and $\Omega_2$ be positive
diagonal matrices, and let $\Gamma$ be a positive diagonal matrix.
Then the implicit fixed-point equation considered by the authors is
\begin{equation*}
	\begin{aligned}
		(M_1\Gamma+\Omega_1)x
		={}&
		(N_1\Gamma-\Omega_2)x
		+
		(\Omega-M_2\Gamma)|x|
		+
		N_2\Gamma |x|
		-q,
	\end{aligned}
\end{equation*}
where
$
\Omega=\Omega_1+\Omega_2.
$
The solution of this fixed-point equation gives the solution of the
original LCP through
$
	z=\Gamma(|x|+x),
$
and $
	w=\Omega(|x|-x).
$
For computational purposes, by taking
\[
\Omega_1=\Omega,
~
\Omega_2=0,
~
\Gamma=\frac{1}{\gamma}I,
~ \gamma>0,
\]
the above equation is reduced to
\begin{equation*}
	(M_1+\gamma\Omega)x
	=
	N_1x
	+
	(\gamma\Omega-M_2)|x|
	+
	N_2|x|
	-\gamma q.
\end{equation*}
Accordingly, the accelerated modulus-based matrix splitting
(AMMS) iteration is given by
\begin{equation*}
	(M_1+\gamma\Omega)x^{(k+1)}
	=
	N_1x^{(k)}
	+
	(\gamma\Omega-M_2)|x^{(k)}|
	+
	N_2|x^{(k+1)}|
	-\gamma q,
	~ k=0,1,2,\ldots .
\end{equation*}
After convergence, the solution approximation of the LCP is obtained
from
\begin{equation*}
	z^{(k+1)}
	=
	\frac{1}{\gamma}
	\left(
	|x^{(k+1)}|+x^{(k+1)}
	\right).
\end{equation*}
The term $N_2|x^{(k+1)}|$ uses the latest available iterate in the
modulus term and constitutes the acceleration mechanism of the
proposed scheme. By choosing suitable pairs of matrix splittings,
several accelerated modulus-based iterative methods are obtained as
special cases, including the accelerated modulus-based Jacobi (AMJ),
accelerated modulus-based Gauss-Seidel (AMGS), accelerated
modulus-based successive overrelaxation (AMSOR), and accelerated
modulus-based accelerated overrelaxation (AMAOR) methods.

For example, when
$
A=D-L-U,
$
the accelerated modulus-based Gauss-Seidel method is obtained by
choosing
\[
M_1=D-L,~ N_1=U,
~
M_2=D,~ N_2=L+U,
\]
while suitable parameterized choices of the splittings give the
AMSOR and AMAOR methods.

The authors established convergence conditions for the accelerated
modulus-based matrix splitting methods when the system matrix is a
positive-definite matrix or an $H_{+}$-matrix. In particular, the
convergence of the proposed iteration is characterized in terms of
the associated iteration matrix, and for the accelerated modulus-based
AOR method the iteration parameters can be selected by minimizing the
spectral radius of the corresponding iteration matrix.

Numerical experiments were presented for  linear
complementarity problems. The numerical results showed that the
accelerated modulus-based methods can reduce the number of iterations
and CPU time compared with the corresponding modulus-based matrix
splitting methods.

\item 
\textbf{Cvetkovi\'c, Hadjidimos, and Kosti\'c (2014)}~\cite{CvetkovicHadjidimosKostic2014} investigated the convergence and optimal parameter selection of modulus-based matrix splitting methods for solving LCP (\ref{eq:LCP}), where \(A\) is a nonsingular \(H_+\)-matrix. They considered the modified accelerated overrelaxation (MAOR) method based on the single-splitting approach. Introducing the transformations
\[
w=\Omega(|x|-x),\qquad z=|x|+x,
\]
they obtained the fixed-point equation
\[
(\Omega+M)x=Nx+(\Omega-A)|x|-q,
\]
which yields the iterative scheme
\[
(\Omega+M)x^{(k+1)}
=Nx^{(k)}+(\Omega-A)|x^{(k)}|-q.
\]
Here, $\Omega$ is a positive diagonal matrix, and the splitting
$A=M-N$ is defined by
\[
M=\frac{1}{\alpha}(D-\beta L),\qquad
N=\frac{1}{\alpha}
\big[(1-\alpha)D+(\alpha-\beta)L+\alpha U\big],
\]
where $A=D-L-U$, $D=\operatorname{diag}(A)$, and $L$ and $U$
are the strictly lower and upper triangular parts of $A$,
respectively. The parameter $\alpha$ is the acceleration parameter,
whereas $\beta$ is the relaxation parameter.
The MAOR framework gives the following modified iterative methods
for particular choices of the parameters $\alpha$ and $\beta$:

	 Modified Jacobi (MJ) method:
	$\alpha=1,\ \beta=0$.
	 Modified Extrapolated Jacobi (MEJ) method:
	$\beta=0$ with general $\alpha$.
	 Modified Gauss--Seidel (MGS) method:
	$\alpha=1,\ \beta=1$.
	 Modified Extrapolated Gauss--Seidel (MEGS) method:
	$\beta=1$ with general $\alpha$.
	 Modified SOR (MSOR) method:
	$\beta=\alpha$.
	 Modified Extrapolated SOR (MESOR) method:
	general choices of $\alpha$ and $\beta$.

The authors further showed that the MGS method, corresponding to
$\alpha=\beta=1$, has the smallest spectral radius among the
considered MAOR methods.

 For a nonsingular $H_+$-matrix $A=D-L-U$ and a diagonal
parameter $\Omega\geq D$, they derived sufficient convergence
ranges for the acceleration parameter $\alpha$ and the relaxation
parameter $\beta$. These ranges are expressed in terms of
$\rho=\rho(D^{-1}|L+U|)<1$ and cover two cases: 
$0\leq\beta/\alpha\leq1$ with $0<\alpha<1/\rho$, and
$\beta/\alpha\geq1$ with
$\delta\rho<\alpha\leq\delta\leq\beta$ and $\delta\rho<1$,
where $\delta$ is determined by the corresponding spectral-radius
condition involving $L$ and $U$. The authors further showed that,
for fixed admissible values of $\alpha$ and $\beta$, choosing
$\Omega=D$ minimizes the spectral radius of the iteration matrix
among all diagonal choices satisfying $\Omega\geq D$. Their
analysis therefore provides sufficient convergence conditions
and identifies the optimal diagonal parameter for the MAOR method.

\item
\textbf{Xu and Liu (2014)}~\cite{XuLiu2014} proposed a modified general modulus-based matrix splitting iteration method for solving  LCP (\ref{eq:LCP}) with an \(H_+\)-matrix. Their method extends the general modulus-based matrix splitting framework by allowing the diagonal parameter matrices to be nonnegative rather than strictly positive.

Let
\begin{equation*}
	A\Omega_{1}
	=
	M_{\Omega_{1}}-N_{\Omega_{1}},
\end{equation*}
where $\Omega_{1}$ is a nonnegative diagonal matrix and
$M_{\Omega_{1}}$ is such that the required coefficient matrix is
nonsingular. Let $\Omega_{2}$ be another nonnegative diagonal matrix.
The corresponding implicit fixed-point equation is written in the form
\begin{equation*}
	(\Omega_{2}+M_{\Omega_{1}})x
	=
	N_{\Omega_{1}}x
	+
	(\Omega_{2}-A\Omega_{1})|x|
	-q.
\end{equation*}
Based on this formulation, the modified general modulus-based matrix
splitting iteration is obtained as
\begin{equation*}
	(\Omega_{2}+M_{\Omega_{1}})x^{(k+1)}
	=
	N_{\Omega_{1}}x^{(k)}
	+
	(\Omega_{2}-A\Omega_{1})|x^{(k)}|
	-q,
	~ k=0,1,2,\ldots,
\end{equation*}
and the solution of the original LCP is recovered from
\begin{equation*}
	z^{(k+1)}
	=
	\Omega_{1}
	\left(
	|x^{(k+1)}|+x^{(k+1)}
	\right) ~~\mbox{and}~~ w^{(k+1)}
	={\Omega_2}
	\left(
	|x^{(k+1)}|-x^{(k+1)}
	\right).
\end{equation*}
The authors established convergence results for the proposed modified
general modulus-based matrix splitting method when the coefficient
matrix is an $H_{+}$-matrix and showed that the resulting convergence
theory improves corresponding existing results. Numerical examples
were presented to illustrate the effectiveness of the proposed method.
\item
\textbf{Xu and Zhu (2019)}~\cite{XuZhu2019} published a corrigendum to the modified general modulus-based matrix splitting method proposed by Xu and Liu~\cite{XuLiu2014}. They identified a flaw in the original \emph{Method 2.1} and its convergence proof and provided a corrected proof for the same iterative method. Their work focused on correcting the theoretical analysis rather than introducing a new modulus-based iteration scheme. Specifically, let \(\Omega_1\) and \(\Omega_2\) be diagonal matrices associated with LCP (\ref{eq:LCP}). For the splitting
\begin{equation*}
	A\Omega_1=M_{\Omega_1}-N_{\Omega_1},
\end{equation*}
the modulus-based fixed-point equation considered in the original
method is
\begin{equation*}
	(\Omega_2+M_{\Omega_1})x
	=
	N_{\Omega_1}x
	+
	(\Omega_2-A\Omega_1)|x|
	-q.
\end{equation*}
Accordingly, the corresponding iteration is
\begin{equation*}
	\begin{aligned}
		(\Omega_2+M_{\Omega_1})x^{(k+1)}
		={}&
		N_{\Omega_1}x^{(k)}
		+
		(\Omega_2-A\Omega_1)|x^{(k)}|
		-q,
		~ k=0,1,2,\ldots .
	\end{aligned}
\end{equation*}
The solution of the original LCP is recovered from
\begin{equation*}
	z^{(k+1)}
	=
	\Omega_1
	\left(
	|x^{(k+1)}|+x^{(k+1)}
	\right) ~~\mbox{and}~~ w^{(k+1)}
	={\Omega_2}
	\left(
	|x^{(k+1)}|-x^{(k+1)}
	\right).
\end{equation*}
The corrigendum did not modify this iterative scheme. Instead, it
provided a corrected proof for the convergence result associated with
\emph{Method 2.1}. Hence, the principal contribution of Xu and Zhu
(2019) is the correction and validation of the convergence analysis of
the modified general modulus-based matrix splitting iteration
originally proposed by Xu and Liu (2014).
\item
\textbf{Zheng and Yin (2014)}~\cite{Zheng2014} further investigated the convergence of accelerated modulus-based matrix splitting iteration methods for solving  LCP (\ref{eq:LCP}). The accelerated framework had been introduced by the authors
using two splittings of the system matrix,
\[
A=M_1-N_1=M_2-N_2.
\]
For a positive diagonal matrix $\Omega$ and $\gamma>0$, the accelerated
modulus-based matrix splitting (AMMS) iteration is
\begin{equation*}
	(\Omega+M_1)x^{(k+1)}
	=
	N_1x^{(k)}
	+
	(\Omega-M_2)|x^{(k)}|
	+
	N_2|x^{(k+1)}|
	-\gamma q,
\end{equation*}
where the solution of the original LCP is recovered from
\begin{equation*}
	z^{(k+1)}
	=
	\frac{1}{\gamma}
	\left(
	|x^{(k+1)}|+x^{(k+1)}
	\right) ~~\mbox{and}~~ w^{(k+1)}
	=\frac{\Omega}{\gamma}
	\left(
	|x^{(k+1)}|-x^{(k+1)}
	\right).
\end{equation*}
The term $N_2|x^{(k+1)}|$ incorporates the latest available
information into the modulus term and is responsible for the
acceleration of the iteration.

By choosing suitable pairs of splittings, the general AMMS framework
gives several accelerated modulus-based methods. In particular, let
\[
A=D-L-U.
\]
The accelerated modulus-based AOR (AMAOR) method is obtained by taking
\[
M_1=\frac{1}{\alpha}(D-\beta L),
\]
\[
N_1=
\frac{1}{\alpha}
\left[
(1-\alpha)D
+
(\alpha-\beta)L
+
\alpha U
\right],
\]
together with
\[
M_2=D-U,
~
N_2=L.
\]
Consequently, the AMAOR iteration is
\begin{align*}
	\left[
	\Omega+\frac{1}{\alpha}(D-\beta L)
	\right]x^{(k+1)}
	={}&
	\frac{1}{\alpha}
	\left[
	(1-\alpha)D
	+
	(\alpha-\beta)L
	+
	\alpha U
	\right]x^{(k)}
	\nonumber\\
	&+
	(\Omega-D+U)|x^{(k)}|
	+
	L|x^{(k+1)}|
	-\gamma q .
\end{align*}
Equivalently, after multiplying the equation by $\alpha$,
\begin{align*}
	(D+\alpha\Omega-\beta L)x^{(k+1)}
	={}&
	\left[
	(1-\alpha)D
	+
	(\alpha-\beta)L
	+
	\alpha U
	\right]x^{(k)}
	\nonumber\\
	&+
	\alpha(\Omega-D+U)|x^{(k)}|
	+
	\alpha L|x^{(k+1)}|
	-\alpha\gamma q .
\end{align*}
The corresponding accelerated modulus-based SOR (AMSOR), Gauss-Seidel
(AMGS), and Jacobi (AMJ) methods are obtained by taking
\[
\alpha=\beta,~
\alpha=\beta=1,~
\text{and}~
\alpha=1,\ \beta=0,
\]
respectively.

In this paper, the authors further investigated the convergence
properties of these accelerated modulus-based matrix splitting
methods. In particular, convergence conditions were established when
the system matrix is either positive definite or an $H_{+}$-matrix.
For the accelerated modulus-based AOR method, the iteration parameters
were further analyzed by minimizing the spectral radius of the
corresponding iteration matrix, which provides optimal choices of the
parameters within the considered class. Numerical experiments
confirmed the theoretical results and showed that the accelerated
methods can reduce both the number of iterations and the CPU time
relative to the corresponding modulus-based matrix splitting methods.

\item
\textbf{Fang and Wei (2015)}~\cite{FangWei2015} proposed the 
\emph{general modulus-based Jacobi (GMJ) iteration method} for solving 
LCP (\ref{eq:LCP}). Let 
$
A=M-N, 
$
where $M$ is a symmetric matrix. Let $\Omega$ be a positive diagonal 
matrix and let $\gamma$ be a positive constant. Based on the 
modulus-based fixed-point formulation
\begin{equation*} 
	(M+\Omega)x 
	= 
	Nx+(\Omega-A)|x|-\gamma q,
\end{equation*} 
the authors defined the GMJ method as follows. Given an initial vector 
$x^{(0)}\in\mathbb{R}^{n}$, compute $x^{(k+1)}$ from
\begin{equation*} 
	(M+\Omega)x^{(k+1)} 
	= 
	Nx^{(k)} 
	+ 
	(\Omega-A)|x^{(k)}| 
	-\gamma q,
	~ k=0,1,2,\ldots,
\end{equation*} 
and then set
\begin{equation*} 
	z^{(k+1)} 
	= 
	\frac{1}{\gamma}
	\left(
	|x^{(k+1)}|+x^{(k+1)}
	\right) ~~\mbox{and}~~ w^{(k+1)}
	=\frac{\Omega}{\gamma}
	\left(
	|x^{(k+1)}|-x^{(k+1)}
	\right).
\end{equation*} 
Here, $\Omega$ is a selected positive diagonal matrix and $\gamma>0$. 
The authors noted that, when the GMJ iteration converges, its limit 
gives a solution of the original LCP. If the LCP has a unique solution, 
then the limit is necessarily that solution.

For the convergence analysis, let $\|\cdot\|$ denote the 
matrix norm induced by a monotonous vector norm and define
\begin{equation*} 
	\delta(\Omega)
	=
	2\left\| |
	(M+\Omega)^{-1}N |
	\right \|
	+
	\left\| |
	(M+\Omega)^{-1}(M-\Omega)|
	\right\|.
\end{equation*} 
The authors proved that, if the LCP has a unique solution and 
$\delta(\Omega)<1$, then the GMJ iteration converges to the solution 
for any initial vector. They further considered the case in which $M$ 
is symmetric positive definite and $\Omega=\omega I$, and established 
convergence for
\begin{equation*}
	\omega\in(\tau\mu_{\max},+\infty)
\end{equation*}
under the conditions
\begin{equation*}
	\tau<1,
	~
	\tau\mu_{\max}<\mu_{\min},
\end{equation*}
with the optimal parameter
\begin{equation*}
	\omega_{\mathrm{opt}}
	=
	\sqrt{\mu_{\min}\mu_{\max}}.
\end{equation*}
For the case in which $A$ is an $H_{+}$-matrix, taking
$M=\operatorname{diag}(A)$ and $\Omega=\omega I$, the authors also 
derived a corresponding convergence condition with
\begin{equation*}
	\tau<1,
	~
	\tau^{2}\mu_{\max}<\mu_{\min},
\end{equation*}
where $\mu_{\min}$ and $\mu_{\max}$ are the smallest and largest 
diagonal elements of $M$, respectively. 

Numerical experiments were carried out to compare the GMJ method with
the modulus (M), modified modulus (MM), and modulus-based
Gauss-Seidel (MGS) methods. The numerical results showed that the
GMJ method can be superior to some of the modulus-related methods in
certain situations. In particular, when the diagonal elements of the
system matrix are different, the GMJ method with a suitable diagonal
parameter matrix $\Omega$ was shown to perform well.
\item
\textbf{Zheng and Li (2015)}~\cite{ZhengLi2015}
developed a modulus-based nonsmooth Newton method for LCP (\ref{eq:LCP}). They started from the equivalent modulus
equation
\begin{equation*}
	(A+I)x+(A-I)|x|+q=0,
\end{equation*}
where, once the solution $x^{*}$ of this equation is obtained, the
solution of the original LCP is recovered from
$
	z^{*}=x^{*}+|x^{*}|.
$
Define
\begin{equation*}
	F(x)
	=
	(A+I)x+(A-I)|x|+q.
\end{equation*}
Since the modulus function is nonsmooth at zero, the authors applied
the nonsmooth Newton method to $F(x)=0$. Let
$
	D(x)
	=
	\operatorname{diag}
	\left(
	\operatorname{sgn}(x_{1}),
	\ldots,
	\operatorname{sgn}(x_{n})
	\right).
$
A corresponding element of the generalized Jacobian of $F$ is
\begin{equation*}
	V(x)
	=
	A+I+(A-I)D(x).
\end{equation*}
Consequently, the modulus-based nonsmooth Newton iteration is
\begin{equation*}
	x^{(k+1)}
	=
	x^{(k)}
	-
	V(x^{(k)})^{-1}
	F(x^{(k)}),
\end{equation*}
that is,
\begin{equation*}
	x^{(k+1)}
	=
	x^{(k)}
	-
	V(x^{(k)})^{-1}
	\left[
	(A+I)x^{(k)}
	+
	(A-I)|x^{(k)}|
	+
	q
	\right].
\end{equation*}
After convergence, the corresponding LCP solution is obtained from
\begin{equation*}
	z^{(k+1)}
	=
	|x^{(k+1)}|+x^{(k+1)}.
\end{equation*}
The authors established locally quadratic convergence of the
modulus-based nonsmooth Newton method under suitable conditions.
They further showed that, with an appropriately chosen initial vector
and generalized Jacobian, the method achieves local one-step
convergence. In particular, when the initial vector satisfies the
required sign-pattern conditions with respect to the solution, the
next Newton iterate coincides with the solution. These results
characterize the local convergence behavior of the proposed
modulus-based nonsmooth Newton method.

	\item
	\textbf{Liu, Zheng, and Li (2016)}
	~\cite{LiuZhengLi2016} proposed a \emph{general accelerated
		modulus-based matrix splitting iteration method} for solving LCP(\ref{eq:LCP}). The proposed framework generalizes the
	general modulus-based matrix splitting iteration methods as well as
	the accelerated modulus-based matrix splitting iteration methods.
	
	Let $\Omega_1$ and $\Omega_2$ be two $n\times n$ positive diagonal
	matrices. Consider two splittings of $A\Omega_1$,
	\begin{equation*}
		A\Omega_1
		=
		M_{\Omega_1}^{(1)}-N_{\Omega_1}^{(1)}
		=
		M_{\Omega_1}^{(2)}-N_{\Omega_1}^{(2)}.
	\end{equation*}
	The authors introduced the following fixed-point equation:
	\begin{equation*}
		\begin{aligned}
			\left(\Omega_2+M_{\Omega_1}^{(1)}\right)x
			={}
			N_{\Omega_1}^{(1)}x
			+
			\left(\Omega_2-M_{\Omega_1}^{(2)}\right)|x|
			+
			N_{\Omega_1}^{(2)}|x|
			-q.
		\end{aligned}
	\end{equation*}
	Based on this equation, the authors proposed the following
	\emph{general accelerated modulus-based matrix splitting iteration
		method} (Method~3.1). Given an initial vector
	$x^{(0)}\in\mathbb{R}^{n}$, for $k=0,1,2,\ldots$, compute
	$x^{(k+1)}$ from
	\begin{equation*}
		\begin{aligned}
			\left(\Omega_2+M_{\Omega_1}^{(1)}\right)x^{(k+1)}
			={}
			N_{\Omega_1}^{(1)}x^{(k)}
			+
			\left(\Omega_2-M_{\Omega_1}^{(2)}\right)|x^{(k)}|
			+
			N_{\Omega_1}^{(2)}|x^{(k+1)}|
			-q,
		\end{aligned}
	\end{equation*}
	where $N_{\Omega_1}^{(2)}$ is required to be a strictly lower
	triangular matrix. The corresponding solution of the original LCP is
	recovered by
	\begin{equation*}
		z^{(k+1)}
		=
		\Omega_1
		\left(
		|x^{(k+1)}|+x^{(k+1)}
		\right) ~~\mbox{and}~~ w^{(k+1)}
		=\Omega_2
		\left(
		|x^{(k+1)}|-x^{(k+1)}
		\right).
	\end{equation*}
	
	For the particular splitting
$
		A\Omega_1
		=
		D_{\Omega_1}
		-
		L_{\Omega_1}
		-
		U_{\Omega_1},
	$
	where $D_{\Omega_1}$, $L_{\Omega_1}$, and $U_{\Omega_1}$ denote the
	diagonal, strictly lower triangular, and strictly upper triangular
	parts of $A\Omega_1$, respectively, the authors obtained a general
	accelerated modulus-based accelerated overrelaxation (GAMAOR) method.
	Taking
	\begin{equation*}
		M_{\Omega_1}^{(1)}
		=
		\frac{1}{\alpha}
		\left(
		D_{\Omega_1}-\beta L_{\Omega_1}
		\right),
	\end{equation*}
	\begin{equation*}
		N_{\Omega_1}^{(1)}
		=
		\frac{1}{\alpha}
		\left[
		(1-\alpha)D_{\Omega_1}
		+
		(\alpha-\beta)L_{\Omega_1}
		+
		\alpha U_{\Omega_1}
		\right],
	\end{equation*}
	and
	$
		M_{\Omega_1}^{(2)}
		=
		D_{\Omega_1}-U_{\Omega_1},
		~
		N_{\Omega_1}^{(2)}
		=
		L_{\Omega_1},
	$
	the GAMAOR iteration is
	\begin{align*}
		\left[
		\Omega_2+
		\frac{1}{\alpha}
		\left(
		D_{\Omega_1}-\beta L_{\Omega_1}
		\right)
		\right]x^{(k+1)}
		={}&
		\frac{1}{\alpha}
		\left[
		(1-\alpha)D_{\Omega_1}
		+
		(\alpha-\beta)L_{\Omega_1}
		+
		\alpha U_{\Omega_1}
		\right]x^{(k)}
		\nonumber\\
		&+
		\left(
		\Omega_2-D_{\Omega_1}+U_{\Omega_1}
		\right)|x^{(k)}|
		+
		L_{\Omega_1}|x^{(k+1)}|
		-q.
	\end{align*}
	The GAMAOR scheme contains the following special cases.
	
	For $\alpha=\beta$, the \emph{general accelerated modulus-based
		successive overrelaxation (GAMSOR)} method is obtained:
	\begin{align*}
		\left[
		\Omega_2+
		\frac{1}{\alpha}
		\left(
		D_{\Omega_1}-\alpha L_{\Omega_1}
		\right)
		\right]x^{(k+1)}
		={}&
		\frac{1}{\alpha}
		\left[
		(1-\alpha)D_{\Omega_1}
		+
		\alpha U_{\Omega_1}
		\right]x^{(k)}
		\nonumber\\
		&+
		\left(
		\Omega_2-D_{\Omega_1}+U_{\Omega_1}
		\right)|x^{(k)}|
		+
		L_{\Omega_1}|x^{(k+1)}|
		-q.
	\end{align*}
	For $\alpha=\beta=1$, the \emph{general accelerated modulus-based
		Gauss-Seidel (GAMGS)} method is obtained:
	\begin{equation*}
		\begin{aligned}
			\left(
			\Omega_2+D_{\Omega_1}-L_{\Omega_1}
			\right)x^{(k+1)}
			={}
			U_{\Omega_1}x^{(k)}
			+
			\left(
			\Omega_2-D_{\Omega_1}+U_{\Omega_1}
			\right)|x^{(k)}
			|
			+
			L_{\Omega_1}|x^{(k+1)}|
			-q.
		\end{aligned}
	\end{equation*}
	For $\alpha=1$ and $\beta=0$, the \emph{general accelerated
		modulus-based Jacobi (GAMJ)} method is obtained:
	\begin{equation*}
		\begin{aligned}
			\left(
			\Omega_2+D_{\Omega_1}
			\right)x^{(k+1)}
			={}
			\left(
			L_{\Omega_1}+U_{\Omega_1}
			\right)x^{(k)}
			+
			\left(
			\Omega_2-D_{\Omega_1}+U_{\Omega_1}
			\right)|x^{(k)}|
			+
			L_{\Omega_1}|x^{(k+1)}|
			-q.
		\end{aligned}
	\end{equation*}
	Thus, Method~3.1 provides a more general framework: if
	$
	M_{\Omega_1}^{(2)}=A\Omega_1,
	~
	N_{\Omega_1}^{(2)}=0,
	$
	it reduces to the general modulus-based matrix splitting iteration
	methods. Moreover, with the special choices
	\[
	\Omega_1=\frac{1}{\gamma}I,
	~
	\Omega_2=\frac{1}{\gamma}\Omega,
	\]
	together with the corresponding scaled splittings, it reduces to the
	accelerated modulus-based matrix splitting iteration methods proposed
	previously.
	
	For the convergence analysis, the authors considered the case in
	which $A$ is an $H_{+}$-matrix. They assumed that
	\[
	A\Omega_1
	=
	M_{\Omega_1}^{(1)}-N_{\Omega_1}^{(1)}
	=
	M_{\Omega_1}^{(2)}-N_{\Omega_1}^{(2)}
	\]
	are two $H$-splittings satisfying the comparison condition
	$
	\left\langle M_{\Omega_1}^{(1)}\right\rangle
	\geq
	\left\langle M_{\Omega_1}^{(2)}\right\rangle.
	$
	Under the corresponding sufficient conditions on the positive
	diagonal matrix $\Omega_2$, the iteration sequence generated by
	Method~3.1 converges to the exact solution of the LCP for any initial
	vector $x^{(0)}\in\mathbb{R}^{n}$. The authors also discussed the
	choice of the parameter matrix $\Omega_1$ and provided a practical
	choice for a special case.
	
	Numerical experiments were presented to compare the proposed
	GAMSOR method with the corresponding accelerated modulus-based SOR
	(AMSOR) method. The numerical results showed that the GAMSOR method
	can be more stable and can require fewer iteration steps and less CPU
	time for the test problems considered.
\item
\textbf{Wu and Li (2016)}~\cite{WuLi2016} proposed the
\emph{two-sweep modulus-based matrix splitting (TMMS) iteration
	method} for solving LCP(\ref{eq:LCP}).
Their approach extends the standard modulus-based matrix splitting
framework by using the information from the two preceding iteration
points in the modulus term.

Let
$
	A=M-N,
$
where $M$ is nonsingular, and let $\Omega$ be a positive diagonal
matrix and $\gamma>0$. The LCP is associated with the implicit
modulus-based fixed-point equation
\begin{equation*}
	(\Omega+M)x
	=
	Nx
	+
	(\Omega-A)|x|
	-\gamma q.
\end{equation*}
Instead of using the modulus of the current iterate in every term,
Wu and Li introduced the two-sweep scheme. Starting from two initial
vectors
$
x^{(0)},x^{(1)}\in\mathbb{R}^{n},
$
the \emph{TMMS iteration} is defined, for $k=1,2,\ldots$, by
\begin{equation*}
	(\Omega+M)x^{(k+1)}
	=
	Nx^{(k)}
	+
	(\Omega-A)|x^{(k-1)}|
	-\gamma q.
\end{equation*}
After each iteration, the approximation to the solution of the
original LCP is recovered from
\begin{equation*}
	z^{(k+1)}
	=
	\frac{1}{\gamma}
	\left(
	|x^{(k+1)}|+x^{(k+1)}
	\right) ~~\mbox{and}~~ w^{(k+1)}
	=\Omega
	\left(
	|x^{(k+1)}|-x^{(k+1)}
	\right).
\end{equation*}
Thus, the principal feature of the TMMS method is that the linear
part uses the current iterate $x^{(k)}$, whereas the modulus term uses
the previous iterate $|x^{(k-1)}|$. Hence, the method combines
information from two successive iteration points and constitutes a
two-sweep extension of the one-sweep modulus-based matrix splitting
iteration.

For the convergence analysis, the authors considered the iteration
matrix associated with the two-sweep recurrence. If $x^{*}$ denotes
the solution of the corresponding fixed-point equation, then the
errors satisfy a two-step relation of the form
\begin{equation*}
	\begin{bmatrix}
		|x^{(k+1)}-x^{*}|\\[2mm]
		|x^{(k)}-x^{*}|
	\end{bmatrix}
	\leq
	\mathcal{L}
	\begin{bmatrix}
		|x^{(k)}-x^{*}|\\[2mm]
		|x^{(k-1)}-x^{*}|
	\end{bmatrix},
\end{equation*}
where $\mathcal{L}$ is the corresponding nonnegative
two-by-two block iteration matrix. Consequently, a sufficient
condition for convergence is
\begin{equation*}
	\rho(\mathcal{L})<1, 
\end{equation*}
  where \begin{equation*}
  	\mathcal{L}=
  	\left|
  	\begin{matrix}
  		(\Omega+M)^{-1}N & (\Omega+M)^{-1}\Omega A\\[2mm]
  		I & O
  	\end{matrix}
  	\right|.
  \end{equation*} 
Wu and Li established convergence results for the TMMS method when
the system matrix is positive definite and when it is an
$H_{+}$-matrix, under suitable assumptions on the matrix splitting
and the positive diagonal matrix $\Omega$. In the $H_{+}$-matrix
case, the convergence analysis was developed by means of
comparison matrices and matrix splitting theory.

The authors also investigated the practical performance of the
proposed TMMS method through numerical experiments. The numerical
results showed the efficiency of the two-sweep scheme for the test
LCPs considered.
 \item
\textbf{Zheng, Li, and Vong (2017)}~\cite{ZhengLiVong2017} proposed a relaxation modulus-based matrix splitting (RMMS) iteration method for solving LCP(\ref{eq:LCP}). Let \(A\Omega_1=M_{\Omega_1}-N_{\Omega_1}\) be a splitting of \(A\Omega_1\), where \(\Omega_1\) and \(\Omega_2\) are positive diagonal matrices, and let \(\theta\neq 0\) be a relaxation parameter. Given an initial vector \(x^{(0)}\), the method computes an intermediate iterate \(x^{(k-\frac12)}\) by solving
\[
\begin{aligned}
	(\Omega_2+M_{\Omega_1})x^{(k-\frac12)}
	={}N_{\Omega_1}x^{(k-1)}
	+(\Omega_2-A\Omega_1)|x^{(k-1)}|-q,
\end{aligned}
\]
and then performs the relaxation step
$
x^{(k)}=(1-\theta)x^{(k-1)}
+\theta x^{(k-\frac12)}.
$
The solution of the original LCP is recovered through
$
z^{(k)}=\Omega_1\bigl(|x^{(k)}|+x^{(k)}\bigr)
$ $~~\mbox{and}~~ w^{(k)}
=\Omega_2
\left(
|x^{(k)}|-x^{(k)}
\right).$
By choosing
\[
\begin{aligned}
	M_{\Omega_1}
	&=\frac{1}{\alpha}
	(D_{A\Omega_1}-\beta L_{A\Omega_1}),\\
	N_{\Omega_1}
	&=\frac{1}{\alpha}\bigl[
	(1-\alpha)D_{A\Omega_1}
	+(\alpha-\beta)L_{A\Omega_1}
	+\alpha U_{A\Omega_1}\bigr],
\end{aligned}
\]
the authors derived the relaxation modulus-based accelerated overrelaxation (RMAOR) method. Its special cases include the relaxation modulus-based successive overrelaxation (RMSOR) method for \(\alpha=\beta\), the relaxation modulus-based Gauss--Seidel (RMGS) method for \(\alpha=\beta=1\), and the relaxation modulus-based Jacobi (RMJ) method for \(\alpha=1,\ \beta=0\). The
authors established convergence conditions for the proposed
relaxation modulus-based matrix splitting method and investigated
strategies for selecting the relaxation parameter. In particular,
they derived sufficient conditions for convergence when the system
matrix is an $H_{+}$-matrix. Numerical experiments were presented to
show that suitable relaxation parameters can accelerate the
convergence of the modulus-based matrix splitting iteration, with
fewer iteration steps and lower CPU time for the test problems
considered.
\item
\textbf{Wen, Zheng, Li, and Peng (2018)}~\cite{WenZhengLiPeng2018} proposed a relaxation modulus-based matrix splitting iteration method for solving LCP(\ref{eq:LCP}). Let \(A=M-N\) be a splitting of \(A\), and let \(\Omega_1\) and \(\Omega_2\) be positive diagonal matrices. The proposed iterative method is given by
\[
\begin{aligned}
	x^{(k)}={}&(I-\Omega)x^{(k-1)}
	+\Omega(\Omega_2+M_{\Omega_1})^{-1}
	\Big[N_{\Omega_1}x^{(k-1)}
	+(\Omega_2-A\Omega_1)|x^{(k-1)}|-q\Big],
\end{aligned}
\]
where \(\Omega\) is a positive diagonal matrix and
$
A\Omega_1=M_{\Omega_1}-N_{\Omega_1}.
$
The approximate solution of the original LCP is recovered by
\[
z^{(k)}
=\Omega_1\left(|x^{(k)}|+x^{(k)}\right)
\]
and
\[
w^{(k)}
=\Omega_2\left(|x^{(k)}|-x^{(k)}\right).
\]

To establish convergence, the authors considered the case in which \(A\) is positive definite. Let \((z^*,w^*)\) be a solution of \(\operatorname{LCP}(q,A)\). The associated vector
\[
x^*=\frac{1}{2}
\left(\Omega_1^{-1}z^*-\Omega_2^{-1}w^*\right)
\]
satisfies the implicit fixed-point equation
\[
\begin{aligned}
	(\Omega_2+M_{\Omega_1})x^*
	={}&N_{\Omega_1}x^*
	+(\Omega_2-A\Omega_1)|x^*|-q.
\end{aligned}
\]

Using the identity \(x^*=(I-\Omega)x^*+\Omega x^*\), the authors obtained
\[
\begin{aligned}
	x^*={}&(I-\Omega)x^*
	+\Omega(\Omega_2+M_{\Omega_1})^{-1}\\
	&\quad\times
	\Big[N_{\Omega_1}x^*
	+(\Omega_2-A\Omega_1)|x^*|-q\Big].
\end{aligned}
\]

Subtracting this equation from the iterative scheme yields the error relation
\[
\begin{aligned}
	x^{(k)}-x^*
	={}&\Big[(I-\Omega)
	+\Omega(\Omega_2+M_{\Omega_1})^{-1}
	N_{\Omega_1}\Big](x^{(k-1)}-x^*)\\
	&+\Omega(\Omega_2+M_{\Omega_1})^{-1}
	(\Omega_2-A\Omega_1)
	\big(|x^{(k-1)}|-|x^*|\big).
\end{aligned}
\]

This error relation is used to derive sufficient convergence conditions for the proposed method. In particular, to establish convergence of the recovered solution sequence \(z^{(k)}\) to \(z^*\), it is sufficient to prove that \(x^{(k)}\to x^*\).

Numerical experiments were conducted for LCPs with positive-definite coefficient matrices to evaluate the performance of the proposed method. The results indicated that, with a suitable choice of the relaxation parameter, the relaxation modulus-based matrix splitting method can improve convergence performance, reducing the number of iterations and CPU time compared with the corresponding non-relaxed modulus-based matrix splitting method.

\item
Li \cite{Li2013} introduced a general modulus-based matrix splitting (GMS) method for solving LCP (\ref{eq:LCP}). Let \(\Omega_1\) and \(\Omega_2\) be positive diagonal matrices. Using the variable transformation
\[
z=\Omega_1(|x|+x),\qquad
r=\Omega_2(|x|-x),
\]
the LCP (\ref{eq:LCP}) can be equivalently reformulated as
\[
(\Omega_2+A\Omega_1)x
=(\Omega_2-A\Omega_1)|x|-q.
\]
Let
\[
A\Omega_1=M\Omega_1-N\Omega_1
\]
be a matrix splitting. The corresponding GMS iteration is given by
\[
\begin{aligned}
	(\Omega_2+M\Omega_1)x^{(k+1)}
	={}&N\Omega_1x^{(k)}+(\Omega_2-A\Omega_1)|x^{(k)}|-q,
\end{aligned}
\]
where \(k=0,1,2,\ldots\).

\textbf{Wu, Peng, and Li (2018)} \cite{WuPengLi2018} proposed a preconditioned general modulus-based matrix splitting (PGMS) iteration method to accelerate the solution of LCP (\ref{eq:LCP}) when the system matrix \(A\) is an \(H_{+}\)-matrix. Let \(P\) be a nonsingular preconditioning matrix, and consider the splitting
\[
PA\Omega_1=M_{P,\Omega_1}-N_{P,\Omega_1}.
\]
The resulting PGMS iteration is defined as
\[
\begin{aligned}
	(P\Omega_2+M_{P,\Omega_1})x^{(k+1)}
	={}&N_{P,\Omega_1}x^{(k)}+P(\Omega_2-A\Omega_1)|x^{(k)}|-Pq.
\end{aligned}
\]
After computing \(x^{(k+1)}\), the LCP variables are recovered as
\[
z^{(k+1)}
=\Omega_1\bigl(|x^{(k+1)}|+x^{(k+1)}\bigr)
\]
and
\[
r^{(k+1)}
=\Omega_2\bigl(|x^{(k+1)}|-x^{(k+1)}\bigr).
\]

Sufficient convergence conditions for the PGMS method are established under appropriate assumptions on the matrix splitting and the \(H_{+}\)-matrix \(A\). The authors also discuss related preconditioned iterative schemes, including the preconditioned general modulus-based accelerated over-relaxation (PGMAOR) method and its special cases.

In their numerical experiments, four preconditioners are considered:
\[
P_1=I+t_1C_1,\qquad
P_2=I+t_2(C_1+C_m),
\]
and
\[
\widehat P_1=I+\widehat t_1C_1^T,\qquad
\widehat P_2=I+\widehat t_2(C_1+C_1^T),
\]
where \(t_1,t_2,\widehat t_1,\widehat t_2>0\). The entries of \(C_i\) are defined by
\[
c_{kj}=
\begin{cases}
	1, & k=j+i,\\
	0, & \text{otherwise}.
\end{cases}
\]
The numerical results indicate that suitable preconditioners can improve the computational performance of the considered iterative methods on the tested problems.

\item
\textbf{Fang and Zhu (2019)}~\cite{FangZhu2019} proposed the \emph{modulus-based matrix double splitting (MB-DS) iteration method} for solving LCP(\ref{eq:LCP}). The method is developed from the accelerated modulus-based matrix splitting framework by employing two splittings of the coefficient matrix,
\[
A=M_1-N_1=M_2-N_2.
\]
Let \(\Omega\) be a positive diagonal matrix and \(\gamma>0\). The authors first obtained the implicit fixed-point equation
\[
\begin{aligned}
	(M_1+\Omega)x
	={}&N_1x+(\Omega-M_2)|x|+N_2|x|-\gamma q.
\end{aligned}
\]
Based on this equation, they proposed the MB-DS iteration method
\[
\begin{aligned}
	(M_1+\Omega)x^{(k+1)}
	={}N_1x^{(k)}
	+(\Omega-M_2)|x^{(k)}|+N_2|x^{(k+1)}|-\gamma q,
\end{aligned}
\]
where \(N_2\) is a nonzero strictly lower triangular matrix and \(M_1\) is a lower triangular matrix. The approximate solution of the original LCP is recovered by
\[
z^{(k+1)}
=\frac{1}{\gamma}
\left(|x^{(k+1)}|+x^{(k+1)}\right)  ~~\mbox{and}~~ w^{(k+1)}
=\Omega
\left(
|x^{(k+1)}|-x^{(k+1)}
\right).
\]
The presence of \(N_2|x^{(k+1)}|\) distinguishes the proposed MB-DS method from the standard modulus-based matrix splitting iteration, in which the modulus term depends only on the previous iterate. Since \(N_2\) is strictly lower triangular and \(M_1\) is lower triangular, the new iterate can be computed successively, component by component.

The authors established convergence conditions for the proposed method when \(A\) is a \(P\)-matrix. In particular, defining
\[
\begin{aligned}
	\delta_S(\Omega)
	={}\|(M_1+\Omega)^{-1}N_1\|+\|(M_1+\Omega)^{-1}N_2\|+\|(M_1+\Omega)^{-1}(\Omega-M_2)\|,
\end{aligned}
\]
and
\[
\begin{aligned}
	\delta_L(\Omega)
	={}2\|(M_1+\Omega)^{-1}N_1\|+2\|(M_1+\Omega)^{-1}N_2\|+\|(M_1+\Omega)^{-1}(\Omega-M_1)\|,
\end{aligned}
\]
they showed that the iteration converges to the unique solution of the LCP if either \(\delta_S(\Omega)<1\) or \(\delta_L(\Omega)<1\). They also proved convergence for \(H_+\)-matrices under the conditions that the two splittings are \(H\)-compatible, \(|N_1|\geq |N_2|\), and \(\Omega\geq\operatorname{diag}(M_2)\).

Furthermore, for particular splittings with \(\Omega=\omega I\), the authors derived sufficient convergence regions for the parameter \(\omega\) and identified optimal parameter values. For a positive diagonal \(M_1\neq kI\), one such optimal value is
\[
\omega^*=\sqrt{\lambda_{\min}\lambda_{\max}},
\]
where \(\lambda_{\min}\) and \(\lambda_{\max}\) are the smallest and largest eigenvalues of \(M_1\), respectively. When \(M_1=kI\), the corresponding optimal value is \(\omega^*=k\), under the conditions of the relevant corollary. Numerical experiments demonstrated that the proposed MB-DS method can be competitive with the modulus-based Gauss--Seidel (MGS) method, sometimes requiring less CPU time and achieving comparable accuracy.

\item {Wu, Peng, and Li }\cite{WuPengLi2018}
proposed a \emph{preconditioned general modulus-based matrix
	splitting (PGMMS) iteration method} for solving LCP(\ref{eq:LCP}). Subsequently, Wu et al.'s PGMMS
method was used as the basis for further development.
In particular, \textbf{Ren, Wang, Tang, and Wang (2019)}
~\cite{RenWangTangWang2019} proposed a
\emph{preconditioned general two-step modulus-based matrix
	splitting (PGTMMS) iteration method} for solving linear
complementarity problems whose coefficient matrix is an
$H_{+}$-matrix.

Let $P$ be a suitable nonsingular preconditioner, and let
$A=M-N$ be a splitting of the coefficient matrix. Let
$\Omega_1$ and $\Omega_2$ be diagonal matrices satisfying
the assumptions of the methods. The PGMMS iteration
proposed by Wu et al. is given by
\begin{equation*}
	\begin{aligned}
		(P\Omega_2+MP\Omega_1)x^{(k+1)}
		={}NP\Omega_1x^{(k)}
		+P(\Omega_2-A\Omega_1)|x^{(k)}|-Pq.
	\end{aligned}
\end{equation*}
The corresponding approximation to the LCP solution is
recovered using
$$
	z^{(k+1)}
	=\Omega_1\bigl(|x^{(k+1)}|+x^{(k+1)}\bigr) ~~\mbox{and}~~ w^{(k+1)}
	=\Omega_2
	\left(
	|x^{(k+1)}|-x^{(k+1)}
	\right).
$$

Building on this existing PGMMS method, Ren et al. developed
the PGTMMS method, which performs two successive matrix
splitting steps in each iteration. First, an intermediate
iterate $x^{(k+\frac12)}$ is computed by solving
\begin{equation*}
	\begin{aligned}
		(P\Omega_2+M_1P\Omega_1)x^{(k+\frac12)}
		={}N_1P\Omega_1x^{(k)}
		+P(\Omega_2-A\Omega_1)|x^{(k)}|-Pq,
	\end{aligned}
\end{equation*}
where $M_1-N_1$ denotes the splitting used in the first
stage. Next, the intermediate iterate is used to compute
$x^{(k+1)}$ from
\begin{equation*}
	\begin{aligned}
		(P\Omega_2+M_2P\Omega_1)x^{(k+1)}
		={}N_2P\Omega_1x^{(k+\frac12)}
		+P(\Omega_2-A\Omega_1)
		|x^{(k+\frac12)}|-Pq,
	\end{aligned}
\end{equation*}
where $M_2-N_2$ denotes the splitting used in the second
stage. The solution approximation is then obtained as
\begin{equation*}
	z^{(k+1)}
	=\Omega_1\bigl(|x^{(k+1)}|+x^{(k+1)}\bigr).
\end{equation*}
Ren et al. established sufficient convergence conditions
for the proposed PGTMMS method under suitable assumptions
when $A$ is an $H_{+}$-matrix. Numerical experiments were
also conducted to evaluate the performance of the proposed
two-step method.

\item
\textbf{Wu and Li (2020)}~\cite{WuLi2020} proposed a
\emph{preconditioned modulus-based matrix splitting iteration method}
for solving LCP(\ref{eq:LCP}) with a
symmetric positive-definite a  Toeplitz matrix. For a given parameter $\alpha>0$, the authors transformed
the LCP into the following equivalent nonlinear equation:
\begin{equation*}
	(\alpha I+A)x=(\alpha I-A)|x|-q.
\end{equation*}
Starting from an initial vector $x^{(0)}\in\mathbb{R}^n$, the
outer iteration is defined by
\begin{equation*}
	(\alpha I+A)x^{(k+1)}
	=(\alpha I-A)|x^{(k)}|-q,
	\qquad k=0,1,2,\ldots.
\end{equation*}
At each outer iteration, the resulting linear system is solved
approximately by the preconditioned conjugate gradient (PCG)
method. To accelerate the solution of these symmetric
positive-definite Toeplitz systems, the authors considered two
circulant preconditioners, namely Strang's circulant
preconditioner and T.~Chan's circulant preconditioner.

For a Toeplitz matrix $A$, Strang's circulant approximation
$S(A)$ is constructed by selecting the central diagonals of
$A$ and arranging them in circulant form. In particular, for
$A=(a_{i-j})$, the entries of $S(A)$ are determined by the
corresponding central Toeplitz diagonals, with the precise
construction depending on whether the matrix order is even
or odd. The second preconditioner, T.~Chan's circulant
preconditioner $C^{F}(A)$, is defined as the optimal circulant
approximation to $A$ in the Frobenius norm:
\begin{equation*}
	C^{F}(A)
	=\underset{C\in\mathcal{C}_n}{\operatorname{arg\,min}}
	\,\|A-C\|_F,
\end{equation*}
where $\mathcal{C}_n$ denotes the set of $n\times n$
circulant matrices. Accordingly, for the matrix
$\alpha I+A$, the two preconditioners used in the PCG method
are
\begin{equation*}
	P_S=S(\alpha I+A),
	\qquad
	P_C=C^{F}(\alpha I+A),
\end{equation*}
respectively. Thus, the proposed method has two implementations,
one using Strang's circulant preconditioner and the other using
T.~Chan's circulant preconditioner.

The convergence analysis is based on the modulus-based
iteration and the spectral norm of its iteration matrix.
A sufficient convergence condition established in the paper is
\begin{equation*}
	\left\|(\alpha I+A)^{-1}(\alpha I-A)\right\|_2<1.
\end{equation*}
Under this condition, the outer iteration converges to the solution of the equivalent modulus-based equation. After computing \(x^{(k+1)}\), the corresponding LCP variables are recovered by
\begin{equation*}
	z^{(k+1)}=|x^{(k+1)}|+x^{(k+1)},
\end{equation*}
and
\begin{equation*}
	w^{(k+1)}=Az^{(k+1)}+q.
\end{equation*}
The iteration terminates when the residual
\begin{equation*}
	\operatorname{RES}_{k+1}
	=\left\|\begin{pmatrix}
		z^{(k+1)}\\
		w^{(k+1)}
	\end{pmatrix}\right\|_2
\end{equation*}
satisfies the prescribed tolerance. Numerical experiments were conducted for large-scale Toeplitz problems to compare the unpreconditioned implementation with the Strang \cite{Strang1986} and T. Chan \cite{Chan1988} circulant preconditioners. The results showed that these preconditioners reduce the number of inner PCG iterations and the computational time in the tested examples, with T. Chan's preconditioner generally requiring fewer inner iterations than Strang's preconditioner.
	\item
\textbf{Wu, Zhang, and Yang (2022)}~\cite{WuZhangYang2022}
proposed a modulus-based inexact non-alternating preconditioned
splitting (MINPS) iteration method for solving LCP
(\ref{eq:LCP}) with a positive definite system matrix. By
introducing two positive diagonal parameter matrices $\Omega_1$
and $\Omega_2$ and employing the variable transformations
\[
z=\frac{1}{2}\Omega_1(|x|+x), \qquad
w=\frac{1}{2}\Omega_2(|x|-x),
\]
the original LCP is equivalently reformulated as the
modulus-based fixed-point equation
\[
(\Omega_2+A\Omega_1)x
=(\Omega_2-A\Omega_1)|x|-2q.
\]
Let $\mathcal{A}=\Omega_2+A\Omega_1$ and
\[
b^{(k)}=(\Omega_2-A\Omega_1)|x^{(k)}|-2q.
\]
The authors developed an inner--outer iteration scheme in which
the outer iteration computes an approximate solution to
$\mathcal{A}x=b^{(k)}$, while the inner iteration solves the
associated linear system approximately. To construct the inner
iteration, the coefficient matrix is decomposed into its symmetric
and skew-symmetric parts:
\[
H=\frac{\mathcal{A}+\mathcal{A}^{T}}{2},
\qquad
S=\frac{\mathcal{A}-\mathcal{A}^{T}}{2}.
\]
Thus, $\mathcal{A}=(\alpha P+H)-(\alpha P-S)$, where
$\alpha>0$ is an iteration parameter and $P$ is a symmetric
positive definite preconditioning matrix. Based on this splitting,
the inner iteration is expressed as
\[
(\alpha P+H)y^{(k,\ell+1)}
=(\alpha P-S)y^{(k,\ell)}+r^{(k)},
\]
where
\[
r^{(k)}=b^{(k)}-\mathcal{A}x^{(k)}.
\]
The preconditioned conjugate gradient (PCG) method is used to
solve the inner system inexactly. The inner iteration is terminated
when the residual
\[
p^{(k)}=r^{(k)}-(\alpha P+H)y^{(k,\ell+1)}
\]
satisfies
\[
\|p^{(k)}\|_2\leq \eta_k\|r^{(k)}\|_2,
\qquad \eta_k\to0
\quad\text{as }k\to\infty.
\]
After obtaining the approximate inner solution $y^{(k)}$, the
outer iterate is updated according to
\[
x^{(k+1)}=x^{(k)}+y^{(k)}.
\]
The outer iteration terminates when the prescribed residual-based
stopping criterion is met or the maximum number of iterations is
reached. The LCP solution is then recovered using
\[
z^{(k)}=\frac{1}{2}\Omega_1(|x^{(k)}|+x^{(k)}).
\]
The authors established sufficient convergence conditions involving
the system matrices, the preconditioning matrix, the iteration
parameter, and the inner-iteration tolerance. They also derived
sufficient conditions on $\alpha$ that characterize a convergence
region. Numerical experiments were conducted to assess the
performance of the proposed method and compare its iteration
counts and computational times with those of the PGMSOR, PMSOR,
and PMSA methods.
	\item
\textbf{Li, Wang, and Liu (2022)}~\cite{LiWangLiu2022}
proposed a relaxation general two-sweep modulus-based matrix splitting
(RGTMMS) iteration method for solving LCP(\ref{eq:LCP}).
For two positive diagonal matrices $\Omega_1$ and $\Omega_2$, let
\[
A\Omega_1=M_{\Omega_1}-N_{\Omega_1},
\]
be a matrix splitting of $A\Omega_1$. The authors first recalled the
two-sweep modulus-based matrix splitting (TMMS) iteration method and the
general two-sweep modulus-based matrix splitting (GTMMS) iteration
method, and then introduced two nonnegative relaxation parameters
$\theta$ and $\omega$ to obtain the new RGTMMS iteration method. The
proposed iteration is given by
\begin{equation*}
	\begin{aligned}
		(\Omega_2+M_{\Omega_1})x^{(k+1)}
		={}&
		N_{\Omega_1}
		\left[
		\theta x^{(k)}
		+
		(1-\theta)x^{(k-1)}
		\right]
		\\
		&+
		(\Omega_2-A\Omega_1)
		\left[
		\omega |x^{(k)}|
		+
		(1-\omega)|x^{(k-1)}|
		\right]
		-q.
	\end{aligned}
\end{equation*}
The corresponding solution of the LCP is obtained from
\begin{equation*}
	z^{(k+1)}
	=
	\Omega_1
	\left(
	|x^{(k+1)}|+x^{(k+1)}
	\right)
\end{equation*} and \begin{equation*}
w^{(k+1)}
=
\Omega_2
\left(
|x^{(k+1)}|-x^{(k+1)}
\right).
\end{equation*}
When $\theta=1$ and $\omega=1$, the RGTMMS iteration reduces to the
GTMMS iteration, while for $\theta=1$ and $\omega=0$, it becomes the
generalized modulus-based matrix splitting iteration method.

Furthermore, by choosing an AOR-type splitting
\begin{equation*}
	M_{\Omega_1}
	=
	\frac{1}{\alpha}
	\left(
	D_{\Omega_1}-\beta L_{\Omega_1}
	\right),
\end{equation*}
and
\begin{equation*}
	N_{\Omega_1}
	=
	\frac{1}{\alpha}
	\left[
	(1-\alpha)D_{\Omega_1}
	+
	(\alpha-\beta)L_{\Omega_1}
	+
	\alpha U_{\Omega_1}
	\right],
\end{equation*}
the proposed RGTMMS framework yields the relaxation general
two-sweep modulus-based accelerated overrelaxation (RGTMAOR)
iteration method:
\begin{equation*}
	\begin{aligned}
		\left(
		\Omega_2+
		\frac{1}{\alpha}
		(D_{\Omega_1}-\beta L_{\Omega_1})
		\right)x^{(k+1)}
		={}&
		\frac{1}{\alpha}
		\left[
		(1-\alpha)D_{\Omega_1}
		+
		(\alpha-\beta)L_{\Omega_1}
		+
		\alpha U_{\Omega_1}
		\right]
		\left[
		\theta x^{(k)}
		+
		(1-\theta)x^{(k-1)}
		\right]
		\\
		&+
		(\Omega_2-A\Omega_1)
		\left[
		\omega |x^{(k)}|
		+
		(1-\omega)|x^{(k-1)}|
		\right]
		-q.
	\end{aligned}
\end{equation*}
For particular choices of $(\alpha,\beta)$, the following iterative
methods are obtained from the proposed RGTMMS framework:

	 for $(\alpha,\beta)=(\alpha,\alpha)$, the
	relaxation general two-sweep modulus-based successive
	overrelaxation (RGTMSOR) iteration method is obtained;
 for $(\alpha,\beta)=(1,1)$, the
	relaxation general two-sweep modulus-based Gauss-Seidel
	(RGTMGS) iteration method is obtained;
	 for $(\alpha,\beta)=(\alpha,0)$, the
	relaxation general two-sweep modulus-based JOR
	(RGTMJOR) iteration method is obtained;
 for $(\alpha,\beta)=(1,0)$, the
	relaxation general two-sweep modulus-based Jacobi
	(RGTMJ) iteration method is obtained.

The authors established convergence conditions for the proposed RGTMMS
method when $A$ is an $H_{+}$-matrix and derived corresponding
convergence results for the RGTMAOR-type scheme. Numerical experiments
were also presented to investigate the effect of the relaxation
parameters and to compare the proposed methods in terms of iteration
numbers and computational time.

\item

\textbf{Wu and Li (2022)}~\cite{WuLi2022} proposed a class of new
modulus-based matrix splitting methods for solving the LCP (\ref{eq:LCP}).

The authors established a new equivalent fixed-point formulation
of the LCP without using a change of variables. Let
$\Omega_1$ be a positive diagonal matrix. Let
\[
a=\Omega_1 z,\qquad b=w=Az+q.
\]
Using Equation (\ref{abequivalent}), the equivalent form of the LCP is given by
\begin{equation}
	(\Omega_1+A)z
	=
	\left|(A-\Omega_1)z+q\right|-q.
	\label{eq:WuLi_fixed_point}
\end{equation}

Furthermore, consider the matrix splitting
\begin{equation*}
	A=M-N.
	\label{eq:WuLi_splitting}
\end{equation*}
Using this splitting, the authors proposed the new
modulus-based matrix splitting (NMS) iteration method:
\begin{equation*}
	(\Omega_1+M)z^{(k+1)}
	=
	Nz^{(k)}
	+
	\left|(A-\Omega_1)z^{(k)}+q\right|-q.
	\label{eq:WuLi_NMS}
\end{equation*}
Provided that $\Omega_1+M$ is nonsingular, each new iterate
is obtained by solving the corresponding linear system. The
complementary variable is recovered using
\begin{equation*}
	w^{(k+1)}=Az^{(k+1)}+q.
	\label{eq:WuLi_w}
\end{equation*}

To construct relaxation-type methods, the system matrix is
expressed as
\begin{equation*}
	A=D-L-U,
\end{equation*}
where $D=\operatorname{diag}(A)$, and $L$ and $U$ are the
strictly lower and upper triangular parts of $A$, respectively.
Choosing the AOR-type splitting
\begin{equation*}
	M=\frac{1}{\alpha}(D-\beta L),
\end{equation*}
and
\begin{equation*}
	N=\frac{1}{\alpha}
	\left[(1-\alpha)D+(\alpha-\beta)L+\alpha U\right],
\end{equation*}
yields the new modulus-based AOR (NMAOR) iteration method:
\begin{equation*}
	\begin{aligned}
		(\alpha\Omega_1+D-\beta L)z^{(k+1)}
		={}&
		\left[(1-\alpha)D+(\alpha-\beta)L+\alpha U\right]
		z^{(k)}
		\\
		&+\alpha\left(
		\left|(A-\Omega_1)z^{(k)}+q\right|-q
		\right).
	\end{aligned}
	\label{eq:WuLi_NMAOR}
\end{equation*}

For particular choices of the parameters $(\alpha,\beta)$,
the NMAOR method reduces to the new modulus-based successive
overrelaxation (NMSOR), Gauss--Seidel (NMGS), and Jacobi
(NMJ) methods. Specifically, these methods are obtained for
$(\alpha,\beta)=(\alpha,\alpha)$, $(1,1)$, and $(1,0)$,
respectively.

The authors established sufficient convergence conditions for
the proposed NMS method and its NMAOR variant. In particular,
they proved convergence results under suitable assumptions
when $A$ is an $H_+$-matrix. Numerical experiments were also
presented to demonstrate the effectiveness of the proposed
methods. The results showed that the new modulus-based
Gauss--Seidel method can be more computationally efficient
than several existing modulus-based methods.

\item 
\textbf{Huang and Cui (2023)} \cite{HuangCui2023} proposed a double-relaxation two-sweep modulus-based matrix splitting (DRTMMS) iteration method for solving LCP(\ref{eq:LCP}). The method is developed by incorporating the relaxation technique and the relaxation two-sweep strategy into the modulus-based matrix splitting (MMS) framework. In contrast to the existing MMS-type methods, the proposed scheme introduces two relaxation parameters together with an additional relaxation matrix, thereby providing a more general iterative framework. The authors showed that several existing MMS-like methods can be obtained as special cases of the proposed DRTMMS scheme.

More precisely, let $\Omega_1$ and $\Omega_2$ be positive diagonal matrices, and let
$$
A\Omega_1=M_{\Omega_1}-N_{\Omega_1}
$$
be a splitting of $A\Omega_1$. Let $P$ be a nonsingular relaxation matrix, while $\omega_1$ and $\omega_2$ denote the two relaxation parameters. Starting from initial vectors $x^{(0)}$ and $x^{(1)}$, the DRTMMS method first computes an intermediate vector $x^{(k+\frac12)}$ according to
$$
\begin{aligned}
	(\Omega_2+M_{\Omega_1})x^{(k+\frac12)}
	={}&N_{\Omega_1}
	\left[
	\omega_1x^{(k)}
	+(1-\omega_1)x^{(k-1)}
	\right]
	\\
	&+(\Omega_2-A\Omega_1)
	\left[
	(1-\omega_2)|x^{(k)}|
	+\omega_2|x^{(k-1)}|
	\right]
	-q,
\end{aligned}
$$
and then applies the relaxation step
$$
x^{(k+1)}
=(I-P)x^{(k)}
+P x^{(k+\frac12)}.
$$
The corresponding LCP solution is recovered through
$$
z^{(k+1)}
=\Omega_1
\left(
|x^{(k+1)}|+x^{(k+1)}
\right) ~~\mbox{and}~~ w^{(k+1)}
=\Omega_2
\left(
|x^{(k+1)}|-x^{(k+1)}
\right).
$$
Thus, the proposed method combines the two-sweep mechanism, controlled by $\omega_1$ and $\omega_2$, with the additional relaxation matrix $P$.

An important feature of the DRTMMS framework is that it contains several previously developed MMS-type methods as special cases. Specifically, taking $\omega_1=1$ and $\omega_2=0$ gives the RMMS method; setting $P=I$ gives the GRTMMS method; choosing $P=I$, $\Omega_1=\frac{1}{\gamma}I$, and $\Omega_2=\frac{1}{\gamma}\Omega$ gives the RTMMS method. Further, with $P=I$, $\omega_1=1$, and $\omega_2=0$, the method reduces to the GMMS scheme, while the additional choices $\Omega_1=\frac{1}{\gamma}I$ and $\Omega_2=\frac{1}{\gamma}\Omega$ recover the classical MMS method. Hence, the authors presented DRTMMS as a unified and more general framework for several existing modulus-based matrix splitting methods rather than treating these methods as newly proposed schemes.

Furthermore, by choosing the splitting
$$
M_{\Omega_1}
=\frac{1}{\alpha}
\left(D_{\Omega_1}-\beta L_{\Omega_1}\right),
$$
and
$$
N_{\Omega_1}
=\frac{1}{\alpha}
\left[
(1-\alpha)D_{\Omega_1}
+(\alpha-\beta)L_{\Omega_1}
+\alpha U_{\Omega_1}
\right],
$$
the authors derived a double-relaxation two-sweep modulus-based accelerated overrelaxation (DRTMAOR) method. With appropriate choices of $\alpha$ and $\beta$, this framework further yields the DRTM Jacobi, DRTM Gauss-Seidel, and DRTM successive overrelaxation methods.

For convergence, Huang and Cui established sufficient conditions for the DRTMMS method when the system matrix $A$ is positive definite or an $H_+$-matrix. In the positive definite case, they showed that the generated sequence converges to the exact solution of the LCP if
the convergence condition for the DRTMMS iteration is given by
\[
\rho(F+G)<1,
\]
where \(F\) and \(G\) are the nonnegative matrix bounds associated with the two-step error propagation of the DRTMMS iteration, defined as
\[
\begin{aligned}
	F ={}& \left|(I-P)+\omega_1 P(\Omega_2+M\Omega_1)^{-1}N\Omega_1\right|+|1-\omega_2|
	\left|P(\Omega_2+M\Omega_1)^{-1}(\Omega_2-M\Omega_1)\right|\\
	&+|1-\omega_2|
	\left|P(\Omega_2+M\Omega_1)^{-1}N\Omega_1\right|,
\end{aligned}
\]
and
\[
\begin{aligned}
	G ={}& |1-\omega_1|
	\left|P(\Omega_2+M\Omega_1)^{-1}N\Omega_1\right|+\omega_2
	\left|P(\Omega_2+M\Omega_1)^{-1}(\Omega_2-M\Omega_1)\right|\\
	&+\omega_2
	\left|P(\Omega_2+M\Omega_1)^{-1}N\Omega_1\right|.
\end{aligned}
\]
Here, \(\rho(F+G)\) denotes the spectral radius of \(F+G\). They also derived norm-based sufficient convergence conditions and investigated practical parameter-selection strategies for the relaxation parameters.

For practical parameter selection, the authors considered $P=\theta I$ and $\omega_2=0$ and obtained a quasi-optimal choice of $\theta$ by minimizing an upper bound of the spectral radius. In particular,
$$
\theta^*
=
\min\left\{
1,\,
\frac{
	\left\|(\Omega_2+M_{\Omega_1})^{-1}
	(\Omega_2-M_{\Omega_1})\right\|
	+
	2\left\|(\Omega_2+M_{\Omega_1})^{-1}N_{\Omega_1}\right\|
}{2}
\right\}.
$$
This provides a practical way of selecting the relaxation parameter without directly minimizing the spectral radius, which is generally difficult to determine.

Finally, numerical experiments reported by Huang and Cui demonstrated that, with suitable choices of the relaxation parameters and parameter matrices, the DRTMMS method can achieve better numerical performance than several existing MMS-like methods and can outperform projected splitting methods for some test problems.
\item 
\textbf{Kumar, Deepmala, and Das (2023)}
~\cite{KumarDeepmalaDas2023opsearch} proposed a family of new modified
modulus-based matrix splitting methods for solving the LCP
(\ref{eq:LCP}). Using Equation (\ref{abequivalent}),
they introduced two positive diagonal matrices $\Omega_1$ and
$\Omega_2$ and defined
\[
a=\Omega_1(Az+q),\qquad b=\Omega_2z.
\]
Consequently, the LCP is equivalent to the fixed-point equation
\[
(\Omega_1A+\Omega_2)z
=
\left|(\Omega_1A-\Omega_2)z+\Omega_1q\right|
-\Omega_1q.
\]
Using the matrix splitting
\[
A=(M+\Omega)-(N+\Omega),
\]
where $M\Omega_1+\Omega\Omega_1+\Omega_2$ is nonsingular, the
resulting new modified modulus-based matrix splitting iteration is
\[
\begin{aligned}
	z^{(k+1)}
	={}&(M\Omega_1+\Omega\Omega_1+\Omega_2)^{-1}
	\Big[
	(N\Omega_1+\Omega\Omega_1)z^{(k)}\\
	&\quad+
	\left|(A\Omega_1-\Omega_2)z^{(k)}+\Omega_1q\right|
	-\Omega_1q
	\Big].
\end{aligned}
\]
This general framework yields several proposed methods as special cases.
For $M=A$, $N=0$, $\Omega_1=I$, and $\Omega_2=I$, the new modified
modulus-based iterative method (NMMOD) is obtained as
\[
(A+\Omega+I)z^{(k+1)}
=
\Omega z^{(k)}
+\left|(A-I)z^{(k)}+q\right|-q.
\]
For $M=A$, $N=0$, $\Omega_1=I$, and $\Omega_2=\alpha I$, the new modified,
modified modulus-based iterative method (NMMMOD) is obtained as
\[
(A+\Omega+\alpha I)z^{(k+1)}
=
\Omega z^{(k)}
+\left|(A-\alpha I)z^{(k)}+q\right|-q.
\]
Furthermore, by taking
\[
M=\frac{1}{\alpha}(D-\beta L),~
N=\frac{1}{\alpha}
\left[(1-\alpha)D+(\alpha-\beta)L+\alpha U\right],
\]
where $A=D-L-U$, the new modified modulus-based accelerated
overrelaxation (NMMAOR) method is obtained:
\[
\begin{aligned}
	(D\Omega_1-\beta L\Omega_1
	+\alpha\Omega\Omega_1+\alpha\Omega_2)z^{(k+1)}
	&=
	\left[(1-\alpha)D\Omega_1
	+(\alpha-\beta)L\Omega_1
	+\alpha U\Omega_1
	+\alpha\Omega\Omega_1\right]z^{(k)}
	\\
	&~\quad
	+\alpha\left(
	\left|(A\Omega_1-\Omega_2)z^{(k)}
	+\Omega_1q\right|-\Omega_1q
	\right).
\end{aligned}
\]
The NMMAOR method contains the new modified modulus-based successive
overrelaxation (NMMSOR), new modified modulus-based Gauss-Seidel
(NMGSD), and new modified modulus-based Jacobi (NMMJ) methods as
special cases for appropriate choices of the parameters $\alpha$ and
$\beta$. Thus, the paper develops a unified family consisting of NMMOD,
NMMMOD, NMMAOR, NMMSOR, NMGSD, and NMMJ methods.

For a $P$-matrix, the proposed general iteration converges to the unique
solution for any nonnegative initial vector if
\[
\rho\!\left(
\left|(M\Omega_1+\Omega\Omega_1+\Omega_2)^{-1}\right|
\left(
|N\Omega_1+\Omega\Omega_1|
+
|A\Omega_1-\Omega_2|
\right)
\right)<1.
\]
The authors also established sufficient convergence conditions for
$H_{+}$-matrices, including conditions for the general new modified
modulus-based method and the NMMAOR method. In particular, for an
$H_{+}$-matrix, convergence of the general method is guaranteed when
$
\Omega_2\geq D\Omega_1,
$
under the assumed $H$-splitting conditions. For the NMMAOR method,
additional conditions involving $\alpha$, $\beta$, and
\[
\rho\!\left(D_{\Omega_1}^{-1}|B\Omega_1|\right)<1,
~ B=L+U,
\]
were established. Numerical experiments compared NMMOD, NMMMOD, and
NMGSD with the corresponding MOD, MMOD, and GSD methods and showed that
the proposed methods generally require fewer iterations and achieve
higher computational efficiency. 
\item 
\textbf{Kumar, Deepmala and Das (2024)}
~\cite{KumarDeepmalaDas2024} proposed a class of new accelerated
modulus-based iteration methods for solving LCP (\ref{eq:LCP}).
The proposed methods are developed by reformulating the problem into an
equivalent fixed-point equation. Using the transformations
\[
z=\frac{1}{\gamma}(|x|+x),
\qquad
w=\Omega(|x|-x),
\]
the LCP can be expressed in an equivalent form in terms of the variable
\(x\). In particular, they considered the matrix splitting
\[
A=(M+I-L)-(N+I-L),
\]
where \(I\) is the identity matrix and \(L\) is the strictly lower
triangular part of \(A\). The resulting fixed-point formulation is
\[
\begin{aligned}
	x^{(k+1)}
	={}&(M+\Omega+I-L)^{-1}
	\big[(N+I-L)x^{(k)}
	+(\Omega-A)|x^{(k)}|-\gamma q\big].
\end{aligned}
\]
Here, \(\Omega\) is a positive diagonal matrix and \(\gamma>0\).
Based on this general formulation, several accelerated modulus-based
methods are obtained by selecting different matrix splittings and
parameters. For \(M=A\), \(N=0\), \(\Omega=I\), and \(\gamma=1\), the new
accelerated modulus-based iteration method is
\[
\begin{aligned}
	x^{(k+1)}
	={}&(A+2I-L)^{-1}
	\big[(I-L)x^{(k)}+(I-A)|x^{(k)}|-q\big].
\end{aligned}
\]
For \(M=A\), \(N=0\), \(\Omega=\alpha I\), and \(\gamma=1\), the new
accelerated modified modulus-based iteration method is
\[
\begin{aligned}
	x^{(k+1)}
	={}&(A+(\alpha+1)I-L)^{-1}
	\big[(I-L)x^{(k)}+(\alpha I-A)|x^{(k)}|-q\big].
\end{aligned}
\]
For \(M=D\), \(N=L+U\), and \(\gamma=2\), the new accelerated modulus-based
Jacobi iteration method is
\[
\begin{aligned}
	x^{(k+1)}
	={}&(D+\Omega+I-L)^{-1}
	\big[(U+I)x^{(k)}+(\Omega-A)|x^{(k)}|-2q\big].
\end{aligned}
\]
For \(M=D-L\), \(N=U\), and \(\gamma=2\), the new accelerated
modulus-based Gauss-Seidel (NAMGS) method is
\[
\begin{aligned}
	x^{(k+1)}
	={}&(D-2L+\Omega+I)^{-1}
	\big[(U+I-L)x^{(k)}+(\Omega-A)|x^{(k)}|-2q\big].
\end{aligned}
\]
For \(M=\frac{1}{\alpha}D-L\),
\(N=(\frac{1}{\alpha}-1)D+U\), and \(\gamma=2\), the new accelerated
modulus-based successive over-relaxation (NAMSOR) method is
\[
\begin{aligned}
	x^{(k+1)}
	={}&(D-2\alpha L+\alpha\Omega+\alpha I)^{-1}
	\big[((1-\alpha)D+\alpha U+\alpha I-L)x^{(k)}\\
	&\qquad+(\alpha\Omega-\alpha A)|x^{(k)}|
	-2\alpha q\big].
\end{aligned}
\]
Furthermore, by choosing
\(M=\frac{1}{\alpha}(D-\beta L)\) and
$
N=\frac{1}{\alpha}
[(1-\alpha)D+(\alpha-\beta)L+\alpha U],
$
the authors obtained the new accelerated modulus-based accelerated
over-relaxation (NAMAOR) method:
\[
\begin{aligned}
	x^{(k+1)}
	={}&(D-(\beta+\alpha)L
	+\alpha\Omega+\alpha I)^{-1}\\
	&\quad\big[((1-\alpha)D
	+(2\alpha-\beta)L+\alpha U+\alpha I)x^{(k)}\\
	&\qquad+(\alpha\Omega-\alpha A)|x^{(k)}|
	-2\alpha q\big].
\end{aligned}
\]
The NAMAOR scheme further reduces to the NAMSOR, NAMGS, and accelerated
Jacobi methods for suitable choices of the parameters. The authors also
established convergence conditions when \(A\) is a \(P\)-matrix or an
\(H_+\)-matrix and demonstrated through numerical experiments that the
proposed NAMGS and NAMSOR methods require fewer iterations and less CPU
time than the corresponding MGS and MSOR methods.

\item \textbf{Fang (2024)}~\cite{Fang2024TSMAOR} proposed a two-step
simplified modulus-based matrix splitting iteration in the standard LCP (\ref{eq:LCP}).
Consider two positive diagonal matrices $\Omega_1$ and $\Omega_2$ and
two matrix splittings
\[
A=M_1-N_1=M_2-N_2.
\]
Let
\[
a=\Omega_1 z,\qquad b=\Omega_2 w,
\]
where \(w=Az+q\). Using the identity \(a+b=|a-b|\), we obtain
\[
(\Omega_1+\Omega_2 A)z
=
|(\Omega_1-\Omega_2 A)z-\Omega_2 q|
-\Omega_2 q.
\]
Based on this formulation, the two-step simplified modulus-based matrix
splitting (TSMMS) method is defined as
\[
\begin{aligned}
	(\Omega_1+M_1)z^{(k+1)}
	={}&N_1z^{(k+\frac12)}
	+\left|(A-\Omega_1)z^{(k+\frac12)}+q\right|-q,\\
	(\Omega_2+M_2)z^{(k+\frac12)}
	={}&N_2z^{(k)}
	+\left|(A-\Omega_2)z^{(k)}+q\right|-q.
\end{aligned}
\]
The authors further considered the accelerated over-relaxation (AOR)
splittings given by
\[
M_1=\frac{1}{\omega_1}(D-\gamma_1L),
\qquad
N_1=\frac{1}{\omega_1}
\left[(1-\omega_1)D+(\omega_1-\gamma_1)L+\omega_1U\right],
\]
and
\[
M_2=\frac{1}{\omega_2}(D-\gamma_2U),
\qquad
N_2=\frac{1}{\omega_2}
\left[(1-\omega_2)D+(\omega_2-\gamma_2)U+\omega_2L\right],
\]
where \(A=D-L-U
\), with \(D\), \(L\), and \(U\) denoting the diagonal, strictly
lower triangular, and strictly upper triangular parts of \(A\),
respectively.
Substituting these splittings into the TSMMS framework gives the
two-step simplified modulus-based accelerated over-relaxation (TSMAOR)
method:
\begin{equation*}
	\begin{aligned}
		(\omega_1\Omega_1+D-\gamma_1L)z^{(k+1)}
		&=\left[(1-\omega_1)D+(\omega_1-\gamma_1)L+\omega_1U\right]
		z^{(k+\frac12)}\\
		&\quad+\omega_1\left|(A-\Omega_1)z^{(k+\frac12)}+q\right|
		-\omega_1q,\\[6pt]
		(\omega_2\Omega_2+D-\gamma_2U)z^{(k+\frac12)}
		&=\left[(1-\omega_2)D+(\omega_2-\gamma_2)U+\omega_2L\right]
		z^{(k)}\\
		&\quad+\omega_2\left|(A-\Omega_2)z^{(k)}+q\right|
		-\omega_2q.
	\end{aligned}
\end{equation*}
The proposed TSMAOR method contains several two-step simplified
modulus-based schemes as special cases. In particular, choosing
$\omega_i=\gamma_i$ gives the two-step simplified modulus-based
successive over-relaxation (TSMSOR) method, while
$\omega_i=\gamma_i=1$ gives the two-step simplified modulus-based
Gauss-Seidel (TSMGS) method. Furthermore,
$\omega_i=1$ and $\gamma_i=0$ gives the two-step simplified
modulus-based Jacobi (TSMJ) method, for $i=1,2$.
The paper establishes convergence results for the proposed schemes
under appropriate assumptions, including convergence results for
$H_+$-matrices, and numerical experiments are provided to demonstrate
their effectiveness.
	\item 
	\textbf{Bashirizadeh and Hajarian (2024)}
	~\cite{BashirizadehHajarian2024} proposed two classes of modulus-based
	matrix splitting iteration methods for solving LCP(\ref{eq:LCP}). Let
	\[
	A=M_1-N_1=M_2-N_2
	\]
	be two splittings of $A$, where $\Omega$ is a positive diagonal
	matrix. Starting from two initial vectors
	$x^{(0)},y^{(0)}\in\mathbb{R}^n$, the first proposed method, referred
	to as modulus-based matrix splitting iteration method I, is given by
	\begin{equation*}
		\begin{aligned}
			(\Omega+M_1)y^{(k+1)}
			&=
			N_1y^{(k)}
			+(\Omega-A)|x^{(k)}|
			-\gamma q,
			\\
			(\Omega+M_2)x^{(k+1)}
			&=
			N_2y^{(k+1)}
			+(\Omega-A)|x^{(k)}|
			-\gamma q.
		\end{aligned}
	\end{equation*}
	The LCP solution is recovered through
	\[
	z^{(k+1)}
	=
	\frac{1}{\gamma}
	\left(|x^{(k+1)}|+x^{(k+1)}\right).
	\]
	The second proposed scheme, referred to as modulus-based matrix
	splitting iteration method II, uses the updated vector
	$y^{(k+1)}$ in the modulus term of the second step:
	\begin{equation*}
		\begin{aligned}
			(\Omega+M_1)y^{(k+1)}
			&=
			N_1y^{(k)}
			+(\Omega-A)|x^{(k)}|
			-\gamma q,
			\\
			(\Omega+M_2)x^{(k+1)}
			&=
			N_2y^{(k+1)}
			+(\Omega-A)|y^{(k+1)}|
			-\gamma q.
		\end{aligned}
	\end{equation*}
	Again, the LCP solution is obtained from
	\[
	z^{(k+1)}
	=
	\frac{1}{\gamma}
	\left(|x^{(k+1)}|+x^{(k+1)}\right)
	\] and 	$w^{(k+1)}
	=
	\Omega
	\left(
	|x^{(k+1)}|-x^{(k+1)}
	\right).$
	The authors showed that the two proposed schemes provide general
	frameworks from which different modulus-based iterative methods can be
	obtained by appropriate choices of the matrix splittings. In the
	numerical experiments, AOR-type splittings lead to the corresponding
	MSOR I and MSOR II schemes, which are compared with the conventional
	MSOR and TSMSOR methods. The convergence of the proposed methods is
	established when $A$ is an $H_{+}$-matrix, under the stated assumptions
	on the splittings and the positive diagonal matrix $\Omega$. Numerical
	experiments demonstrate that the proposed MSOR I and MSOR II schemes
	can outperform the existing MSOR and TSMSOR methods in terms of
	iteration number and CPU time~\cite{BashirizadehHajarian2024}.
	\item \textbf{Guo, Zheng, Lu, and Zhang (2024)}
	~\cite{GuoZhengLuZhang2024} investigated the convergence of two-stage
	multisplitting iteration methods for solving LCP(\ref{eq:LCP}). Unlike a paper proposing a new iteration scheme, this work
	focuses on obtaining improved convergence results for the existing
	modulus-based synchronous two-stage multisplitting (MSTM) method and
	the two-stage multisplitting modulus-based matrix splitting (TMMMS)
	method. The authors begin with the modulus-based equivalent equation
	\[
	(\Omega+M)x
	=
	Nx+(\Omega-A)|x|-\gamma q,
	\]
	where $\Omega$ is a positive diagonal matrix and
	$A=M-N$. For a family of two-stage multisplittings
	\[
	A=M_s-N_s,~
	M_s=F_s-G_s,~ s=1,2,\ldots,\ell,
	\]
	with nonnegative diagonal weighting matrices $E_s$ satisfying
	\[
	\sum_{s=1}^{\ell}E_s=I,
	\]
	the MSTM method reviewed in the paper computes local iterates through
	the two-stage splitting and combines them using the weighting matrices.
	In particular, with positive integers $v_s$ and $v_s^{(m)}$, the
	MSTM iteration is based on the local computations
	\[
	(\Omega+F_s)y_s^{(k+1)}
	=
	G_sy_s^{(k)}
	+(\Omega-A)|x^{(k)}|-\gamma q,
	\]
	followed by the multisplitting combination
	\[
	x^{(k+1)}
	=
	\sum_{s=1}^{\ell}E_sy_s^{(k+1)},
	\]
	and the LCP solution is recovered as
	\[
	z^{(k+1)}
	=
	\frac{1}{\gamma}
	\left(|x^{(k+1)}|+x^{(k+1)}\right)
	\]
	and 	$w^{(k+1)}
	=
	\Omega
	\left(
	|x^{(k+1)}|-x^{(k+1)}
	\right).$
	The paper also reviews the TMMMS method, in which the inner iteration
	for each worker is taken to be the modulus-based matrix splitting (MMS)
	iteration. Specifically, starting from $x^{(0)}$ and
	\[
	z^{(0)}
	=
	\frac{1}{\gamma}
	\left(x^{(0)}+|x^{(0)}|\right),
	\]
	the local inner iterations are performed according to
	\[
	(\Omega+F_s)x^{(k,s,j+1)}
	=
	G_sx^{(m,s,j)}
	+(\Omega-M_s)|x^{(k,s,j)}|
	-\gamma\left(q-N_sz^{(m)}\right),
	\]
	for
	$
	j=0,1,\ldots,v_s^{(k)}-1,
	$
	with
	$
	x^{(k,s,0)}=x^{(k)}.
	$
	The local results are then combined as
$
	x^{(k+1)}
	=
	\sum_{s=1}^{\ell}
	E_sx^{(k,s,v_s^{(k)})},
$
	and
	$
	z^{(k+1)}
	=
	\frac{1}{\gamma}
	\left(x^{(k+1)}+|x^{(k+1)}|\right).
	$
	The main contribution of Guo, Zheng, Lu, and Zhang ~\cite{GuoZhengLuZhang2024} is the derivation of new convergence conditions for the existing multisplitting two-stage method (MSTM) and two-stage multisplitting modulus-based matrix splitting method (TMMMS), rather than the proposal of these methods themselves. Assuming that the coefficient matrix \(A\) is an \(H_+\)-matrix, they established sufficient convergence conditions that enlarge the admissible domains of the parameter matrices. In particular, their analysis allows the positive diagonal matrix \(\Omega\) to satisfy conditions weaker than the commonly used \(\Omega\geq D\), thereby extending the applicability of both methods. The theoretical results were validated through two numerical examples and parallel computations using OpenACC. The numerical experiments demonstrated the effectiveness of the enlarged convergence domains and showed that both MSTM and TMMMS can converge faster than the modulus-based splitting method (MSM). Moreover, TMMMS incorporates the modulus-based matrix splitting (MMS) iteration as an inner iteration within the two-stage multisplitting framework to reduce computational time, highlighting the potential of these methods for efficient parallel solution of LCPs.
\item
\textbf{Liu, Wang, and Huang (2025)}~\cite{LiuWangHuang2025} proposed improved
modulus-based and projected iterative methods for solving LCP(\ref{eq:LCP})  involving $Z$-matrices whose diagonal
entries may be positive, zero, or negative. For the LCP 
they considered the modulus-based reformulation
\[
z=\frac{1}{\gamma}\left(|x|+x\right),~
w=\frac{\Omega}{\gamma}\left(|x|-x\right),
\]
where $\gamma>0$ and $\Omega$ is a positive diagonal matrix. By using
the splitting
$
A=M-N,
$
the corresponding implicit fixed-point equation can be written as
\[
(M+\Omega)x
=
Nx+(\Omega-A)|x|-\gamma q.
\]

For a $Z$-matrix $A$, Liu, Wang, and Huang considered the special
choice
\[
\Omega=\delta I,
~
\delta=\max_{1\leq i\leq n}|a_{ii}|.
\]
Based on the corresponding matrix splitting, they developed an
improved modulus-based Gauss-Seidel method. The resulting iteration
takes the form
\[
x^{(k+1)}
=
\frac{1}{2}\Omega^{-1}
\left[
L\left(x^{(k+1)}+|x^{(k+1)}|\right)
+
(U+\widehat{\Omega})
\left(x^{(k)}+|x^{(k)}|\right)
-\gamma q
\right].
\]
They further established the relationship between the improved
modulus-based Gauss-Seidel method and the corresponding improved
projected Gauss-Seidel method. In particular, the latter can be
expressed componentwise as
\[
z_i^{(k+1)}
=
\max\left\{
0,\,
z_i^{(k)}
-\frac{1}{\delta}
\left(
\sum_{j=1}^{i-1}a_{ij}z_j^{(k+1)}
+
\sum_{j=i}^{n}a_{ij}z_j^{(k)}
+
q_i
\right)
\right\},
\]
where
$
\Omega=\delta I.
$

Similarly, Liu, Wang, and Huang proposed an improved modulus-based
Jacobi method by taking
$
M=\Omega,~
N=\Omega-A.
$
The resulting iteration is
\[
x^{(k+1)}
=
\frac{1}{2}\Omega^{-1}
\left[
N\left(x^{(k)}+|x^{(k)}|\right)-\gamma q
\right].
\]
The corresponding improved projected Jacobi method is obtained by
choosing
\[
\Omega=\delta I,
~
\delta=\max_{1\leq i\leq n}|a_{ii}|.
\]
Liu, Wang, and Huang established the equivalence between the proposed
improved modulus-based and projected iterative methods and provided
convergence results for the resulting algorithms. These methods
therefore extend modulus-based and projected Gauss-Seidel and Jacobi
approaches to LCPs involving $Z$-matrices with diagonal entries that
are not necessarily positive.

	\item
 \textbf{Gander, Liao, and Lu (2026)}
~\cite{GanderLiaoLu2026} proposed a new Schwarz modulus-based matrix
splitting method for the LCPs arising from
American option pricing and further accelerated the resulting iteration
using Modified Polynomial Extrapolation (MPE). After transforming and
discretizing the Black-Scholes complementarity problem, the authors
obtain a sequence of standard LCP(\ref{eq:LCP}).
For the modulus-based formulation, they consider a matrix splitting
$
A=M-N
$
and a positive diagonal matrix $\Omega$, and reformulate the LCP as the
implicit fixed-point problem
\[
(M+\Omega)x
=
Nx+(\Omega-A)|x|-\gamma q,
\]
where $\gamma>0$. Once the fixed-point problem is solved, the LCP solution
is recovered from
$$
z=\frac{|x|+x}{\gamma}
~~ \mbox{and} 	~~w
=
\Omega(|x|-x)
.$$

The corresponding stationary modulus-based iteration is therefore
\[
(M+\Omega)x^{(k+1)}
=
Nx^{(k)}
+(\Omega-A)|x^{(k)}|
-\gamma q.
\]
The main new method proposed by the authors is an
\emph{alternating Schwarz modulus-based matrix splitting method}.
The spatial domain is decomposed into overlapping subdomains with a
minimal overlap of one mesh size, leading to a block Gauss-Seidel type
splitting. Although an arbitrary number of overlapping subdomains can be
considered, the authors focus on a two-subdomain decomposition
because it gives the best performance in their numerical experiments.
For the two-subdomain case, the splitting is defined by
\[
A=
\begin{pmatrix}
	A_{11}&A_{12}\\
	A_{21}&A_{22}
\end{pmatrix},
~
M=
\begin{pmatrix}
	A_{11}&0\\
	A_{21}&A_{22}
\end{pmatrix},
~
N=
\begin{pmatrix}
	0&-A_{12}\\
	0&0
\end{pmatrix}.
\]
This Schwarz splitting gives the new block Gauss-Seidel (BGS)
modulus-based method. For comparison, the authors also consider the
classical pointwise Gauss-Seidel splitting
\[
M=D+L,~ N=-U,
\]
where $D$, $L$, and $U$ denote the diagonal, strictly lower triangular,
and strictly upper triangular parts of $A$, respectively.

To further accelerate the convergence of the modulus-based iteration,
the authors introduce \emph{Modified Polynomial Extrapolation (MPE)}.
For the sequence
\[
x^{(0)},x^{(1)},\ldots,x^{(k+1)}
\]
generated by the underlying modulus-based iteration, the difference
vectors are defined by
\[
u^{(i)}=x^{(i+1)}-x^{(i)}.
\]
The MPE approximation is constructed as a linear combination
\[
s^{(k)}=\sum_{i=0}^{k}\gamma_i x^{(i)},
~
\sum_{i=0}^{k}\gamma_i=1,
\]
where the coefficients $\gamma_i$ are determined from the corresponding
difference vectors through a least-squares problem. In this way, MPE
accelerates the convergence of the nonlinear vector sequence generated
by the modulus-based splitting iteration.

The numerical study considers both the classical pointwise
Gauss-Seidel splitting and the proposed Schwarz/block Gauss-Seidel
splitting, together with their MPE-accelerated versions. In particular,
the methods compared are GS, MPE-GS, MPECycle-GS, BGS, MPE-BGS, and
MPECycle-BGS. The MPECycle variants apply MPE only after a prescribed
number of underlying GS or BGS iterations, thereby reducing the
computational cost associated with repeatedly solving the least-squares
problem required by MPE. The numerical experiments show that the
proposed Schwarz/BGS splitting requires fewer iterations than the
classical pointwise GS method, while MPE and MPECycle acceleration can
further reduce the iteration count and computational time.
\end{enumerate}
\paragraph{}
The development of modulus-based methods for LCPs has followed a progressive path from the early modulus algorithms to increasingly general, flexible, and computationally efficient iterative schemes. The principal developments can be summarized in Table~\ref{tab:modulus_chronology}.
\newpage
\begin{table*}[ht]
	\centering
	\caption{Chronological development of modulus-based iterative methods for LCPs.}
	\label{tab:modulus_chronology}
	\begin{tabular}{ccc}
		\hline
		\textbf{Year} & \textbf{Authors} & \textbf{Main contribution} \\
		\hline
		1980 & \textbf{Van Bokhoven} & Classical modulus algorithm for LCPs. \\
		
		1986 & \textbf{Kappel and Watson} & Block modulus algorithms. \\
		
		2009 & \textbf{Hadjidimos and Tzoumas} & Modulus and extrapolated methods. \\
		
		2009 & \textbf{Dong and Jiang} & Modified modulus  method. \\
		
		2010 & \textbf{ Bai} & Modulus-based matrix splitting methods. \\
		
		2012 & \textbf{Hadjidimos et al.} & Scaled extrapolated block-modulus methods. \\
		2013 & \textbf{
			Li} & General modulus-based method. \\
		2013 & \textbf{Zhang and Ren} & Convergence analysis. \\
		
		2013 & \textbf{Zheng and Yin} & Accelerated MMS methods. \\
		
		2014 & \textbf{Cvetkovi\'c et al.} & MAOR-type modulus methods. \\
		
		2014 & \textbf{Xu and Liu} & Modified general  MMS methods. \\
		
		2019 & \textbf{Xu and Zhu} & Corrigendum of modified general  MMS methods. \\
		2014 & \textbf{Zheng and Yin} & Convergence of accelerated MMS method. \\
		
			2015 & \textbf{Fang and Wei} & General modulus Jacobi method. \\
			
		2015 & \textbf{Zheng and Li} & Modulus-based non smooth newton method. \\
		
		2016 & \textbf{Liu et al.} & Accelerated MMS methods. \\
		
		2016 & \textbf{Wu and Li} & Two-sweep MMS method. \\
			2017 & \textbf{Zheng et al.} &  Relaxation MMS method. \\
		2018 & \textbf{Wen et al.} & Relaxed MMS methods. \\
		
		2018 & \textbf{Wu et al.} & Preconditioned MMS methods. \\
		
		2019 & \textbf{Fang and Zhu} & Modulus-based double splitting method. \\
		
		2019 & \textbf{Ren et al.} &  Preconditioned general two-sweep MMS method. \\
		
		2020 & \textbf{Wu and Li} & Preconditioned  MMS method. \\
		
		2022 & \textbf{Wu et al.} & Inexact and generalized two-sweep methods. \\
		
		2022 & \textbf{Li et al.} & Relaxation general two-sweep methods. \\
		
			2022 & \textbf{Wu and Li} & New MMS method. \\
		
			2023 & \textbf{Huang and Cui} & Double relaxation two sweep MMS  method. \\
			
			2023 & \textbf{Kumar et al.} & New modified MMS method. \\
			
				2024 & \textbf{Kumar et al.} & New accelerated MMS  method. \\
				
					2024 & \textbf{Fang} & Simplified TSMAOR method. \\
					
					2024 & \textbf{Bashirzadeh and Hajarian } & MMS  method. \\
					
						2024 & \textbf{Guo et al.} & Convergence of two stage multi-splitting  method. \\
		
		2025 & \textbf{Liu et al.} & Improved Jacobi and Gauss-Seidel methods. \\
		
		2026 & \textbf{Gander et al.} & Schwarz and accelerated methods. \\
		\hline
	\end{tabular}
\end{table*}
Taken together, these developments indicate a gradual shift from basic modulus formulations toward more flexible and computationally efficient iterative frameworks. The literature has progressively incorporated acceleration, relaxation, preconditioning, advanced splitting strategies, and inexact techniques to improve convergence and computational performance. Despite these advances, the effectiveness of individual methods remains dependent on the matrix structure and algorithmic parameters, leaving scope for further development of adaptive, robust, and scalable modulus-based methods for large-scale LCPs.

\section{A Critical Comparative Assessment of Projection-Based and Modulus-Based Methods for LCPs}
\label{sec:comparative_assessment}
The numerical solution of LCP
has motivated the development of several iterative frameworks, among
which projection-based and modulus-based methods constitute two
important and complementary approaches. Although both frameworks
seek to satisfy the same complementarity requirements, they differ
substantially in the manner in which these conditions are incorporated
into an iterative process. Projection-based methods exploit the
geometric structure of the feasible region through projection
operators, whereas modulus-based methods employ an algebraic
transformation of the complementarity conditions that facilitates the
construction of matrix-splitting and relaxation schemes. These
differences influence not only the convergence analysis but also the
computational cost, parameter dependence, implementation complexity,
and suitability of the resulting algorithms for different classes
of problems.

A comprehensive comparison must therefore extend beyond the formal
structure of the underlying reformulations. In particular, the
performance of an iterative method depends on the interaction between
the coefficient matrix, the selected algorithmic parameters, the
linear algebra operations performed at each iteration, and the
accuracy required from the computed solution. A method that is
attractive from a theoretical perspective may not necessarily be
the most efficient in a practical computation, while a comparatively
simple scheme may perform well when the problem possesses a favorable
matrix structure. Accordingly, this section examines the two
frameworks from several interconnected perspectives, with particular
emphasis on convergence mechanisms, computational efficiency,
structural suitability, and the methodological considerations that
should guide their application to LCPs.

\subsection{Convergence Mechanisms and Theoretical Considerations}
Convergence is one of the principal criteria for assessing the
reliability of an iterative method for an LCP. However, convergence
cannot generally be attributed to the reformulation alone; rather,
it is determined by the properties of the coefficient matrix, the
equivalent or transformed problem, and the iterative operator
generated by the chosen algorithm. In projection-based methods,
convergence analysis commonly relies on properties of the projection
operator, the associated fixed-point mapping, and suitable conditions
on the step-size or relaxation parameters. When the iteration can
be expressed as a fixed-point process, establishing that the
corresponding operator is contractive, or satisfies an appropriate
weaker convergence condition, provides a natural route to proving
convergence. Nevertheless, the validity of such arguments depends
on the assumptions imposed on the problem and should not be
interpreted as an unconditional guarantee for arbitrary LCPs.

For modulus-based methods, convergence is frequently investigated
through the transformed nonlinear system and the matrix splitting
employed to construct the iterative scheme. The decomposition of
the coefficient matrix into suitable components determines the
resulting iteration matrix and, consequently, influences the
convergence conditions. Properties such as positive definiteness,
appropriate matrix regularity, and membership in particular matrix
classes may provide useful sufficient conditions, depending on the
specific method under consideration. The transformation and
splitting parameters can further affect the spectral characteristics
of the iteration. Thus, the theoretical analysis of modulus-based
methods requires careful consideration of both the original matrix
structure and the algebraic properties of the transformed system.

An important distinction should be made between the existence of
a convergent iteration and the speed at which convergence occurs.
A theoretical result may establish convergence under prescribed
assumptions without providing a sharp estimate of the number of
iterations required to attain a specified accuracy. Similarly,
a method with a favorable asymptotic convergence factor may not
necessarily minimize the total computational time, particularly
when its individual iterations are expensive. Consequently, a
meaningful theoretical comparison should distinguish convergence
guarantees, convergence rates, admissible parameter ranges, and
the strength of the assumptions required by each analysis. This
distinction is particularly relevant when comparing methods
developed for different classes of coefficient matrices.

The available literature illustrates the diversity of these
convergence mechanisms. Bai \cite{Bai2010} developed and analyzed
modulus-based matrix-splitting methods for complementarity problems,
demonstrating the importance of suitable splitting strategies and
parameter choices. More recently, projection-based matrix-splitting
approaches have also been investigated for structured classes of
problems, including LCPs associated with $P$- and $H_{+}$-matrices
\cite{DasDeepmalaKumar2025}. These developments suggest that
convergence should be assessed at the level of the specific
iterative algorithm rather than inferred solely from its
classification as projection-based or modulus-based.

\subsection{Computational Complexity and Numerical Efficiency}
The computational efficiency of an iterative LCP solver is governed
by the cumulative cost of its iterations rather than by the cost
of a single update in isolation. Projection-based methods often
combine matrix-vector multiplication, relaxation, and a projection
onto the relevant feasible set. When the projection can be evaluated
componentwise or through an inexpensive closed-form expression, the
resulting iteration may have a relatively low computational overhead.
This feature can be advantageous for large-scale problems in which
the coefficient matrix is sparse and the dominant operations can
be implemented efficiently. However, the actual cost depends on
the specific projection, the representation of the feasible set,
and the additional operations required by the chosen algorithm.

Modulus-based methods introduce a different computational structure.
Depending on the formulation, an iteration may involve a modulus
transformation, the evaluation of split matrix components, and
the solution of one or more auxiliary linear systems. If the
splitting produces a matrix with a structure that permits efficient
solution, such as a triangular or otherwise easily invertible
component, the additional operations may remain computationally
manageable. Conversely, an unsuitable splitting can increase
computational overhead or lead to expensive linear algebra
operations. Therefore, the computational cost of a modulus-based
method cannot be characterized solely by the presence of a modulus
mapping; the structure and implementation of the associated
linear systems are equally important.

A useful assessment must distinguish between per-iteration cost
and total time-to-solution. Let $k$ denote the number of iterations
required to satisfy a prescribed stopping criterion, and let
$T_{\mathrm{iter}}$ denote the representative cost of one iteration.
The total computational time can be approximated by
\[
T_{\mathrm{total}}
\approx k\,T_{\mathrm{iter}},
\]
provided that the variation in iteration cost and initialization
overhead is sufficiently small. This simple relation highlights
why a method with more expensive iterations can still be preferable
if it requires substantially fewer iterations. Conversely, a
low-cost iteration offers limited practical benefit when convergence
is slow or when the parameter selection process requires repeated
computations.

For a reliable numerical comparison, computational time should
therefore be reported alongside iteration counts and an explicitly
defined residual or complementarity measure. Memory consumption,
sparsity preservation, preprocessing costs, and the expense of
parameter tuning should also be considered when relevant. In
particular, timings should be obtained under comparable hardware,
software, initialization, and stopping conditions. Without such
controls, apparent differences in computational performance may
reflect implementation choices rather than intrinsic differences
between the two methodological frameworks.

\subsection{Structural Properties of the Coefficient Matrix}
The coefficient matrix plays a central role in determining whether
a particular projection-based or modulus-based method is suitable
for a given LCP. Matrix properties influence the existence and
uniqueness of solutions, the validity of convergence conditions,
and the effectiveness of the operations used within an iterative
scheme. Accordingly, the structural characteristics of the matrix
should be examined before drawing conclusions about the relative
merits of competing algorithms.

For projection-based methods, favorable structural properties may
facilitate the establishment of suitable fixed-point estimates
and the selection of admissible step sizes. Depending on the
formulation, bounds associated with matrix norms, monotonicity,
or spectral properties can be useful in deriving convergence
conditions. The practical advantages of these methods may be
particularly pronounced when matrix-vector products dominate
the computational workload and the projection operation is
inexpensive. Nevertheless, a simple update formula does not by
itself ensure rapid convergence, since the conditioning and
structural characteristics of the underlying problem can still
strongly influence the behavior of the iteration.

For modulus-based methods, the effectiveness of a scheme is
closely related to the suitability of the selected matrix
splitting. Different decompositions may yield different iteration
matrices and, hence, different convergence characteristics even
when applied to the same LCP. The availability of a computationally
convenient splitting can be particularly important for sparse
problems, where preserving sparsity and avoiding expensive
factorizations may substantially reduce the overall cost.
However, a splitting that is inexpensive to evaluate is not
necessarily the one that produces the fastest convergence.
The design of an effective modulus-based algorithm must therefore
balance algebraic convenience, convergence requirements, and
computational expense.

The distinction between matrix classes is equally important.
For example, theoretical guarantees established for a particular
class of $P$-matrices or $H_{+}$-matrices should not automatically
be extended to general nonsingular matrices or to arbitrary
complementarity problems. Likewise, favorable behavior observed
on a collection of structured test problems does not establish
uniform superiority over other methods. The matrix assumptions,
problem dimensions, sparsity patterns, and numerical properties
of the test instances must be taken into account when interpreting
published results.

Consequently, matrix structure should be treated as a primary
selection criterion rather than as a secondary descriptive feature.
A method that exploits the available structure effectively may
outperform an alternative with a more general formulation, whereas
a theoretically attractive method may become less competitive
when its structural assumptions are not satisfied or its required
linear algebra operations are costly.

\subsection{Parameter Sensitivity, Robustness, and Implementation}
Algorithmic parameters provide an important means of controlling
the convergence behavior of both projection-based and modulus-based
methods. In projection-based schemes, the step size or relaxation
parameter may influence the stability of the iteration and the
rate at which the computed sequence approaches a solution.
Parameters that are too conservative can lead to slow progress,
whereas values outside the admissible range may invalidate the
convergence guarantee or cause numerical instability. The
practical performance of these methods therefore depends not
only on the existence of a suitable parameter range but also
on the ease with which effective values can be identified.

In modulus-based methods, parameter dependence can arise through
the transformation, relaxation, and matrix-splitting components
of the algorithm. Their interaction may produce a larger design
space in which different parameter combinations lead to
substantially different numerical behavior. Although this
flexibility can be exploited to improve convergence, it may also
increase the difficulty of selecting parameters that perform
consistently across different problem instances. In particular,
a parameter configuration that performs well for one matrix
may not remain effective when the matrix dimension, sparsity
pattern, conditioning, or complementarity structure changes.

Robustness should therefore be distinguished from peak performance.
A method may achieve an excellent iteration count on a selected
test problem but exhibit substantial deterioration under modest
changes in problem data or parameter values. A more robust method
may deliver slightly less impressive performance on the easiest
instances while maintaining predictable behavior over a broader
range of problems. For practical applications, the latter
characteristic can be important when the solver must be used
repeatedly on related systems or on problems whose structural
properties are not known precisely in advance.

Implementation considerations further influence this comparison.
Projection-based schemes may be easier to implement when the
projection admits a simple componentwise representation.
Modulus-based schemes may require additional routines for matrix
splitting, auxiliary linear solves, or parameter-dependent
preconditioning. Nevertheless, neither framework is inherently
restricted to simple or sophisticated implementations; the
actual complexity depends on the particular algorithm and the
computational environment. Automated parameter selection,
adaptive relaxation, efficient sparse linear algebra, and
reusable factorizations can help reduce these practical
difficulties in either framework.

A fair assessment should consequently examine parameter sensitivity,
the cost of tuning, numerical stability, and performance across
multiple problem instances. Reporting results for only one
carefully selected parameter configuration may conceal important
limitations. Where feasible, sensitivity studies and reproducible
parameter-selection procedures provide stronger evidence of
robustness than isolated numerical comparisons.

\subsection{Advantages, Limitations, and Practical Trade-Offs}
The principal attraction of projection-based methods lies in
their direct treatment of complementarity through projection
operations and their potential for computationally economical
iterations. When the projection is inexpensive and the update
can be expressed using efficient matrix-vector operations, these
methods provide a natural framework for constructing simple
iterative solvers. Their formulations can also be combined with
relaxation and matrix-splitting ideas, allowing algorithmic
extensions without abandoning the underlying projection
mechanism. These features make projection-based approaches
particularly appealing when implementation simplicity, low
iteration overhead, and efficient sparse computations are
important objectives.

Their limitations arise primarily from the convergence conditions
associated with the specific iteration and the influence of
algorithmic parameters. A low-cost projection step does not
guarantee a low iteration count, and convergence may be difficult
to establish when the required assumptions on the coefficient
matrix or iterative mapping are not satisfied. Furthermore,
performance can deteriorate when the selected parameters are
poorly matched to the problem structure. Projection-based
methods should therefore be evaluated in terms of both their
computational simplicity and the conditions under which their
convergence can be established.

Modulus-based methods offer a flexible algebraic framework for
transforming complementarity conditions into nonlinear systems
that can be treated through a variety of iterative procedures.
This framework has supported the development of matrix-splitting,
double-splitting, two-sweep, preconditioned, and inexact methods
\cite{Bai2010,WuLi2016,DaiLiBaiQiu2019}. Such developments
provide opportunities to exploit matrix structure, modify
relaxation strategies, and construct algorithms tailored to
particular classes of LCPs. When an effective splitting is
available, this flexibility may translate into improved
convergence or better computational performance.

At the same time, the additional algebraic flexibility can
introduce implementation overhead and more involved parameter
selection. Auxiliary linear systems may dominate the computational
cost, while the nonsmooth nature of the modulus mapping requires
careful treatment in the relevant convergence analysis. Moreover,
the performance of a particular modulus-based scheme depends
on the splitting and parameter configuration; the advantages
of the general framework do not automatically transfer to every
algorithm constructed within it.

The central trade-off is therefore not simply between a
computationally simple method and a more elaborate one. Rather,
it concerns the balance between the cost of an iteration, the
number of iterations required, the strength of the convergence
guarantees, the sensitivity to algorithmic parameters, and the
extent to which the method exploits the problem structure.
Projection-based methods may be preferable when a simple and
efficient projection yields satisfactory convergence, whereas
modulus-based methods may be attractive when a suitable
transformation and splitting provide a meaningful computational
advantage. These observations should be interpreted as
method-dependent tendencies rather than universal rules.

\subsection{A Framework for Numerical Comparison}
A meaningful comparison of projection-based and modulus-based
methods requires a consistent experimental protocol. Numerical
results obtained from different test problems, stopping
criteria, parameter-selection procedures, or computational
environments cannot be compared reliably without accounting
for these differences. A structured evaluation should therefore
combine accuracy, efficiency, robustness, and implementation
considerations within a common experimental framework.

The following measures provide a useful basis for such an
assessment:

\begin{enumerate}
	\item \textbf{Iteration count:} the number of iterations
	required to satisfy a common stopping criterion, providing
	an indication of the convergence speed.
	
	\item \textbf{CPU time:} the total execution time, including
	relevant preprocessing and initialization costs when these
	are part of the intended computational workflow.
	\item \textbf{Residual accuracy:} the quality of the
	computed solution, measured using a consistently defined
	complementarity residual or an equivalent error measure.
	\item \textbf{Memory requirements:} the storage needed for
	the coefficient matrix, splitting components, auxiliary
	variables, and intermediate computations.
	\item \textbf{Robustness:} the ability to maintain reliable
	convergence across different matrix structures, problem
	dimensions, and admissible parameter choices.
	\item \textbf{Implementation cost:} the additional effort
	associated with parameter tuning, matrix factorization,
	auxiliary linear systems, and method-specific preprocessing.
\end{enumerate}
To compare the numerical performance of different methods fairly, the
same natural residual is used to measure complementarity accuracy.
For an LCP with matrix $A$ and vector $q$, the iteration is terminated
when the following stopping criterion is satisfied:
\[
\left\|\min\left\{z^{(k)},Az^{(k)}+q\right\}\right\|_2
\leq \varepsilon,
\]
where the minimum is taken componentwise and $\varepsilon>0$
is a prescribed tolerance.

In addition to reporting aggregate performance, numerical
experiments should include problems with different dimensions,
sparsity patterns, and matrix structures. Parameter-selection
procedures should be documented, and the computational
environment should be described sufficiently to support
reproducibility. When several instances are tested, reporting
success rates and the distribution of execution times can
provide more informative evidence than reporting only the
best-performing case. Such practices help distinguish
systematic advantages from performance differences that arise
from isolated examples or implementation-specific choices.
\subsection{Guidelines for Method Selection}
The selection of a suitable iterative method should begin
with an examination of the mathematical and computational
characteristics of the target LCP. The first consideration
is the available information about the coefficient matrix,
including its matrix class, sparsity, conditioning, and
any structural properties that can be exploited by the
algorithm. These characteristics determine which convergence
results are applicable and which computational operations
can be performed efficiently. A method should not be selected
solely because it performs well on a different class of
problems or because its formulation appears algebraically
simpler.

Projection-based methods are natural candidates when the
projection operation is inexpensive, the iterative update
can exploit efficient matrix-vector products, and suitable
convergence conditions can be established. They may be
particularly attractive in implementations where minimizing
per-iteration overhead is important. Their suitability
should nevertheless be verified through an assessment of
parameter admissibility, convergence reliability, and the
accuracy attainable within the available computational
budget.

Modulus-based methods merit consideration when the problem
structure supports an effective matrix splitting and the
resulting auxiliary systems can be solved economically.
Their flexibility may be useful when different splittings
or relaxation strategies can be designed to exploit specific
matrix properties. However, the potential benefits should
be weighed against the cost of constructing and solving the
split systems, as well as the sensitivity of the resulting
iteration to parameter choices.

For large-scale sparse LCPs, particular emphasis should be
placed on memory usage, preservation of sparsity, parallel
computability, and the efficiency of repeated linear algebra
operations. For applications involving many related LCPs,
the cost of parameter tuning and the possibility of reusing
matrix factorizations or preconditioners may also influence
the choice. In settings where reliability is more important
than the best performance on a single instance, robustness
across a representative problem set should receive greater
weight than isolated reductions in iteration count.

Ultimately, the most appropriate approach is the one that
satisfies the required accuracy and reliability constraints
at an acceptable total computational cost. This decision
should be supported by the convergence theory relevant to
the problem class and by numerical experiments conducted
under comparable conditions. Such a selection strategy
avoids treating either reformulation as universally preferable
and instead emphasizes the compatibility between the
algorithm and the problem it is intended to solve.
\subsection{Open Challenges and Future Research Directions}
Despite substantial progress in both projection-based and
modulus-based methods, several issues remain relevant to the
development of reliable and computationally efficient LCP
solvers. One important challenge is the construction of
algorithms whose convergence guarantees remain meaningful
for broad problem classes without imposing unnecessarily
restrictive assumptions on the coefficient matrix. Stronger
results concerning admissible parameter ranges, convergence
rates, and sensitivity to perturbations would help clarify
the circumstances under which each framework is most effective.
A second direction concerns adaptive parameter selection.
Fixed parameter values may perform well for one problem
instance but become inefficient when the matrix structure
or numerical conditioning changes. Methods that adjust
relaxation, transformation, or splitting parameters using
computationally inexpensive information could improve
practical robustness. Such strategies should, however,
account for the overhead of adaptation and should be
evaluated against well-defined baseline methods to establish
whether the additional complexity yields a net benefit.
A further challenge is the development of scalable methods
for very large and sparse problems. Efficient sparse
factorizations, matrix-free operations, preconditioning,
parallel implementations, and inexact inner solves may
substantially influence the practical performance of both
frameworks. In this context, the key question is not merely
whether a theoretical iteration converges, but whether it
can be implemented in a way that preserves the relevant
structure and controls the total computational effort.
Comparative studies that incorporate these implementation
aspects would provide a more realistic picture of algorithmic
performance.

Finally, the literature would benefit from more systematic
and reproducible numerical benchmarks. Comparisons based
on common test sets, consistent residual definitions,
transparent parameter-selection procedures, and clearly
specified stopping criteria would make it easier to identify
the strengths and limitations of individual algorithms.
Rather than seeking a universally superior reformulation,
future research should focus on identifying the structural
conditions under which particular methods are effective
and on developing adaptive or hybrid strategies that
combine complementary computational advantages.

\subsection{Concluding Comparative Perspective}

Projection-based and modulus-based methods represent
complementary approaches to the numerical solution of
LCPs. The former provide a
direct and often computationally economical mechanism
for enforcing complementarity through projection, whereas
the latter offer an algebraically flexible framework for
constructing iterative schemes through transformations
and matrix splittings. These distinctions influence the
convergence analysis, parameter dependence, computational
workload, and practical implementation of the resulting
algorithms.

A comparative study by\textbf{ Hadjidimos and Zhang (2018)}
~\cite{HadjidimosZhang2018} provides a useful example of
these differences. The authors compared three classes of
iterative methods through their respective best methods:
the projected matrix analogue Gauss-Seidel (PMAGS)
method, the stationary extrapolated block modulus algorithm
(SEBMA), and the modulus-based matrix splitting
Gauss-Seidel (MGS) method.

Consider the original LCP
(\ref{eq:LCP}),
where \(A\) is an \(H_+\)-matrix. The splitting is 
\[
A=D-L-U,
\]
where \(D=\operatorname{diag}(A)\), let
\(E=D^{-1}\). The authors first scale the LCP using \(E\)
and then choose a positive diagonal matrix \(Y\) such that
\((EA)Y\) is strictly diagonally dominant. The equivalent
transformed problem is defined by
\[
\bar{A}=Y^{-1}EAY
=I-\bar{L}-\bar{U},
\qquad
\bar{q}=Y^{-1}Eq,
\]
where
$
\bar{z}=Y^{-1}z,~
\bar{w}=Y^{-1}Ew.
$
Consequently, the transformed LCP is
\[
\bar{w}
=\bar{A}\bar{z}+\bar{q},
\qquad
\bar{w}\geq0,\quad
\bar{z}\geq0,\quad
\bar{w}^{T}\bar{z}=0.
\]

For this transformed problem, the PMAGS method is given by
\[
\bar{z}^{(k+1)}
=
\left[
\bar{L}\bar{z}^{(k+1)}
+\bar{U}\bar{z}^{(k)}
-\bar{q}
\right]_{+},
\]
where \([\cdot]_+\) denotes projection onto the
nonnegative orthant. The SEBMA method, representing the
block modulus algorithm class, is defined by
\[
(I+\bar{A})\bar{x}^{(k+1)}
=(I-\bar{A})|\bar{x}^{(k)}|
-\bar{q}.
\]
In contrast, the MGS method, belonging to the
modulus-based matrix splitting class, is given by
\[
\begin{aligned}
	(2I-\bar{L})\bar{x}^{(k+1)}
	={}&\bar{U}\bar{x}^{(k)}
	+(\bar{L}+\bar{U})
	|\bar{x}^{(k)}|-\bar{q}.
\end{aligned}
\]
These iterations illustrate three distinct approaches
to solving the transformed LCP.

The authors identified PMAGS as the best method within
the projected matrix analogue AOR class and MGS as the
best method within the modulus-based matrix splitting
AOR class. Their theoretical comparison uses the
associated nonnegative iteration matrices. For PMAGS,
the matrix is
\[
T_{\mathrm{PMAGS}}
=(I-|\bar{L}|)^{-1}|\bar{U}|,
\]
whereas for MGS it is
\[
T_{\mathrm{MGS}}
=(2I-|\bar{L}|)^{-1}
\bigl(|\bar{L}|+2|\bar{U}|\bigr).
\]
The authors established that
\[
\rho(T_{\mathrm{PMAGS}})
\leq \rho(T_{\mathrm{MGS}})<1,
\]
with strict inequality between the two spectral radii
when the original matrix is irreducible. Thus, PMAGS has
the better spectral-radius bound in this theoretical
comparison. However, the numerical experiments showed
that MGS required less CPU time than PMAGS, because
the projected Gauss--Seidel iteration updates components
successively, whereas the matrix inverse used in MGS can
be computed once before the iterations begin. The
experiments also showed that SEBMA could break down for
some tested problems whose original coefficient matrices
were not strictly diagonally dominant, whereas it performed
best among the three methods for reasonably strictly
diagonally dominant matrices. These findings illustrate
that the spectral-radius comparison and the practical
computational performance of iterative methods need not
produce the same ranking.

Nevertheless, the relative performance of the two
frameworks cannot be established independently of the
specific problem and iterative scheme. Convergence
guarantees, matrix structure, per-iteration cost, total
execution time, residual accuracy, and robustness must
be considered together. In particular, a low-cost
iteration is valuable only when the method reaches
the prescribed accuracy efficiently, while a more
elaborate scheme is justified only when its additional
computational requirements lead to a meaningful
improvement in overall performance.

The principal conclusion is therefore that the choice
between projection-based and modulus-based methods
should be guided by mathematical suitability and
empirical evidence rather than by the reformulation
framework alone. A rigorous combination of convergence
analysis, structural assessment, and reproducible
numerical evaluation provides the most reliable basis
for selecting an appropriate LCP solver. This perspective
also motivates continued research into adaptive,
structure-aware, and computationally scalable algorithms
capable of delivering reliable performance across
diverse classes of complementarity problems.

\section{Conclusion}

This review has presented a comprehensive examination of the principal equivalent reformulations of the LCP and the iterative methods developed from these representations. The reformulations have been organized into two broad frameworks: projection-based and modulus-based. The projection-based framework covers minimum-function, maximum-function, projection, and positive-negative part formulations, whereas the modulus-based framework includes scaled modulus and absolute-value formulations. The mathematical foundations of these approaches have been discussed to establish their relationships with the original complementarity conditions and their role in the development of iterative solution methods.

Furthermore, the associated iterative algorithms have been reviewed with respect to their mathematical construction, matrix-splitting strategies, convergence conditions, and computational characteristics. The comparative discussion has highlighted the significance of coefficient matrix properties and algorithmic parameters in determining the convergence behavior and computational effectiveness of different methods. It has also emphasized the importance of considering iteration requirements, computational cost, numerical accuracy, and robustness when assessing the suitability of an iterative scheme. The overall analysis demonstrates that both projection-based and modulus-based approaches offer valuable strategies for solving LCPs, with their applicability depending on the mathematical and computational characteristics of the problem.
The review thus brings together the theoretical developments and algorithmic approaches within a unified classification, providing a systematic understanding of the available reformulation-based methods for solving LCPs. By examining their mathematical properties and comparative characteristics, it offers a useful perspective for understanding the development and applicability of these iterative techniques.

Future research may focus on developing adaptive parameter-selection procedures, more effective matrix-splitting strategies, and hybrid algorithms that integrate the complementary strengths of projection-based and modulus-based methods. Further investigations into convergence analysis, advanced preconditioning techniques, and scalable algorithms for large-scale sparse LCPs may lead to improved computational performance. The exploration of parallel implementations and systematic numerical benchmarking under standardized computational settings may also contribute to the development and practical application of efficient, reliable, and robust iterative methods for linear complementarity problems.

\section*{Acknowledgments} Abhishek Kumar Singh Sengar acknowledges the support of the Ministry of Education, Government of India, through the GATE fellowship (Registration No. MA22S25023148).
\bibliographystyle{plain}
\bibliography{Review3}
\end{document}